\documentclass[11pt,leqno]{article}
\usepackage{amssymb, amsmath, amsthm, amsxtra}

\theoremstyle{plain}
\newtheorem{theorem}{Theorem}
\newtheorem{lemma}[theorem]{Lemma}

\theoremstyle{remark}
\newtheorem{remark}{Remark}

\theoremstyle{remark}
\newtheorem{conject}[theorem]{Conjecture}

\newcommand{\lf}{\left}
\newcommand{\rt}{\right}

\newcommand{\be}{\begin{equation}}
\newcommand{\ee}{\end{equation}}

\newcommand{\al}{\alpha}
\newcommand{\ep}{\epsilon}
\newcommand{\eit}{e^{i\tht}}
\newcommand{\hde}{\hat{\de}}
\newcommand{\hka}{\hat{\chi}}

\newcommand{\wht}{\widehat{\tht}}
\newcommand{\wh}{\widehat}
\newcommand{\wbe}{\widehat{\beta}}
\newcommand{\wtbe}{\tilde{\beta}}
\newcommand{\wpi}{\widetilde{\pi}}
\newcommand{\whpi}{\widehat{\pi}}
\newcommand{\tvph}{\widetilde{\varphi}}
\newcommand{\tf}{\widetilde{f}}
\newcommand{\de}{\delta}
\newcommand{\De}{\Delta}
\newcommand{\lb}{\lambda}
\newcommand{\La}{\Lambda}

\newcommand{\tht}{\theta}

\newcommand{\sg}{\sigma}
\newcommand{\Sg}{\Sigma}
\newcommand{\p}{\partial}
\newcommand{\oc}{\ovl{c}}
\newcommand{\ovl}{\overline}
\newcommand{\CC}{\mathbb{C}}
\newcommand{\DD}{\mathbb{D}}
\newcommand{\EE}{\mathbb{E}}
\newcommand{\NN}{\mathbb{N}}
\newcommand{\RR}{\mathbb{R}}
\newcommand{\ZZ}{\mathbb{Z}}

\newcommand{\CM}{\mathcal{M}}
\newcommand{\MCH}{\mathcal{H}}
\newcommand{\MCU}{\mathcal{U}}
\newcommand{\tr}{\,{\mathrm{tr}}\,}

\newcommand{\Bl}{|||\beta|||}

\newcommand{\tit}{\textit}
\newcommand{\tbf}{\textbf}

\newcommand{\const}{\text{const.}}
\newcommand{\tons}{\text{ons}}
\newcommand{\tdg}{\text{diag}}
\newcommand{\tph}{\text{phys}}
\newcommand{\half}{\frac{1}{2}}

\newcommand{\Cara}{{Carath\'eodory} }

\numberwithin{remark}{section}

\begin{document}
\title{Toeplitz matrices and Toeplitz determinants under the
impetus of the Ising model.  Some history and some recent results} 
\author{Percy Deift\\
\tit{Courant Institute, New York} \\
\and
Alexander Its\\
\tit{IUPUI, Indianapolis}\\
\and Igor Krasovsky\\
\tit{Imperial College, London}}
\date{}
\maketitle

\bigskip
\bigskip
\bigskip
\centerline{\tit{In memory of Bruria Kaufman, 1918--2010}    }

\newpage
\tableofcontents

\bigskip\bigskip\bigskip
\noindent
Acronyms: RHP (p. 4), SSLT (p. 5), OPUC (p. 14), FH (p. 29), CMV (p. 56)

\newpage

\section{Setting the problem}
An $n\times n$ \textit{Toeplitz matrix} $T_n$ is a matrix with coefficients of the form $(T_n)_{jk} = c_{j-k},
\;0\le j,\; k \le n-1$, for some given sequence $\{c_\ell\}_{\ell \in \ZZ}$.  An $n\times n$ \tit{Toeplitz determinant} $D_n(T)$
is the determinant of some given Toeplitz matrix $T=T_n$.
If $\varphi = \varphi\!\lf(\eit\rt)$ is an integrable
function on the unit circle $S^1=\{z=\eit\}$ (oriented in the positive direction) with Fourier coefficients
$$
\varphi_{\ell} = \int^{\pi}_{-\pi} e^{-i\ell \tht} \, \varphi \!\lf(\eit\rt) \,\frac{d\tht}{2\pi},
\qquad \ell \in \ZZ
$$
the Toeplitz matrices $T_n(\varphi)$ and Toeplitz determinants $D_n(\varphi)$ associated with $\varphi$ are given by
\be
T_n (\varphi) = \lf\{ \varphi_{j-k} \rt\}_{0 \le j,\; k \le n-1}
\ee
and
\be
D_n (\varphi) = \det \, T_n (\varphi)
\ee
respectively. Note that if $c_{-\ell}= \oc_{\ell}$, then $T_n$ is self-adjoint.  This is true, in particular, if
$\varphi$ is real valued, and it follows that for such functions $\varphi$, $D_n(\varphi)$ and the eigenvalues
$\lb^{(n)}_1, \dots \lb^{(n)}_n$ of $T_n(\varphi)$ are real.

Toeplitz matrices and determinants are named for Otto Toeplitz, who, in his Habilitationsschrift in 1907
(see \cite{Toep1} \cite{Toep2}), initiated the study of quadratic forms $\Sg \,\varphi_{jk}\, x_j y_k$ with coefficients
of special type $\varphi_{jk}= \varphi_{j-k}$.  At the time Toeplitz was a young researcher in
G\"{o}ttingen where Hilbert was developing his general and abstract theory of functional analysis.
Toeplitz introduced the forms $\Sg\, \varphi_{j-k}\, x_j y_k$, together with the associated matrices
$T_n(\varphi)= \{ \varphi_{j-k}\}_{0\le j, \; k\le n-1}$ and determinants $D_n(\varphi)= \det
(T_n (\varphi))$, in order to give concrete examples of Hilbert's general theory which could be analyzed in
great detail.  Toeplitz called the forms $\Sg\, \varphi_{j-k}\, x_j\, y_k$ ``$L$-forms'' because of their
connection, via Fourier theory, with Laurent series $\sum^\infty_{-\infty} \varphi_k \,z^k$.

A striking example of the connection between Laurent series, and more broadly analytic function theory,
on the one hand, and $L$-forms on the other, is given by  the following problem studied by Carath\'eodory:
For which sequences $c_0, c_1, \dots, c_n$ of complex numbers does the polynomial
$$
p(z) = c_0 + c_1 z+ \dots + c_n z^n
$$
have an extension to an analytic function $\varphi(z)$ in the unit disk $\DD=\{ |z| < 1\}$
\be
\begin{aligned}
\varphi(z) = &d_0 + d_1 z + \dots + d_n\, z^n + d_{n+1} \, z^{n+1} + \dots \quad , \\
& \quad d_j = c_j, \qquad 0\le j \le n
\end{aligned}
\ee
such that $\Re  \varphi(z) \ge 0$ in $\DD$? In \cite{Car1} in 1907, \Cara characterized such
sequences in terms of Minkowski's theory of convex bodies, but a few years later in 1911 Toeplitz
\cite{Toep3} gave the following algebraic characterization (see also \cite{Car2}, \cite{CarFej}): $p(z) =
c_0 + c_1\, z + \dots + c_n \, z^n $ has an extension $\varphi(z)$ as above with $\Re  \varphi(z) \ge 0$
in $\DD$ if and only if the Toeplitz matrix
\be
\left(
\begin{array}{lllcccll}
2c_0        &c_1     &c_2     &&\dots&     &c_{n-1}  &c_n \\
c_{-1}      &2c_0    &c_1     &&\dots&     &&c_{n-1}  \\
c_{-2}      &c_{-1}  & 2c_0   &&&          &&\\
\phantom{c_{-n}}  \vdots&      &        &&&          & &\\
c_{-n+1} &  \dots           &        &&&          & \\
c_{-n}       &c_{-n+1}      &\dots &&&   & &2c_0
\end{array}
\rt)
\ee
with $c_{-k} \equiv \oc_k$, $0\le k\le n$, is non-negative definite.  Over the years, up till
the present, there have been many applications of this beautiful result, both theoretical
and applied (see e.g.\ \cite{Weg}).  Many of the results and applications
in this paper, however, rely on the connection
between Toeplitz matrices/determinants and another problem in analytic function theory, the
Riemann-Hilbert Problem (RHP) --- see, in particular, \cite{CG} \cite{LS} and the references therein. We note
that Wiener-Hopf factorization theory and the Wiener-Hopf method fall under the general rubric of RHP.  

An extraordinary variety
of problems in mathematics, physics and engineering, can be expressed in terms of Toeplitz matrices
and determinants.
This is a remarkable fact that could hardly have been anticipated in 1907, and there is a vast literature on the subject.
A classical reference is the text
by Grenander and Szeg\H{o}\footnote{For many years ``Szeg\H{o}'' was mis-spelled in the literature
as ``Szeg\" o''. Knuth first introduced the Hungarian umlaut into TeX specifically so that names like 
Szeg\H{o}, Erd\H{o}s and K\H{o}nig would be spelled correctly. (Private communication to the authors from 
D. Knuth.)}   
 \cite{GreSz} in which the authors describe applications of Toeplitz
matrices to problems in analytic function theory, stationary stochastic processes and linear
estimation theory in statistics.  In addition to problems in probability theory and statistics,
Fisher and Hartwig \cite{FisHart1} present a long and varied  list of applications of Toeplitz determinants to problems
in statistical mechanics, including the ground-breaking work of Kaufman and Onsager on the Ising model.
Many results in the theory of Toeplitz determinants can be found, for example, in the seminal papers of
Harold Widom (see \cite{Wid1} \cite{Wid2} \cite{Wid3} \cite{Wid4}, amongst many others) and in the books of
B\"{o}ttcher, Silbermann, and Grudsky  (see \cite{BottSilb1} \cite{BottSilb2} \cite{BottSilb3} \cite{BottGr},
amongst others). For a review of some more recent developments in the theory of Toeplitz determinants,
see \cite{Kr}.

Closely related to Toeplitz matrices are Toeplitz operators, defined as follows (see e.g.\ \cite{Doug},
\cite{BottSilb3}).  Let $S^1$ denote the unit circle as above with Lebesgue measure $\frac{d\tht}{2\pi}$,
and let $L^2(S^1)$ denote the associated Hilbert space.  The Hardy space $H^2$ is defined as the closed
subspace
$$
\lf\{ f\in L^2(S^1) : \frac{1}{2\pi} \int^{2\pi}_0 f\!\lf(\eit\rt)  e^{in\tht} \, d\tht =0,
\qquad n =1,2,\dots \rt\}.
$$
The orthogonal projection $P_+$ of $L^2(S^1)$ onto $H^2$ is given by
$$
L^2(S^1) \ni f = \sum^\infty_{-\infty} f_k \, e^{ik\tht} \longmapsto  P_+ f = \sum^\infty_0 f_k \,
e^{ik\tht} \in H^2.
$$
Now let $\varphi \in L^\infty(S^1)$.  Then the \tit{Toeplitz operator} $T(\varphi): H^2 \to H^2$
associated with $\varphi$ is given by
\be\label{eq1}
T(\varphi) f = P_+ \, \varphi f, \qquad f \in H^2 .
\ee
The operator $T(\varphi)$ is bounded in $H^2$, and with respect to the standard basis $\lf\{e^{ik\tht}:
k\ge 0\rt\}$, $T(\varphi)$ is represented by the semi-infinite matrix
$\lf\{ \varphi_{j-k}\rt\}_{j,\,k \ge 0}\:$
acting in $\ell^2_+ = \{u= (u_0, u_1, \dots ): \sum^\infty_0| u_j|^2 < \infty\}$. Thus the
Toeplitz matrices $T_n(\varphi)$ are \tit{finite sections} of the Toeplitz operator $T(\varphi)$.  
The function $\varphi$ is frequently referred to as the {\it symbol} of the operator  $T(\varphi)$ and the sequences
$\{T_n(\varphi)\}_{n\ge 1}$, $\{D_n(\varphi)\}_{n\ge 1}$.

Matrices with coefficients of the form $\lf\{c_{j+k} \rt\}_{0\le j, \, k\le n-1}$ are called \tit{Hankel matrices}.
These matrices are in turn finite sections of the semi-infinite \tit{Hankel operators} with matrix
coefficients $\lf\{c_{j+k}\rt\}_{j,\, k\ge 0}$  acting in $\ell^2_+$.  The theory of Hankel matrices
and operators is parallel to, and closely related to the theory of Toeplitz matrices and operators
(see e.g. \eqref{eq108} et seq.\ below).  Operators acting in $L^2(0, \infty)$ with kernels of the 
form $K(x-y)$, or
of the form $K(x+y),\, x, y \ge 0$, may be viewed as continuum limits of Toeplitz matrices, or Hankel matrices,
and are called Wiener-Hopf and Hankel operators respectively (see e.g.\ \cite{BottSilb3} \cite{Pel}).

Our goal in this paper is not to review the above applications and developments in a systematic way,
but rather to focus on \tbf{one central problem} in the theory of Toeplitz determinants, viz., the
Szeg\H{o} Strong Limit Theorem (SSLT) (see Theorem 3 below).  We will, in particular, trace the history
of this theorem and its later generalizations, in response to a succession of questions raised in the
analysis of the Ising model.
An earlier, and very informative, description of the influence of the Ising model on SSLT and related
developments, was given by B\"{o}ttcher in \cite{Bott1}.

The story begins in 1915 when Szeg\H{o} \cite{Sz1}, at the very start of his career, proved the following
result conjectured earlier by Polya.

\begin{theorem}\label{theorem1}
Let $\varphi\!\lf(\eit\rt)>0$ be a continuous, positive function on the unit circle $S^1$ and let $D_n(\varphi)$
be the associated Toeplitz determinant.  Then
\be\label{eq2}
\lim_{n\to \infty}\, \frac{1}{n} \,\log D_n(\varphi) = \frac{1}{2\pi} \int^\pi_{-\pi} \log \varphi
\!\lf(\eit\rt) d\tht.
\ee
\end{theorem}
Note that $D_n(\varphi)>0$.  Indeed for $u= (u_0, u_1, \dots, u_{n-1})^T \in \CC^{n}$, a simple calculation shows that
$$
\lf(u, T_n(\varphi) \, u\rt) = \int^\pi_{-\pi} \varphi\!\lf(\eit\rt)
\lf| \sum^{n-1}_{k=0} u_k \: e^{ik\tht}\rt|^2 \frac{d\tht}{2\pi}
$$
and as $\varphi\!\lf(\eit\rt)>0$, we conclude that $T_n(\varphi)$ is strictly positive definite,
$T_n(\varphi) >0$.  Thus $D_n(\varphi) = \det \, T_n (\varphi)>0$.
In addition, the eigenvalues $\lb^{(n)}_1, \dots, \lb^{(n)}_n$ of $T_n(\varphi)$ are positive and so  \eqref{eq2}
can be rewritten in the form
\be\label{eq3}
\lim_{n \to \infty} \, \frac{\log \lb^{(n)}_1 + \dots + \log\lb^{(n)}_n}{n} =
\frac{1}{2\pi} \int^\pi_{-\pi} \log \varphi\!\lf(\eit\rt)d\tht.
\ee
Szeg\H{o} recognized that \eqref{eq3} is a special case of equidistribution for the $\lb^{(n)}_j$  in the sense
of Weyl, and in 1920 \cite{Sz2} Szeg\H{o} proved the following result (see also \cite{GreSz}).
\begin{theorem}\label{theorem2}
Let $\varphi\!\lf(\eit\rt)$ be a real-valued function in $L^\infty (S^1)$ with $m \le \varphi\!\lf(\eit\rt)
\le M$
a.e.  If $F(\lb)$ is any continuous function defined on the interval $m\le \lb \le M$, we have
\be\label{eq4}
\lim_{n \to \infty} \frac{F\lf(\lb^{(n)}_1\rt) + \dots + F\lf(\lb^{(n)}_n \rt)}{n} = \frac{1}{2\pi}
\int^\pi_{-\pi} F\lf( \varphi\!\lf(\eit\rt)\rt) d\tht.
\ee
\end{theorem}
Taking $F(\lb) = \log \lb$, we clearly recover \eqref{eq3}.  We consider more detailed asymptotic properties
of the eigenvalues below (see Section \ref{sec9}).

Relation \eqref{eq2} may be rewritten in the form
\be\label{eq5}
D_n (\varphi) = \exp\lf[ n \int^\pi_{-\pi} \log\varphi\!\lf(\eit\rt) \frac{d\tht}{2\pi} + o(n)\rt]
\ee
as $n\to \infty$.  The issue of determining the precise nature of the term $o(n)$ came to the fore
in the following way.  

The two-dimensional Ising model, the central model in statistical mechanics 
(see \cite{McWu1} \cite{Pal} \cite{Lie} \cite{Dom1}), named for E.~Ising
(see \cite{Is}), concerns the interaction of random spins $\sg_{i,j} = \pm 1$ at sites $(i,j)$ in $\ZZ^2$.
Of great interest is the situation in which only  nearest-neighbor spins interact and the interaction energy is
given by
$$
-J_1\, \sg_{i,j} \, \sg_{i,j+1} - J_2 \, \sg_{i,j} \, \sg_{i+1, j}
$$
where the vertical and horizontal interaction constants, $J_1$ and $J_2$ respectively, are translation invariant, and in addition, $J_1$ and $J_2$ are positive, so that the system is \tit{ferromagnetic}, i.e.\ parallel spins have lower
energy than anti-parallel spins.  As is standard in statistical mechanics, one analyzes the infinite system on
$\ZZ^2$ by first considering the spins in a finite rectangular box $\La \subset \ZZ^2$ of size $M\times N$
(see below for a discussion of boundary conditions)
and then letting $M, \, N \to \infty$.
For such a box the total interaction energy of the spins $\sg = (\sg_{i, j})$ is given by
\be\label{eq6}
E_\La (\sg) = -\sum _{\La} 
\lf[ J_1\, \sg_{i,j} \, \sg_{i,j+1} + J_2 \, \sg_{i,j} \, \sg_{i+1, j} \rt]
\ee
where $ J_1\, \sg_{i,j} \, \sg_{i,j+1}$ is included in the sum if $(i,j)$ or  $(i,j+1)\in \La$
and $ J_2 \, \sg_{i,j} \, \sg_{i+1, j}$ is included if $(i,j)$ or  $(i+1,j)\in \La$.
The associated normalized \tit{Gibbs measure} is then
\be\label{eq7}
Pr_\La(\sg) = \frac{1}{Z_\La} \, e^{-E_\La (\sg)/k_B \,T}
\ee
where $k_B$ is Boltzmann's constant and $T$ is the temperature.  Here the \tit{partition function} $Z_\La$ is
\be\label{eq8}
Z_\La = \sum_\sg e^{-E_\La (\sg)/k_B \,T}
\ee
where the sum is taken over all possible spin configurations $\sg$ in $\La$.  For fixed finite sets
$A\subset \La$,
one defines the \tit{correlation functions}
\be\label{eq9}
\lf\langle \prod_{(i,j) \in A} \,\sg_{i,j} \rt\rangle_\La = \sum_\sg \prod_{(i,j) \in A} \sg_{i,j}
\:Pr_\La (\sg).
\ee
The \tit{thermodynamic limits} of these correlation functions as $\La \uparrow \ZZ^2$,  i.e.
\mbox{$M, N\to \infty$},
\be\label{eq10}
\lf\langle \prod_{(i,j)\in A} \sg_{i,j}\rt\rangle \equiv \lim_{\La \uparrow \ZZ^2}
\lf\langle \prod_{(i,j) \in A} \sg_{i,j} \rt \rangle_{\La}
\ee
are the objects of principal physical interest.  In particular such correlations $\lf\langle \prod_{(i,j)\in A}
\sg_{i,j} \rt \rangle$ should model physical phenomena such as phase transitions.  For example, a bar magnet
has a critical temperature $T_c$ called the Curie point:  for temperatures $T < T_c$, the magnet exhibits
spontaneous magnetization, but for $T>T_c$ the magnetization is zero (in zero external field). This behavior cannot be
described by a function such as $\lf\langle\prod_{(i,j)\in A} \, \sg_{i,j}\rt\rangle_\La$
in \eqref{eq9}, which is real analytic for  $T>0$: in the thermodynamic limit, however,
analyticity can be, and indeed is,
destroyed.

In the presence of an external magnetic field, one must add a term of the form
$h \sum_{(i, j)\in \La} \sg_{i,j}$ to $E_\La (\sg)$,  $h\in \RR$,
but the most complete results have been obtained only when the field is absent, and we
will always assume $h=0$.  It remains an outstanding problem in the theory of the Ising model
to find the analogs
of the explicit formulae below for the correlation functions in the case $h\neq 0$.

We note the following:  In order to compute $E_\La (\sg)$ in \eqref{eq6}, we need to know
$\sg_{i,j}$ not only at the points $(i,j)\in \La$, but also at the points adjacent to $\La$.  In
particular, if $M$ and $N$ are integers and $\La=\lf\{(i,j) \in \ZZ^2: 0 \le i \le M, \:
0\le j \le N\rt\}$, then we need to know  $\sg_{M+1,\,j}$, $\sg_{-1,\, j}$, $0 \le j \le N$ and $\sg_{i,\, N+1}$,
$\sg_{i,\,-1}$, $0\le i \le M$.  In practice, this issue is addressed by imposing boundary conditions on $\La$.
One common choice is, for example, periodic boundary conditions, $\sg_{i+M+1,\, j} = \sg_{i,j}$,
$\sg_{i,\, j+N+1}= \sg_{i,\,j}$.  Another choice is the fully magnetized boundary conditions where $\sg_{i,j}
=+1$ (alternatively $\sg_{i,j}=-1$) at all points $(i,j)$ adjacent to $\La$.  General boundary
conditions are discussed, for example, in \cite{Gal}, and also from a different point of view,
in \cite{LebM-Lof}.  For each choice of boundary conditions, provided the thermodynamic limits
in \eqref{eq10} exist, we obtain a family of functions
\be\label{eq11}
\lf\langle \sg_{i,j}\rt\rangle, \quad \lf\langle \sg_{i,j} \; \sg_{i', j'} \rt\rangle,\quad
\dots
\ee
In the spirit of the classical moment problem, one says that the family in  \eqref{eq11}
determines a (Gibbs) equilibrium state (see e.g.\ \cite{Sim1} \cite{Ru}).  A priori, different boundary
conditions can give rise to different families \eqref{eq11}, and hence to different equilibrium
states.  For example (cf \cite{Pal}), for closed boundary conditions $\sg_{ij} = +1$ adjacent
to $\La$, one finds $\langle \sg_{0, 0} \rangle >0$, but for periodic boundary conditions,
by symmetry, $\langle \sg_{0, 0} \rangle =0$.  Nevertheless, it suffices for our purposes to note
that for the two-point correlation function that we discuss below, the thermodynamic limit is
independent of the boundary conditions (cf \cite{BGJ-LS}).  The same is true for the free
energy $F$ per unit spin in \eqref{eq13} below.

In 1 dimension, Ising \cite{Is} showed that the Ising model does not exhibit a phase transition
for any temperature $T>0$.  What about 2 or 3 dimensions?  In a landmark paper in 1936,
Peierls \cite{Pei} gave a simple argument asserting that in 2 or 3 dimensions, the Ising model
does indeed exhibit spontaneous magnetization at some temperature $T_c>0$.  It turns out that
Peierls' argument involved an incorrect step\footnote{First discovered by M. Fisher (see \cite{Grif}).},
which was corrected only many years later in 1964
by Griffiths \cite{Grif}. Griffiths confirmed Peierls' conclusion and ``Peierls argument'' remains
a standard tool in many  statistical situations.  The first exact quantitative result for
the 2-dimensional Ising model was obtained by Kramers and Wannier in 1941 \cite{KraWan} when they
wrote down the following formula for $T_c$.  In the case $J_1 =J_2 =J$,
\be\label{eq12}
\sinh \lf( \frac{2J}{k_B \, T_c} \rt) =1.
\ee
In 1942, Wannier gave a talk on this work at a meeting of the New York Academy of Sciences.
In an extraordinary denouement, at the end of the talk Onsager announced that he had obtained
an exact solution for the two-dimensional Ising model (without magnetic field).  
Onsager\footnote{All of Onsager's published papers, including his PhD dissertation, can be found 
in ``The Collected Works of Lars Onsager,
Eds. P. C. Hemmer, H. Holden and S. K. Ratkje, World Sci. Ser. in 20'th Century Phys., Vol. 17, World Scientific, 
Singapore, 1996''. An extensive collection of Onsager's hand-written notes, unpublished manuscripts, etc., 
can be found on the
Onsager archive at Trondheim, http://www.ntnu.no/ub/spesialsamlingene/tekark/tek5/arkiv5.php: much of this material was scanned from the ``Guide to the Lars Onsager Papers MS794'' at the Yale library.\label{foot3}}
published the proof of his result two years later in \cite{Ons1}.  More precisely, Onsager had obtained a
formula for the partition function $Z_\La$, and, in particular, in the case of a square lattice
$\La_N$ of size $N\times N$, with $J_1=J_2=J$, he showed that in the thermodynamic limit the free
energy per unit spin $F$ exists and is given by
\be\label{eq13}
\begin{aligned}
F&=-k_B \,T \lim_{N\to\infty} \frac{1}{N^2}\: \log Z_{\La_N} \\
&= -k_B \,T \lf[\log \lf(2 \cosh\frac{2J}{k_B T}\rt) + \frac{1}{2\pi^2} \int^\pi_0 d\phi_1 \int^\pi_0 d\phi_2
\log\lf(1-\frac{\kappa}{2}(\cos\phi_1+\cos\phi_2)\rt)\rt]\\
&= -k_B \,T \lf[\log \lf(2 \cosh \frac{2J}{k_B T}\rt) + \frac{1}{2\pi} \int^\pi_0
\log \frac{1}{2} \lf( 1 + \sqrt{
1-\kappa^2 \, \sin ^2 \phi}
\rt)d\phi\rt]
\end{aligned}
\ee
where $\kappa=2 \sinh \lf(2J/k_B T\rt)/\cosh^2 (2J/k_B T)$.
It follows from this formula that the specific heat $C=-T \,\frac{\p^2 F}{\p T^2}$
diverges logarithmically as $T\to T_c\pm 0$, where $T_c$ is precisely the critical temperature
obtained by Kramers and Wannier (note that for $T_c$ we have $\kappa=1$, and the argument 
of the logarithm in the double integral vanishes at the integration limit $\phi_1=\phi_2=0$).
Onsager's analysis allows for different interaction constants
$J_1 \neq J_2$, in which case equation \eqref{eq12} for $T_c$ must be replaced by
\be\label{eq14}
\sinh \lf( \frac{2J_1}{k_B\, T_c}\rt) \sinh \lf(\frac{2J_2}{k_B\, T_c}\rt) =1.
\ee

For the description of the events that followed, we refer to the informative, and entertaining,
recent papers by Baxter \cite{Bax} and McCoy \cite{McC1}.
In the mid-1940's, Kaufman, a doctoral student at Columbia University, began working with Onsager, who
was then at Yale, on the 2-D Ising model.  Onsager's computation of the partition function had
stunned the physics community by its ingenuity and complexity, and Kaufman was determined to find a
simpler approach.  As Onsager relates \cite{Ons2}, Kaufman `` \dots decided to look for a possible
connection with spinor theory.  Why not?  In fact, it seemed like very good sense, and so it was
\dots.   By the summer of 1946 she had a beautifully compact computation of the partition function,
bypassing all tedious details.''

The methods of Onsager and Kaufman for the partition function
are based on the calculation of the eigenvalues of the transfer matrix \cite{KraWan}
associated with the problem. Another simplified solution based on this approach was given in 1964 
by Schultz, Mattis, and Lieb \cite{ScML} .
On the other hand, it was observed in 1941 by van der Waerden 
\cite{vanderW} that the calculation of the Ising partition function is equivalent to the problem of 
counting closed polygons on the square
lattice. This observation was the origin of an alternative (and more direct) combinatorial approach to the Ising model.
A solution using such an approach was given by Kac and Ward in 1952 \cite{KacW} (see also \cite{PotWar}
where the calculation of the correlation functions we discuss below is addressed by this method). 
Kac and Ward did not intend to
present a rigorous derivation; the missing proofs were given later by Sherman \cite{Sher} and Burgoyne \cite{Bur} following a conjecture by Feynman. A rigorous combinatorial solution by a related method was given by
Hurst and Green \cite{HurGr}. Another combinatorial solution was obtained by Kasteleyn \cite{Kas1}.
A useful modification of Kasteleyn's approach is due to Fisher \cite{Fis3} (see \cite{McWu1}). 
One more simple solution based on the method of Kac and Ward was given by  
Vdovichenko \cite{Vdov} (see \cite{LanLif}) although \cite{Vdov} omits the discussion of the boundary conditions.
In \cite{Kas1}, Kasteleyn maps the Ising 
problem onto the dimer covering problem (i.e. the problem of counting all possible coverings by dimers)
on a special lattice, slightly more complex than the square lattice. It is interesting to note that the 
dimer problem on the square lattice, which thus can be regarded as a simpler version of the Ising problem,
was solved by Kasteleyn himself \cite{Kas2} and, independently, by Temperley and Fisher \cite{TemFis} in 1961
(see also \cite{LieLoss} for a later alternative, particularly simple, solution of the dimer problem on the square lattice).

We now return to 1946.      
Kaufman's method and results were eventually published three years later in \cite{Kau1}.  Onsager
continues in \cite{Ons2}: ``By itself that was only a more elegant derivation of an old result;
but the approach looked powerful enough to produce a few new ones.  Very well, how about correlations?''
As noted
by Baxter, the history of the Ising model from that time forth has been the study of these correlations.

In 1948, at a conference at Cornell, and again in 1949 at a conference in Florence, Onsager astounded
the audience by announcing that he and Kaufman had recently obtained an exact formula for the spontaneous
magnetization $M_0$ of the 2-D Ising model.  He gave the result as
\be\label{eq15}
M_0 = \lf(1-k^2_{\tons}\rt)^{1/8}, \qquad k_{\tons} \equiv \lf( \sinh \frac{2J_1}{k_B T}\;
\sinh \frac{2J_2}{k_B T}\rt)^{-1}
\ee
for $0 < k_{\tons} <1$,
corresponding to  $T< T_c$.  And for $k_{\tons} >1$, corresponding to
$T> T_c$, $M_0 =0$.

Onsager and Kaufman never published a proof of this result.  The situation was as follows.  In
1949 they published their famous paper \cite{KauOns},  based on Kaufman's spinor (free-fermion)
approach, in which they found, in particular, an explicit expression (see \eqref{eq17} below) for the
two-point correlation function along a row in the thermodynamic limit, i.e.\
$\lf\langle \sg_{1,1} \sg_{1,1+n}\rt\rangle$.
By physical arguments
\be\label{eq16}
M_0 = \lim_{n\to \infty}  \lf\langle \sg_{1,1} \, \sg_{1, 1+n}\rt\rangle^{\frac{1}{2}}
\ee
so that in order to derive \eqref{eq15}, what they had to do was to control the limit in
\eqref{eq16}.  As Onsager describes the situation \cite{Ons3}: ``And then, finally, it became
apparent that we had a last problem, of the degree of the order, and that turned out to need
more invention.''

\section{Toeplitz meets Ising}
What precisely was this ``last problem'' that Kaufman and Onsager had to face in order to compute the limit in 
\eqref{eq16}?  It is here that we
make contact with the theory of Toeplitz matrices.  In \cite{KauOns} Kaufman and Onsager expressed
$\langle \sg_{1,1} \, \sg_{1,1+n}\rangle$ as a sum of two Toeplitz determinants as follows
(in the notation of \cite{MPW}):
\be\label{eq17}
\begin{aligned}
(-1)^n \lf\langle \sg_{1,1} \,\sg_{1, 1+n}\rt\rangle &= c^{*\,2}_2 \det
\begin{pmatrix}
b_{-1}  & b_{-2} & b_{-3}  & \dots & b_{-n} \\
b_{0}  & b_{-1} & b_{-2}  & \dots & b_{1-n} \\
b_{1}  & b_{0} & b_{-1}  & \dots & b_{2-n} \\
\vdots & \vdots & \vdots & & \vdots\\
b_{n-2}  & b_{n-3} & b_{n-4}  & \dots & b_{-1}
\end{pmatrix} \\[6mm]
&-s^{*2}_2 \det
\begin{pmatrix}
b_1  & b_2 & b_3  & \dots & b_n \\
b_{0}  & b_1 & b_2  & \dots & b_{n-1} \\
b_{-1}  & b_{0} & b_{1}  & \dots & b_{n-2} \\
\vdots & \vdots & \vdots &\vdots\;\vdots\;\vdots & \vdots\\
b_{2-n}  & b_{3-n} & b_{4-n}  &  & b_{1}
\end{pmatrix}
\end{aligned}
\ee
where
\be\label{eq18}
c^*_2 = \cosh\, K^*_2, \quad s^*_2 = \sinh \, K^*_2, \quad e^{-2K^*_2} = \tanh\, \frac{J_2}{k_B T}
\ee
\be\label{eq19}
b_j = \frac{1}{\pi} \int^{\pi}_0 \cos \lf( {\hde} (\tht) - j \tht \rt) d\tht
\ee
\be\label{eq20}
\cos \hde(\tht) =  \sin \lf(\de^*(\tht) \rt) \sin \tht \cosh 2K^*_2 - \cos \lf(\de^* (\tht)\rt)
\cos \tht
\ee
and $\varphi_{\tons}\! \lf(\eit\rt) = e^{i\de^*(\tht)}$ is a function introduced previously  by Onsager in
\cite{Ons1}
\be\label{eq21}
\varphi_{\tons}\! \lf(\eit\rt)= \lf[ \lf( \frac{1-\gamma_1 \, \eit}{1-\gamma_1 \, e^{-i\tht}}\rt)
\lf( \frac{1-\gamma_2 \, e^{-i\tht}}{1-\gamma_2 \, \eit}\rt) \rt]^{\frac{1}{2}}
\ee
\be\label{eq22}
\begin{aligned}
&\text{where }\qquad \qquad \gamma_1 = z_1 \,z^*_2, \qquad \gamma_2 =z^*_2/z_1 \\
&z_1 = \tanh \frac{J_1}{k_B\, T}, \quad z_2 = \tanh \frac{J_2}{k_B \,T}, \quad z^*_2 =
\frac{1-z_2}{1+z_2},
\end{aligned}
\ee
The branch in \eqref{eq21} is chosen so that $\varphi_{\tons}\! \lf(e^{i\pi}\rt) =e^{i\de^*(\pi)} >0$.

\begin{remark}\label{remark1}
As noted in \cite{MPW}, there is a sign error in formula \eqref{eq17} in \cite{KauOns}.
\end{remark}

So the problem that Kaufman and Onsager had to face to compute $M_0$ via \eqref{eq16},
was the purely mathematical problem of computing the asymptotics of $n\times n$ Toeplitz
determinants as $n\to \infty$.  At the time, however, all that was known was Szeg\H{o}'s result
\eqref{eq5} with unknown error term $o(n)$.  But  to compute \eqref{eq16}, this error term is precisely 
what one needs to know! 

Note that
\be\label{eq23}
0 < \gamma_1 < 1 \qquad\text{ and }\qquad \gamma_1 < \gamma_2.
\ee
If $\gamma_2 \neq 1$, we see that
\be\label{eq24}
\varphi_{\tons}\! \lf(\eit\rt) = \frac{\lf(1-\gamma_1 \, \eit\rt) \lf(1-\gamma_2 \,e^{-i\tht}\rt)}{
\lf|\lf(1-\gamma_1 \, e^{-i\tht}\rt) \lf(1-\gamma_2\, \eit\rt)\rt|}
\ee
and hence $\varphi_{\tons}\! \lf(\eit\rt)$ is a smooth, single-valued function on the unit
circle $S^1$.  On the other hand, if $\gamma_2=1$, then
\be\label{eq25}
\varphi_{\tons}\! \lf(\eit\rt) = i\, e^{-i\tht/2} \,\frac{\lf(1-\gamma_1 \, \eit\rt)}{\lf|1-\gamma_1 \, e^{-i\tht}\rt|}, \qquad
\qquad 0< \tht < 2\pi.
\ee
Simple algebra shows that for $\hka_i = J_i/k_BT, \;i=1,2$,
\begin{align*}
\gamma_2 \gtrless 1 \Longleftrightarrow  \sinh \lf(\hka_1 + \hka_2 \rt)
- \cosh \lf(\hka_1-\hka_2\rt)
\lessgtr 0 \\
\intertext{and}
\sinh^2 \lf(\hka_1 + \hka_2 \rt) - \cosh^2 \lf(\hka_1 - \hka_2\rt) = k^{-1}_{\tons} -1
\end{align*}
where $k_{\tons} = \lf(\sinh\, 2\hka_1 \;\sinh\,2\hka_2\rt)^{-1}$ as in \eqref{eq15}.  Thus
\be\label{eq26}
\gamma_2 \gtrless 1 \Longleftrightarrow  k_{\tons} \gtrless 1.
\ee
In other words, sub-critical (resp.\ super-critical) temperatures, $T< T_c$ (resp.\ $T>T_c$)
correspond to $z^*_2 < z_1$ (resp.\ $z^*_2 > z_1$).  And, of course, $T=T_c \Longleftrightarrow
z^*_2 = z_1$.

It turns out that in their initial (unpublished) calculations to compute $M_0$, Kaufman and
Onsager did not actually use \eqref{eq17}.  What happened is described by Onsager in
\cite{Ons2} and
\cite{Ons3} (see also \cite{Bax} for more details, particularly concerning the Wiener-Hopf
calculation below).  In addition to \eqref{eq17}, Kaufman and Onsager had derived, but
apparently did not publish (see private communication to C.~Domb \cite[p.~201]{Dom2}),
an expression for the 2-point correlation function along a diagonal
$\langle \sg_{1,1}, \sg_{1+n,\, 1+n}\rangle$, which is much simpler than \eqref{eq17}.
Only one Toeplitz determinant is involved,
\be\label{eq27}
\langle \sg_{1,1}\, \sg_{1+n,\, 1+n}\rangle = D_n \lf(\varphi_{\tdg}\rt)
\ee
where\footnote{In \cite{Dom2}, the symbol is $\overline{\varphi_{\tdg}}$, not $\varphi_{\tdg}$.
However, the determinant is the same, 
$D_n(\overline{\varphi_{\tdg}})=\det T_n(\overline{\varphi_{\tdg}})=
\det\overline{T_n(\varphi_{\tdg})^T}=\overline{\det T_n(\varphi_{\tdg})^T}=
\det T_n(\varphi_{\tdg})^T=\det T_n(\varphi_{\tdg})=D_n(\varphi_{\tdg})$.
}
\be\label{eq28}
\varphi_{\tdg}\!\lf(\eit\rt) = \lf( \frac{1-k_{\tons} \, e^{-i\tht}}{1-k_{\tons}
\,e^{i\tht}} \rt)^{\frac{1}{2}}
\ee
and again $k_{\tons}$ is given as in \eqref{eq15} and $\varphi_{\tdg} (e^{i\pi})>0$.
The same physical arguments as in \eqref{eq16} yield
\be\label{eq29}
M_0 =\lim_{n\to \infty} \lf\langle \sg_{1,1} \, \sg_{1+n, 1+n} \rt\rangle^{\half}.
\ee
Onsager then realized that the eigenvalue  equation for the Toeplitz matrix $T\!\lf(\varphi_{\tdg}\rt)$,
\[
u_k = \frac{1}{\lambda} \sum^{n-1}_{j=0} \lf(\varphi_{\tdg}\rt)_{k-j} \, u_j,\qquad 0\le k \le n-1
\]
was a discrete analog of the Milne integral equation
\be\label{eq30}
f(x)  = \int^\infty_0 A(x-y)\, f(y)dy, \qquad x>0
\ee
which arises, with a particular, explicit function $A(\cdot)$,  in radiative equilibrium
theory (see, for example \cite{DymMc}).  Milne's equation can be solved using the Wiener-Hopf
method, a technique with which Onsager was familiar from lectures that Wiener had given earlier
in the Mathematics Department at Yale.
So Onsager tried the Wiener-Hopf technique on the above eigenvalue
equation for $T (\varphi_{\tdg})$,
to obtain the eigenvalues $\lb=\lb^{(n)}_\ell$, $1\le \ell \le n $, and then computed
$D_n (\varphi_{\tdg})= \prod^n_{\ell=1} \, \lb^{(n)}_\ell$.  Letting $n\to \infty$, formula \eqref{eq15}
for $M_0$ then emerged.  This was the basis for the first announcements of the result in 1948 and 1949.  So
why didn't Kaufman and Onsager publish the details of their calculation?  The situation was as follows.
As in the case of Milne's equation, the Wiener-Hopf method requires very detailed knowledge of the kernel
$(\varphi_{\tdg})_{j-k}$, in the case of Ising.  Onsager sensed that there was another prize to be had:
What about Toeplitz determinants with general symbols $\varphi$?
And so, before  Kaufman and Onsager could get around to publishing their calculation, Onsager began
looking for a method to evaluate asymptotically Toeplitz determinants with general symbols.
We quote Onsager \cite{Ons3}:
`` \dots and lo and behold I found it.   It was a general formula for the evaluation of Toeplitz matrices.
The only thing I did not know was how to fill out the holes in the mathematics and show the epsilons and deltas
and all of that, and the limiting processes; I did not know just how it should be done and what mathematicians
really knew about limiting processes in that ball park.'' As it turned out, by that point, mathematicians knew
a lot!  Before it was clear what conditions to place on the symbol, Kaufman and Onsager spoke to Kakutani
and then Kakutani spoke to Szeg\H{o}.  Szeg\H{o} then revisited his calculations from 1915, and  in  1952
\cite{Sz3} evaluated the $o(n)$ term in \eqref{eq5} for a very general class of symbols $\varphi$.  This
celebrated result is the Szeg\H{o} Strong Limit Theorem mentioned above.
Faced with Szeg\H{o}'s general result, Kaufman and Onsager published neither their first, nor their
second method.  As Onsager noted \cite{Ons2}, `` \dots the mathematicians got there first.''
In the long history of mathematics and physics, it is most unusual for a physicist to be scooped
out of a formula by a mathematician!

There is, however, more to the story.  On  the web (see \cite{Bax}) there is the draft of a paper titled
``Long-Range Order''
from 1950 (never published), without names attached but almost certainly by Onsager and Kaufman, 
that describes their more
general method and obtains the result \eqref{eq15} for $M_0$.  Here the authors use formulae
\eqref{eq16} and \eqref{eq17}.  There is also an earlier letter \cite{Kau2} written on April 12, 1950,
from Onsager to Kaufman
containing some of the calculations from the above draft.  Presumably to allay any concerns about
priority,  Onsager assures Kaufman at the end of the letter that ``There will be time for all these
things.''  As one reads \cite{Ons2} \cite{Ons3}, one senses the faint hint of regret that Onsager must
have felt for such optimism.  In a letter from Kaufman to Onsager, written on May 12, 1950, a month after the letter
\cite{Kau2}, Kaufman says ``Here is a draft of Crystal Statistics IV''.  In \cite{Bax}, Baxter makes
a strong case that the above paper from 1950 is indeed Kaufman's draft of 
Crystal Statistics IV. However, this may not be so. On the Onsager archive in Trondheim 
(see footnote \ref{foot3}) under 
``Selected research material and writings'', item 17.121, one finds the paper ``Long-Range order'' discussed by Baxter. 
Although, again, no names are attached, this establishes the provenance of the paper beyond any reasonable doubt,
as suspected by Baxter. 
 Item 17.121, however, also contains a second, completely separate paper called ``Crystal statistics IV. Long-Range Order in a Binary Crystal'', and both Kaufman's and Onsager's names are attached.\footnote{The authors thank H. Holden for locating this paper on the Trondheim archive.}
The ``Long-range'' paper describes, as noted above, the second method of Kaufman-Onsager for general symbols, and the second paper describes their first Wiener-Hopf type method. However, it is unfinished and it ends abruptly. So when Kaufman wrote to Onsager 
``Here is the draft of Crystal Statistics IV'', which of these two drafts did she mean? Either way, the evidence that Kaufman/Onsager knew how to  evaluate the asymptotics of Toeplitz determinants (in two ways!) is now clear and written down in their own hand.

\section{Szeg\H{o}'s Theorem}\label{secSzego}
Here is Szeg\H{o}'s result \cite{Sz3} \cite{GreSz}.

\begin{theorem}[Szeg\H{o} Strong Limit Theorem]\label{theorem3}
Let $\varphi\!\lf(\eit\rt)$ be a positive, $C^{1+\varepsilon}$, $\varepsilon>0$, function on $S^1$.  
Let $(\log \varphi)_k = \frac{1}{2\pi} \int^\pi_{-\pi} e^{-ik\tht} \log \varphi
\!\lf(\eit\rt) d\tht$, $k\in \ZZ$, denote the Fourier coefficients of $\log \varphi\!\lf(\eit\rt)$.  Then
\be\label{eq31}
\lim_{n\to \infty} \frac{D_n (\varphi)}{e^{n(\log\varphi)_0}} = e^{E(\varphi)}
\ee
where
\be\label{eq33}
E(\varphi) = \sum^\infty_{k=1} k \lf| \lf(\log \varphi\rt)_k \rt|^2.
\ee
\hfill $\Box$
\end{theorem}

\noindent
Thus the $o(n)$  term in \eqref{eq5} is given by $E(\varphi) + o(1)$.

\begin{remark}
Note that
\be\label{eq32}
e^{(\log\varphi)_0}= \exp\lf[\frac{1}{2\pi} \;\int^\pi_{-\pi} \log\varphi\!\lf(\eit\rt)  d\tht\rt]
\ee
is the geometric mean of $\varphi$.
\end{remark}

In contrast to \eqref{eq31}, Szeg\H{o}'s earlier result Theorem \ref{theorem1},
$\lim_{n\to\infty} \lf(D_n(\varphi)\rt)^{\frac{1}{n}} = e^{(\log\varphi)_0}$, is known simply as
Szeg\H{o}'s First Theorem, or sometimes just Szeg\H{o}'s Theorem.

In the years that followed, there was considerable effort by many mathematicians to weaken the
assumptions in Theorem \ref{theorem3}.  Along the way, many new methods, quite different one
from the other, were discovered to prove \eqref{eq31}.  These developments are described in the
outstanding recent monographs by Barry Simon \cite{Sim2} and 
B\"ottcher and Silbermann \cite{BottSilb3}: also see \cite{McC2} and \cite{Bott1}
for earlier summaries.  In \cite{Sim2}, the contributions of Kac \cite{Kac}, G. Baxter \cite{BaxG1}, 
Hirschman \cite{Hir}, and Devinatz \cite{Dev},
in particular, are discussed, leading up through successively weaker assumptions 
to the definitive result of Ibragimov  \cite{Ibr} in 1968
giving necessary and sufficient conditions on $\varphi\!\lf(\eit\rt)$ for \eqref{eq31} to hold
under the assumption of the positivity of $\varphi$.
\begin{theorem}[Ibragimov]\label{theorem4}
Let $\varphi\!\lf(\eit\rt) \frac{d\tht}{2\pi}$ be a probability measure on $S^1$ and suppose
that $\log\varphi\!\lf(\eit\rt)$ is integrable.  Then \eqref{eq31} is always true in the following
sense:  $\lim_{n\to \infty} \frac{D_n(\varphi)}{e^{n(\log\varphi)_0}}$ always exists and equals $e^{E(\varphi)}$,
including the case where one, and hence both, are infinite.
\hfill $\Box$
\end{theorem}

One can consider Toeplitz matrices $T_n(d\mu)$ for general probability measures $d\mu(\tht)=
\varphi\!\lf(\eit\rt) \frac{d\tht}{2\pi} +d\mu_s(\tht)$ on $S^1$, where $d\mu_s$ denotes the singular
part of $d\mu$.  We have $T_n (d\mu) = \{ \mu_{j-k}\}_{0\le j,\,k\le n-1}$, where $\mu_\ell= \int e^{-i \ell
\tht} \, d\mu(\tht)$, $\ell \in \ZZ$ and $D_n(d\mu)= \det \, T_n(d\mu)$. 
The following remarkable result generalizes Theorem \ref{theorem1}.

\begin{theorem}[\cite{Sz1}\cite{Sz2}\cite{Sz4}\cite{Ver}]\label{theorem5}
Let $d\mu(\tht) = \varphi\! \lf(\eit\rt) \, \frac{d\tht}{2\pi} + d\mu_s(\tht)$ be a probability  measure as above.
Then Szeg\H{o}'s Limit Theorem is always true in the following sense:  
$\lim_{n\to\infty} D_n (d\mu)^{\frac{1}{n}}$ always exists and equals 
$\exp\lf[\int^\pi_{-\pi} \log \varphi\!\lf(\eit\rt)\, \frac{d\tht}{2\pi}\rt]$,
including the case where one, and hence both, are zero.
\hfill $\Box$
\end{theorem}

Szeg\H{o} proved the result when $d\mu_s=0$, whereas Verblunsky was able to handle the case $d\mu_s \neq 0$.
Of course, as $\varphi\!\lf(\eit\rt) \in L^1 \lf(\frac{d\tht}{2\pi}\rt)$, $\int^\pi_{-\pi} \log
\varphi\!\lf(\eit\rt) \, \frac{d\tht}{2\pi}
< \infty$, and the integral can diverge, possibly, only to $-\infty$.  The striking feature of
Theorem~\ref{theorem5} is that $\lim_{n\to\infty} D_n (d\mu)^{\frac{1}{n}}$ is independent of $d\mu_s$.

Is it possible that \eqref{eq31} remains true when $d\mu_s \neq 0$?  Here the definitive result is due to
Golinskii and Ibragimov.

\begin{theorem}[\cite{GolIbr}]\label{theorem6}
Let $d\mu(\tht)=\varphi\!\lf(\eit\rt) \, \frac{d\tht}{2\pi} + d\mu_s(\tht)$ be a probability measure as above,
and suppose $\log \varphi\!\lf(\eit\rt) \in L^1\lf(\frac{d\tht}{2\pi}\rt)$.  If $d\mu_s\neq 0$, then
\be\label{eq34}
\lim_{n\to\infty} \frac{D_n (d\mu)}{\exp\lf[n \int^\pi_{-\pi} \log \varphi\!\lf(\eit\rt)\,\frac{d\tht}{2\pi}\rt]}
=   +\infty
\ee
\hfill $\Box$
\end{theorem}

In other words, if \eqref{eq31} were to hold with $E(\varphi) =\sum^\infty_{k=1} k \lf|(\log \varphi)_k \rt|^2
< \infty$, then necessarily $d\mu_s =0$.

In \cite{Sim2}, Simon presents six different proofs of \eqref{eq31}.  Some of these proofs require stronger
conditions on $\varphi$ (see \cite{Sim2}).  We now describe these methods in some detail 
(but without giving precise conditions on $\varphi$ in each case)
in order to
illustrate the extraordinary variety of mathematical areas and techniques that inter-relate with Szeg\H{o}'s
Strong Limit Theorem.  Our presentation follows \cite{Sim2}. After this presentation, we will also describe 
three additional proofs not covered in \cite{Sim2}, together with some brief comments on Szeg\H{o}'s 
original proof (see Remark \ref{22+++}). 

The first proof uses ideas and results from the theory of orthogonal polynomials on the unit circle
(OPUC's).  These polynomials were introduced by Szeg\H{o} in the early 1920's \cite{Sz4} \cite{Sz5}
in the course of his investigation of the eigenvalues $\lb^{(n)}_1, \dots, \lb^{(n)}_n$ of Toeplitz
matrices $T_n$.   The classical references for OPUC's are \cite{Sz5} \cite{Ge}.  The  OPUC's $p_k(z) = \chi_k \,
z^k + \dots, \chi_k >0$, associated with a probability measure $d\mu(\tht)$ on $S^1$ are formed
by orthonormalizing $1, z, z^2, \dots, z^k, \dots$ with respect to $d\mu(\tht)$ in $L^2\lf(S^1,\,
d\mu(\tht)\rt)$,
\be\label{eq35}
\int_{S^1} \ovl{p_k\! \lf(\eit\rt)} \, p_j\!\lf(\eit\rt) \, d\mu(\tht) = \de_{k,j}, \qquad k,\, j \ge 0 \ .
\ee
OPUC's are intimately related to Toeplitz determinants.  For example
\be\label{eq36}
\frac{D_{n+1} (d\mu)}{D_n (d\mu)} =\frac{1}{\chi_n^2}= \frac{1}{\chi_0^2}
\prod^{n-1}_{j=0} \lf(1-|\xi_j|^2\rt), \qquad n\ge 1
\ee
where $\xi_j \equiv -\lf(\chi_{j+1}\rt)^{-1} \;\ovl{p_{j+1} (0)}$, $\; j\ge 0$, are  the so-called
\tit{Verblunsky} coefficients.   The $\xi_j$'s have modulus less than 1, and so $F(d\mu)\equiv \lim_{n\to \infty}
D_{n+1} (d\mu)/D_n(d\mu)$ always exists and equals 
$\chi_0^{-2}\prod^\infty_{j=0} \lf(1-|\xi_j|^2\rt)$.  But then
$\lim_{n\to\infty} \lf(D_n(d\mu)\rt)^{\frac{1}{n}}$ always exists, and has the same limit.  So the proof
of Szeg\H{o}'s Limit Theorem, boils down to showing that
\be\label{eq37}
\frac{1}{\chi_0^2}\prod^\infty_{j=0} \lf(1-|\xi_j|^2\rt) = \exp\lf[\frac{1}{2\pi}
\int^\pi_{-\pi} \log \varphi\!\lf(\eit\rt) d\tht\rt]
\ee
where $d\mu(\tht)=\varphi\!\lf(\eit\rt)\, \frac{d\tht}{2\pi} + d\mu_s(\tht)$.  In turn, one finds that
\be
\lim_{n\to\infty} \frac{D_n(d\mu)}{\lf(F(d\mu)\rt)^n} = \lim_{n\to\infty} 
\frac{\prod^{n-2}_{j=0} \lf(1-|\xi_j|^2\rt)^{-j-1}}{\prod_{j=n-1}^\infty (1-|\xi_j|^2)^n}
\ee
and the proof of SSLT for $d\mu_s=0$ reduces to showing that the limit on the RHS is just
$e^{E(\varphi)}$.  As we will see below, OPUC's play a crucial role in analyzing the Fisher-Hartwig
conjecture for Toeplitz determinants with singular symbols $\varphi\!\lf(\eit\rt)$.

The second proof utilizes a remarkable identity which gives an exact formula for $D_n(\varphi)$,
\be\label{eq38}
D_n(\varphi) = e^{n(\log\varphi)_0 + E(\varphi)} \, \det \lf(1-Q_n \, H(b)\, H(\widetilde{c})\, Q_n\rt),
\quad n\ge 1
\ee
where $\det$ denotes the determinant in $\ell^+_2 = \ell_2(\ZZ_+)$, $\ZZ_+=\{0,1,\dots\}$, $Q_n$ is the orthogonal
projection in $\ell^+_2$ onto $\ell^+_{2,n} = \{u= (u_0, u_1, \dots ) \in \ell^+_2\,:\; u_i=0
\text{ for } 0\le i <n\}$,  and $H(b)$, $H(\widetilde{c})$ are the Hankel operators with kernels
$\{b_{i+j+1}\}_{i,\, j\ge 0}, \, \{c_{-i -j-1}\}_{i,\, j \ge 0}$.
Here $b\!\lf(\eit\rt) = \sum_j b_j\, e^{ij\tht}$, $c\!\lf(\eit\rt) = \sum_j c_j \, e^{ij\tht}$ are defined in terms of
the \tit{Szeg\H{o} function} ${\cal D}(z) = \exp \lf(\frac{1}{4\pi} \int^\pi_{-\pi}  \frac{e^{i\tht}+z}{e^{i\tht}-z}
\;\log \varphi\!\lf(\eit\rt)d\tht\rt)$,
\be\label{eq39}
b\equiv \frac{\ovl{\cal D}}{\cal D},\qquad c\equiv \frac{\cal D}{\ovl{\cal D}}
\ee
where ${\cal D}={\cal D}(e^{i\tht})=\lim_{z\to e^{i\tht},|z|<1}{\cal D}(z)$.
The point is that $H(b)\, H(\widetilde{c})$ is a trace-class operator in $\ell^+_2 = \ell^+_{2,\,0}$, and hence
$Q_n\, H(b)\, H(\widetilde{c})\,Q_n~\to~0$ in trace-norm by abstract arguments as $n\to\infty$.
But then $\det \lf(1-Q_n \, H(b)\, H(\widetilde{c}) \, Q_n\rt)\to 1$ as $n\to\infty$,  by the continuity
properties of the determinant.  (For more information about the trace class and determinants, see \cite{BottSilb3}
\cite{Sim3}.)
Thus SSLT is immediate, once one has proved \eqref{eq38}.   Formula \eqref{eq38}
is known as the Borodin-Okounkov formula and was proved \cite{BorOk} in 2000 in the context of
combinatorics and random matrix theory.  The formula evoked great interest, and meanwhile several new
proofs and generalizations have been given (see, e.g., \cite{BasWid} \cite{BottWid1} and \cite{Sim2}).
However, it turns out that it
was already proven many years earlier by Geronimo and Case in 1979 in their work on inverse scattering theory
\cite{GerCase}.  The  broader significance of this paper was not appreciated at the time.  But more to the
point, in a separate section in their paper titled ``Szeg\H{o}'s Theorem'', Geronimo and Case actually used
\eqref{eq38} to prove SSLT.  It is unfortunate that these developments were overlooked by the experts in the field.

The third and fourth methods utilize the following Heine-type multi-integral representation for $D_n(\varphi)$ (see, e.g., \cite{BottSilb3} \cite{Sim2})
\be\label{eq40}
D_n(\varphi) = \frac{1}{n!} \int^{2\pi}_0 \dots \int^{2\pi}_0 \prod_{0\le j < k \le n-1}
\lf|e^{i\tht_j} - e^{i\tht_k} \rt|^2 \prod^{n-1}_{j=0} \varphi\lf(e^{i\tht j}\rt)
\,\frac{d\tht_j}{2\pi} \ .
\ee
In the third method, due to Bump and Diaconis \cite{BumpDia}, the authors observe that the RHS of \eqref{eq40}
can be re-written via Weyl's integration formula as
$\int_{\MCU(n)} e^{F_n(g)} \, dg$, where $dg$ denotes Haar measures on the unitary group $\MCU(n)$, and
$F_n(g) = \sum^{n-1}_{j=0} \log\varphi\!\lf(e^{i\tht_j(g)}\rt)$, where  $\lf\{e^{i\tht_j(g)}\rt\}$ are the
eigenvalues of $g$.  Substituting the Fourier expansion $\log \varphi \!\lf(\eit\rt) = \sum^\infty_{k=-\infty}
(\log\varphi)_k \, e^{ik\tht}$, one obtains
\be
D_n(\varphi) = e^{n (\log \varphi)_0} \int_{\MCU(n)} \exp \lf(\sum_{k\neq 0}
(\log \varphi)_k\tr(g^k)\rt) dg
\ee
and after expanding out the exponential, one ends up with the representation
\be\label{eq41}
D_n(\varphi) = e^{n(\log \varphi)_0} \int_{\MCU(n)} \ovl{\lf(\sum_t \eta_t \, T_t(g)\rt)}
\lf(\sum_s \eta_s \, T_s(g)\rt) dg \ .
\ee
Here the sums are over all $k$-tuples of non-negative integers $t=(t_1, \dots, t_k)$,
$s= (s_1, \dots, s_k)$, and for any $k$, $T_t(g)=\prod^k_{j=1} \lf[\tr \lf(g^j\rt)\rt]^{t_j}$,
and the $\eta_t$'s, $\eta_s$'s are explicit constants depending on $\{(\log \varphi)_\ell\}$.
The magic of the method is that one can use the theory of characters together with Frobenius-Schur-Weyl
duality to compute $\int_{\MCU (n)} \ovl{T_t(g)} \, T_s (g) \, dg$ explicitly, obtaining
\be\label{eq42}
\int_{\MCU (n)} \ovl{T_t(g)}\, T_s (g) = \de_{ts}\; W(t)
\ee
where $W(t)= \prod^k_{j=1} \, j^{t_j}\, t_j!$, provided $\sum^k_{j=1} j\,t_j$ and $\sum^k_{j=1} \, j \, s_j$
are $\le n$.  Substituting \eqref{eq42} into \eqref{eq41} and controlling the errors as $n\to \infty$, one
immediately obtains \eqref{eq31}.

In the fourth proof, due to Johansson \cite{Joh1}, one rewrites \eqref{eq40} with $\varphi\! \lf(\eit\rt)
=e^{V(\tht)}$
in the form
\be\label{eq43}
D_n \!\lf(e^{V}\rt) = \frac{1}{n!} \int e^{-2H_n \lf(\tht_{0}, \dots, \tht_{n-1}\rt)} \: \frac{d\tht_0}{2\pi}
\dots \frac{d\tht_{n-1}}{2\pi}
\ee
where
\be\label{eq44}
H_n = -\half\sum^{n-1}_{j=0} V(\tht_j) +\half\sum_{0\le k< j \le n-1}  W\lf(\tht_k-\tht_j\rt)
\ee
and
\be\label{eq45}
W(\tht) = -\log \lf|e^{i\tht}-1\rt|^2 \ .
\ee


This expression exhibits $n! D_n\!\lf(e^{V}\rt)$ as the partition function of a gas of 1-dimensional particles
with 2-dimensional Coulomb interactions --- a ``log-gas'' in the terminology of Dyson. Dyson introduced this terminology in the context of his analysis \cite{Dy1} of Circular Ensembles of $n\times n$ random matrices. For such ensembles Dyson found that the distributions of the eigenvalues $\{e^{i\tht_j}\}_{j=0}^{n-1}$,
$0\le\tht_j<2\pi$, of the matrices were given by 
\be
\frac{1}{Z_{n,\beta}}\prod_{0\le j<k\le n-1} |  e^{i\tht_j}-e^{i\tht_k} |^\beta \frac{d\tht_0}{2\pi}\cdots
\frac{d\tht_{n-1}}{2\pi}.
\ee
Here $\beta=1$ for the orthogonal, $\beta=2$ for the unitary, and $\beta=4$ for the symplectic ensemble,
and 
\be\label{partDy}
Z_{n,\beta}=\int_0^{2\pi}\frac{d\tht_0}{2\pi}\cdots \int_0^{2\pi}\frac{d\tht_{n-1}}{2\pi}
\prod_{0\le j<k\le n-1} |  e^{i\tht_j}-e^{i\tht_k} |^\beta
\ee
is the normalization constant. Dyson then noted that the distribution of the eigenvalues $\{e^{i\tht_j}\}_{j=0}^{n-1}$ was identical with the distribution of the positions $0\le\tht_j<2\pi$, $j=0,\dots,n-1$ of charges in a finite gas with 
2-D Coulomb interactions at inverse temperature $\beta$ and with Gibbs measure
$\frac{1}{Z_{n,\beta}}\exp\{-\beta H_n(\tht_0,\dots,\tht_{n-1})\} \frac{d\tht_0}{2\pi}\cdots\frac{d\tht_{n-1}}{2\pi}$,
where $H_n(\tht_0,\dots,\tht_{n-1})=-\sum_{0\le j<k\le n-1}\log  |  e^{i\tht_j}-e^{i\tht_k} |$, and the partition function is given by \eqref{partDy}. As emphasized by Dyson \cite{Dy1}, a consequence of this identification is that ``\dots the thermodynamic notions of entropy, specific heat, etc., can be transferred from the Coulomb gas to the eigenvalue series. This will prove very useful, as it gives us a precise and well-understood language in which to describe the statistical properties of the eigenvalues''. If the matrices $U$ in the ensembles carry a weight $e^{-\beta\tr Q(U)}$, then the energy in the corresponding Coulomb gas must be replaced by 
\be
H_n(\tht_0,\dots,\tht_{n-1})=-\sum_{0\le j<k\le n-1}\log  |  e^{i\tht_j}-e^{i\tht_k} |+\sum_{j=0}^{n-1}Q(e^{i\tht_j})
\ee
where $Q$ now plays the role of an external force acting on the charges. 

In operational terms, the identification implies, in particular, (see \cite{Dy2}) that the free energy $F\equiv -\beta^{-1}\log Z_{n,\beta}$ is approximated to leading order as $n\to\infty$ by the classical thermodynamic recipe, $F\sim F_T$, where
\be\label{5.1}
\beta F_T=G_2+G_1+G_0.
\ee
Here $G_2$ is the macroscopic Coulomb energy
\be\label{5.2}
G_2=-\frac{\beta}{2}\int\int\sg(e^{i\tht})\sg(e^{i\tht'})\log |e^{i\tht}-e^{i\tht'}|d\tht d\tht'+
\beta\int Q(\eit)\sg(\eit)d\tht
\ee
and $G_1+G_0$ is the local free energy term
\begin{align}
G_1&=\lf(1-\frac{\beta}{2}\rt)\int\sg(\eit)\log\sg(\eit)d\tht,\label{5.3}\\
G_0&=n\lf[ \lf(1-\frac{\beta}{2}\rt) \log {2\pi\over n} +\log\Gamma\lf(1+\frac{\beta}{2}\rt)\rt]
-\log\Gamma\lf(1+\frac{\beta}{2}n\rt).\label{5.4}
\end{align}
The function $\sg(\eit)\ge 0$ is the macroscopic density of the gas and is chosen to minimize the functional on the right-hand side of \eqref{5.1} subject to the constraint
\be\label{6.1}
\int\sg(\eit)d\tht=n.
\ee
Dyson first applied these thermodynamic notions to compute  \cite{Dy2} the probability  $P_s$ that there are no eigenvalues for a random matrix in the interval $(0,2s/\pi)$, in the bulk scaling limit
where the average distance between the eigenvalues is rescaled to 1,
 for the Gaussian Unitary Ensemble (see \cite{Meh}). Using these notions, he derived, for the first time, the leading order of $P_s$,
$P_s\sim \exp\{-s^2/2\}$ as $s\to\infty$. (We return to this problem further on (see \eqref{eq133} et seq).)
A year later in \cite{Dy3}, Dyson showed how these notions play out in the context of SSLT. First Dyson noted, as explained above, that the Heine-type representation \eqref{eq40} immediately displayed 
$n! D_n\!\lf(e^{V}\rt)$ as the partition function of a log-gas as in \eqref{eq43}--\eqref{eq45} above.
This connection between Toeplitz determinants and the Coulomb gas was first pointed out in an unpublished manuscript by G. Kreisel, which was made available to Dyson by courtesy of E. Wigner.\footnote{We note that already in the 1880's Stieltjes gave an electrostatic interpretation of the zeros of certain classical orthogonal polynomials
(see \cite{Sz5}, Section 6.7). For a discussion of logarithmic potentials from the viewpoint of potential theory, see the book of Saff and Totik \cite{SaTo}.}

For $\beta=2$ the term $G_1=0$ and the constant $G_0=-\log n!$, and we have
\be
F_T=-\frac{1}{2}\int\int\sg(e^{i\tht})\sg(e^{i\tht'})\log |e^{i\tht}-e^{i\tht'}|d\tht d\tht'+
\int Q(\eit)\sg(\eit)d\tht-\half \log n!
\ee
The constrained minimum is obtained when
\be
Q(\eit)- \int\sg(e^{i\tht'})\log|e^{i\tht}-e^{i\tht'}| d\tht'=\mbox{const}
\ee
on the unit circle. Solving this equation together with \eqref{6.1}, Dyson finds that for $Q(\eit)$ sufficiently smooth,
\be
F_T=-\half \log n!+nQ_0-\sum_{k=-\infty}^\infty |k|Q_kQ_{-k}
\ee
where $\{Q_k\}$ are the Fourier coefficients of $Q(\eit)$. 
For $V=-\beta Q=-2 Q$, we obtain
\be
-\half\log(n!D_n(e^{V}))=F\sim F_T=-\half \log n!+nQ_0-\sum_{k=-\infty}^\infty |k|Q_kQ_{-k}
\ee
which is precisely SSLT, $D_n(e^{V})\sim\exp\{ nV_0+\sum_1^\infty k |V_k|^2\}$. 
There is an irony in  this argument: on
the one hand, SSLT was proved in response to a problem in statistical mechanics, but on the other hand,
SSLT is obtained here  by the methods of statistical mechanics.
The above calculation is based on the (unproven) validity of the thermodynamic recipe
\eqref{5.1}--\eqref{5.4}, \eqref{6.1}. Johansson's proof in \cite{Joh1} does not rely on the thermodynamic recipe. 
The proof in \cite{Joh1} starts with the identity $D_n(1)=1$ and develops
a perturbation theory around $\varphi(z)=1$. Johansson notices that the minimum of the double sum in \eqref{eq44},
and thus the maximum of the integrand in \eqref{eq43} for $V\equiv 0$ is reached 
when the points $e^{i\tht_j}$, $j=0,\dots,n-1$, are uniformly distributed on the unit circle, i.e. when they are the vertices of a regular $n$-gon. For $\varphi(z)=e^{V(z)}\not\equiv 1$, normalize $\varphi$ by requiring $V_0=0$.
 If one adds the corresponding first sum with $V\not\equiv 0$ to \eqref{eq44}, the position of the maximum of the integrand 
in \eqref{eq43} becomes displaced, but in a sense only slightly if $n$ is large. For large $n$ Johansson introduces the
following change of variables in \eqref{eq43}:
\be\label{chvar}
\tht_j=\psi_j-\frac{1}{n}h(e^{i\psi_j}),\qquad j=0,\dots, n-1,
\ee
where $h(z)=-i\sum_{k=1}^{\infty}(V_k z^k-V_{-k}z^{-k})$ is the conjugate function to $V(z)$.
The main contribution to the integral comes from a small neighborhood of a configuration where $\psi_j$'s are uniformly distributed on the unit circle.
Due to this crucial fact, all the sums which appear (as a result of the change of variables \eqref{chvar})
in the exponent, apart from the main one, $\sum_{0\le k<j\le n-1}W(\psi_k-\psi_j)$,  can be replaced by integrals with good precision and produce, as $n\to\infty$, the constant \eqref{eq33}, while the remaining $\frac{1}{n!}\int\exp(-\sum_{0\le k<j\le n-1}W(\psi_k-\psi_j))\prod_{k=0}^{n-1} \frac{d\psi_k}{2\pi}=D_n(1)=1$.

The fifth method is combinatorial and is due basically to Kac \cite{Kac}.  At the start one renormalizes
$\varphi\!\lf(\eit\rt)$ so that $\sup_\tht \,\varphi\!\lf(\eit\rt)=1$.  Setting $h=1-\varphi$ and $w_\lb
= 1-\lb h$, $0\le \lb \le 1$, we have $w_1= \varphi$.  Under the additional assumption that $\inf_\tht \,\varphi
\!\lf(\eit\rt) >0$, we clearly have $0\le h(\tht) \le \sup_\tht \, h(\tht) <1$.  Using the identity $\log \det A
= \tr \log A$, one expands $\log D_n \lf(1-\lb h\rt)$ for $0\le \lb \le 1$ in a power series
(note that $\|T_n(h)\| \le \|h\|_\infty <1$)
\begin{align}\label{eq45p}
\log D_n \lf(1-\lb h\rt) &= -\sum^\infty_{m=1} \frac{\lb^m}{m} \,\tr \lf( \lf[ T_n(h) \rt]^m\rt)\\
\intertext{and also}
\int \log \lf(1-\lb h\rt)\, \frac{d\tht}{2\pi} &= -\sum^\infty_{m=1} \frac{\lb^m}{m}
\int h^m(\tht)\, \frac{d\tht}{2\pi}\ .
\end{align}
On the other hand
\be\label{eq45n}
\lf(\log \lf(1-\lb h\rt)\rt)_k = - \int e^{-ik\tht} \sum^\infty_{m=1} \frac{\lb^m}{m}
\, h^m(\tht) \, \frac{d\tht}{2\pi}
\ee
and so
\be\label{eq45nn}
\sum^\infty_{k=1} k \lf| \lf(\log \lf(1-\lb h\rt)\rt)_k \rt|^2  = \sum^\infty_{m=1}
s_m \, \frac{\lb^m}{m}
\ee
where
$$
s_m = \sum^\infty_{k=1} k \sum^{m-1}_{\ell=1} \frac{m}{\ell (m-\ell)} \lf(\int e^{-ik\tht}
\, h^\ell (\tht) \, \frac{d\tht}{2\pi} \rt) \lf(\int e^{ik\tht} \, h^{m-\ell}(\tht)\,
\frac{d\tht}{2\pi}\rt).
$$
Thus the strategy is to show
\be\label{eq45pp}
\lim_{n\to\infty} - \lf[ \tr \lf(\lf[T_n (h)\rt]\rt)^m -n \int h^m(\tht) \,\frac{d\tht}{2\pi} \rt]=s_m.
\ee
Combinatorial issues arise when one tries to control the complicated multilinear sums in $\tr \lf[T_n(h)\rt]^m$.

The sixth method is based   on a formula of Baik, Deift, McLaughlin, and Zhou (unpublished) for relative
determinants
\be\label{eq45ppp}
\log \frac{D_n \lf(\varphi \psi\rt)}{D_n \lf(\psi\rt)} = \int^1_0 dt \int_0^{2\pi}
\lf[ \frac{d}{dt} \log \varphi_t\rt] K_t(\tht) \, \varphi_t(\tht) \, \frac{d\tht}{2\pi}
\ee
where $\varphi$ and  $\psi$ are two functions on $S^1$ and
$\varphi_t = (1-t)+t\varphi$ and $K_t(\tht)$ is the \tit{Christoffel-Darboux} kernel
\be\label{eq45-4}
K_t(\tht) = \sum^{n-1}_{j=0} \lf|p_j\!\lf(e^{i\tht}, \, \varphi_t \, \psi\rt)\rt|^2
\ee
and $p_j\!\lf(\eit,\, \varphi_t\, \psi\rt)$ are the OPUC's associated with the weight
$\varphi_t\!\lf(\eit\rt) \psi\!\lf(\eit\rt) \frac{d\tht}{2\pi}$.  SSLT then follows by
computing the (weak) limit of $\frac{1}{n} \, K_t(\tht)$ as $n\to \infty$. Thus, as in the
first method, SSLT is seen here as following from properties of OPUC's $p_j\!\lf(\eit\rt)$.

The mathematical challenge that is posed by SSLT can be viewed as follows.   In functional analysis
one commonly considers continuous maps $F$, say, from one fixed topological space $X$ to another
fixed topological space $Y$.  Then in order to verify that $\lim_{n\to\infty} F(x_n) = F(x)$ all
one has to do is verify that $x_n \to x$ as $n\to\infty$ in the topology of $X$.  In particular, if
$X$ is the space of trace class operators on a separable Hilbert space $\MCH$, and  $F(x)= \det
\lf(1+x\rt)$, then $F:X\to Y= \CC$ and $\lim_{n\to\infty} \,\det \lf(1+x_n\rt) = \det \lf(1+x\rt)$
if $x_n\to x$ in trace norm.  But SSLT is not of this standard type.  Instead, one has $F:M(n, \CC)
\to \CC$ taking $A_n \in M(n, \CC)$ to $F\!\lf(A_n\rt) = \det A_n$.  So the initial space $X=X_n
=M(n,\CC)$ is changing with $n$ and it is not at all clear how to topologize the situation so that
``$\lim_n F\!\lf(A_n\rt) = F\!\lf(\lim A_n\rt)$''.  If, for example, we let $X=M(\infty, \CC)$ be the
inductive limit of the $M(n, \CC)$'s, then $T_n(\varphi)$ converges to the Toeplitz operator $T(\varphi)$
in the associated topology, but $T(\varphi)$ is not of the form $1+ \text{ trace class}$, and so
$F\!\lf(T(\varphi)\rt) = \det\!\lf(T(\varphi)\rt)$ is not defined.  In such situations when $X= X_n$, or
$Y=Y_n$, is varying, there is no general procedure to follow and each problem must be addressed on an
ad hoc basis.  And when there is no preferred path to the solution of a problem, there are many paths,
as we see from the six
very different methods above.  In this regard, the method of Borodin-Okounkov-Geronimo-Case is
singled out: after factoring out $e^{n(\log\varphi)_0+E(\varphi)}$, one has
\be\label{eq45-5}
\frac{\det \lf(T_n (\varphi)\rt)}{e^{n(\log\varphi)_0+E(\varphi)}} = \det\lf(1-B_n\rt)
\ee
where $B_n=Q_n\, H(b)\, H(\widetilde{c})\, Q_n$ is trace class in $\ell^+_2$ and $B_n\to 0$ in trace norm
as $n\to\infty$.  In other words, one of the paths is to renormalize the problem so that it assumes
the standard form.

This pathway, to reduce SSLT to standard form, goes back to Harold Widom, who was the first to apply
operator-theoretic methods to the problem of Toeplitz asymptotics. In \cite{Wid4}, Widom gave yet another proof of 
\eqref{eq31}. This proof, which was found more than 35 years ago, has not lost its charm.
We follow the presentation in \cite{Bott1}.
It turns out that by a simple trick, the determinants of finite sections of the inverse Toeplitz operator 
$T(\varphi)^{-1}$ are easy to compute,
\be\label{eq46}
\det P_n \, T(\varphi)^{-1}\, P_n = e^{-n(\log\varphi)_0}
\ee
where $P_n$ projects $u= \lf(u_0, u_1, \dots, u_{n-1}, u_n, \dots \rt) \in \ell^+_2$ to $\lf(u_0,
u_1, \dots, u_{n-1}\rt) \in \CC^n$.  Simple algebra shows that
\be\label{eq46-1}
\begin{aligned}
\det T_n(\varphi) = & \det P_nT(\varphi)P_n=\det\!\lf( P_n \, T^{-1} \lf(\varphi^{-1}\rt) P_n\rt)\\
& \times \det\!\lf(1+ \lf(P_n \, T^{-1} \lf(\varphi^{-1}\rt) P_n\rt)^{-1} \, P_n\,K\,P_n\rt)
\end{aligned}
\ee
where $K\equiv T(\varphi)- T^{-1} \lf(\varphi^{-1}\rt)$ is trace-class.  By \eqref{eq46},
\be\label{eq46-2}
\det P_n \, T^{-1} \lf(\varphi^{-1}\rt)P_n = e^{n(\log\varphi)_0}.
\ee
But $\lf(P_n \, T^{-1} \lf(\varphi^{-1}\rt)P_n \rt)^{-1} \, P_n $ can be shown to converge strongly
as $n\to\infty$ in $\ell^+_2$ to $\lf(T^{-1} \lf(\varphi^{-1}\rt)\rt)^{-1} = T\lf(\varphi^{-1}\rt)$,
and as $K$ is trace class, $R_n \equiv \lf(P_n\, T^{-1} \lf(\varphi^{-1}\rt)P_n\rt)^{-1} \, P_n \,
K\,P_n \to T\lf(\varphi^{-1}\rt)K$.  We conclude that as $n\to\infty$
\be\label{eq46-3}
\frac{\det T_n(\varphi)}{e^{n(\log\varphi)_0}} = \det\!\lf(1+R_n\rt) \to \det\!\lf(1 + T\lf(\varphi^{-1}\rt)K\rt)
= \det\,T\lf(\varphi^{-1}\rt) T(\varphi).
\ee
Using the remarkable Helton-Howe-Pincus formula (see \cite{HeltHow}),  
$\det\!\lf(e^A\, e^B\, e^{-A}\, e^{-B}\rt)=e^{\mathrm{tr} [A,B]}$,
if $A$, $B$ are bounded, and
$[A,B]= AB-BA$ is trace class, one then shows that $\det \!\lf(T\lf(\varphi^{-1}\rt)T(\varphi)\rt) =e^{E(\varphi)}$,
and \eqref{eq31} is proved.

Of all the methods to prove SSLT, perhaps the simplest and most direct was introduced by Basor and Helton
in \cite{BasHelt}. As in \cite{Wid4}, the method also proceeds by reducing SSLT to standard form. We follow 
\cite{BasHelt}: Suppose $\varphi$ has a Wiener-Hopf factorization $\varphi=\varphi_+\varphi_-$, where 
$\varphi_+$, $\varphi_+^{-1}$ are bounded analytic functions in $\{|z|<1\}$,
and $\varphi_-$, $\varphi_-^{-1}$ are bounded analytic functions in $\{|z|>1\}$. Then by direct calculation, the Toeplitz
operator $T(\varphi)$ has a factorization $T(\varphi)=T(\varphi_-)T(\varphi_+)$, furthermore,
\be
P_nT(\varphi_+)=P_nT(\varphi_+)P_n,\qquad T(\varphi_-)P_n=P_nT(\varphi_-)P_n.
\ee
As $T(\varphi_\pm)T(\varphi_\pm^{-1})=1$, we see that $T(\varphi_\pm)^{-1}$ exist and
$T(\varphi_\pm)^{-1}=T(\varphi_\pm^{-1})$. We thus have
\be
P_nT(\varphi)P_n=P_nT(\varphi_+)P_nT(\varphi_+)^{-1}T(\varphi_-)T(\varphi_+)T(\varphi_-)^{-1}P_nT(\varphi_-)P_n
\ee
so that
\be\begin{aligned}
D_n(\varphi)=&\det\!\lf(P_nT(\varphi_+)P_n\rt)\det\!\lf(P_n\left\{T(\varphi_+)^{-1},T(\varphi_-)\right\}P_n\rt)
\det\!\lf(P_nT(\varphi_-)P_n\rt)\\
=& \lf((\varphi_+)_0(\varphi_-)_0\rt)^n\det\!\lf(P_n\left\{T(\varphi_+)^{-1},T(\varphi_-)\right\}P_n\rt)
\end{aligned}
\ee
where $\left\{T(\varphi_+)^{-1},T(\varphi_-)\right\}= T(\varphi_+)^{-1}T(\varphi_-)T(\varphi_+)T(\varphi_-)^{-1}$ is the 
multiplicative commutator.
However $\left\{T(\varphi_+)^{-1},T(\varphi_-)\right\}=1+T(\varphi_+)^{-1}\lf[T(\varphi_-),T(\varphi_+)\rt]T(\varphi_-)^{-1}$ 
and under mild
smoothness conditions on $\varphi(e^{i\tht})$, the commutator
$\lf[T(\varphi_-),T(\varphi_+)\rt]=T(\varphi_-)T(\varphi_+)-T(\varphi_+)T(\varphi_-)$ is trace class. For such $\varphi$, 
$\det\!\lf(P_n\left\{T(\varphi_+)^{-1},T(\varphi_-)\right\}P_n\rt)$ is equal to the Fredholm determinant
\[
\det\!\lf(1+P_nT(\varphi_+)^{-1}\lf[T(\varphi_-),T(\varphi_+)\rt]T(\varphi_-)^{-1}P_n\rt)
\] 
which converges by abstract theory as $n\to\infty$ to
\be\begin{aligned}
\det\!\lf(1+T(\varphi_+)^{-1}\lf[T(\varphi_-),T(\varphi_+)\rt]T(\varphi_-)^{-1}\rt)=&
\det\!\left\{T(\varphi_+)^{-1},T(\varphi_-)\right\}\\
=&\det\!\lf(T(\varphi_+)^{-1}T(\varphi_-)T(\varphi_+)T(\varphi_-)^{-1}\rt)\\
=&\det\!T(\varphi^{-1})T(\varphi).
\end{aligned}
\ee
Again by the Helton-Howe-Pincus formula, we have
\be
\lim_{n\to\infty} \frac{D_n (\varphi)}{e^{n(\log\varphi)_0}} =\det\!T(\varphi^{-1})T(\varphi)=e^{E(\varphi)}
\ee
which proves SSLT under mild assumptions on $\varphi$ (see \cite{BasHelt} Theorem 2.1.).

In \cite{Dei1}, the author gives yet another proof of \eqref{eq31} by
expressing $D_n(\varphi)$ as a Fredholm determinant of the following form:
\be\label{D1}
D_n(\varphi) = \det \lf(1-V_n\rt)
\ee
where $V_n$ is a trace class operator acting on $L^2(S^1)$
\begin{align}\label{D2}
V_n \, f(z) &= \int_{|z'|=1} V_n \lf(z, z'\rt) f\lf(z'\rt) dz' \ , \qquad |z|=1,
\\
\intertext{with}\label{D3}
V_n \lf(z, z'\rt) &= \frac{z^n \lf(z'\rt)^{-n} - 1}{z-z'} \; \frac{1-\varphi \lf(z'\rt)}{2\pi i}\ .
\end{align}
The operator $V_n$ is oscillatory and does not converge as $n\to\infty$ in the trace-norm topology.  However,
$V_n$ has the special form of an \tit{integrable operator}, that is, operators with kernels of the form
$\sum_{j=1}^m\frac{f_j(z)g_j(z')}{z-z'}$, with $z$, $z'$ on some contour $\Gamma$ in $\CC$, $m<\infty$, and 
$f_j$, $g_j$ are given functions on $\Gamma$.
Such operators were first singled out as a distinguished class by Its, Izergin, Korepin and Slavnov\footnote{
Some of the elements of the theory of integrable operators were already present in the earlier work \cite{Sakh}.} (see \cite{Dei1} for a pedagogic presentation).
Associated with the integrable operators there is a canonical Riemann-Hilbert Problem, which, in this case, is of
oscillatory type.  Such problems can be analyzed asymptotically by applying the nonlinear
steepest descent method for Riemann-Hilbert Problems introduced by Deift and Zhou in \cite{DZ},
and further developed in \cite{DVZ}, and SSLT then follows (see \cite{Dei1}).

\begin{remark}\label{22+++}
In the methods of Widom and Basor-Helton the existence of the limit of $\frac{D_n(\varphi)}{e^{n(\log\varphi)_0}}$ as 
$n\to\infty$ follows directly from general theorems in functional analysis: the identification of the limit
with $e^{E(\varphi)}$ is a separate matter and uses other means (Helton-Howe-Pincus).
The situation is similar regarding
Szeg\H o's original proof of SSLT (see \cite{Sz3} \cite{GreSz}): as $\varphi(e^{i\tht})>0$,
the ratio $\frac{D_n(\varphi)}{e^{n(\log\varphi)_0}}$
is monotone increasing in $n$, so the existence of the limit (which may, a priori, be infinite) is immediate. 
Again the evaluation of the limit is a separate matter. Indeed, Szeg\H o first computes the limit for a dense class of
symbols $\{\widetilde\varphi\}$ for which $D_n(\widetilde\varphi)$ can be evaluated explicitly for $n$ sufficiently large, 
$n>n_\varphi$, and then obtains the general result by an approximation argument $\widetilde\varphi\to\varphi$.

By contrast, all the other methods we have described above prove the existence of the limit, and compute its value,
both at the same time.
\end{remark}

We return now briefly to the 2-dimensional Ising model.  One notes that Szeg\H{o}'s Theorem \ref{theorem3}
does not actually apply to $\varphi_{\tdg}(\tht)$, say, in \eqref{eq28}, because the symbol
is not real valued.  Moreover, Szeg\H{o}'s proof of \eqref{eq31} depends in a crucial way on the
positivity of $\varphi\!\lf(\eit\rt)$ (see Remark \ref{22+++} above), and it was not at all clear at the time how to extend Theorem
\ref{theorem3} to non-real symbols $\varphi$.  So Onsager's lament that ``\dots the mathematicians got
there first'', was in fact not true!  The result in the 1950 unpublished draft of Kaufman and Onsager mentioned above,
without ``the epsilons and deltas and all of that'', but allowing for non-real $\varphi$, was in fact
the first to break the tape!  Now although the symbol $\varphi_{\tdg} (\tht)$ is not real,
for $0 <k_{\tons} <1$ however, $\log\varphi_{\tdg} (\tht)$ is still smooth on $S^1$, in particular,
the winding number of $\varphi_{\tdg} (\tht)$ on $S^1$ is zero. The first result for such
symbols was proved by G.~Baxter in 1963 \cite{BaxG2}.  He showed that under certain technical conditions,
if $\varphi\!\lf(\eit\rt)$ is smooth, non-zero, and has zero winding number on $S^1$, then
\be\label{eq47}
\lim_{n\to\infty} \frac{D_n (\varphi)}{e^{n(\log\varphi)_0}} =e^{\sum^\infty_{k=1}\, k \lf(\log \varphi\rt)_k
\lf(\log \varphi\rt)_{-k}}.
\ee
Note that for $\varphi_{\tdg}(\tht)$ in \eqref{eq28}, one has for $0<k_{\tons} <1$,
$\lf(\log \varphi_{\tdg}\rt)_\ell =  (\mbox{sgn } \ell) \, k^{|\ell|}_{\tons}/ 2 |\ell|$ for
$\ell \neq 0$ and $\lf(\log \varphi_{\tdg}\rt)_0=0$.  Hence $\sum^\infty_{\ell=1} \ell
\lf(\log \varphi\rt)_{\ell} \lf(\log \varphi\rt)_{-\ell} = -\frac{1}{4}
\sum^\infty_1 \frac{\lf(k^2_{\tons} \rt)^\ell}{\ell} = \frac{1}{4} \;\log \lf(1-k^2_{\tons}\rt)$, which
yields \eqref{eq15} via \eqref{eq29}.

In retrospect, many of the proofs of SSLT do extend directly to complex-valued symbols with smooth logarithm.  This is
true, in particular, for Kac's proof of \eqref{eq31} in 1954 \cite{Kac}, as noted by Montroll, Potts
and Ward in \cite{MPW} (see footnote 17, p.~316). Moreover, the proofs in \cite{Wid4} and 
in the text \cite{BottSilb3}, 
for example, are all designed from the very beginning to cover complex-valued functions. 

If $\varphi(\eit)$ is non-zero and continuous, say, on $S^1$, and has zero winding
number, then $\varphi(\eit) = e^{V(\eit)}$ for some continuous, periodic function $V(\eit)
=\log\varphi(\eit)$.  The strongest version of SSLT with complex symbols $\varphi(\eit)$ is the following.

\begin{theorem}\label{theorem7}
Let $V(\eit)\in L^1(S^1)$ be a (possibly complex-valued) function on $S^1$ with Fourier coefficients
$\lf(V_k\rt)_{k\in \ZZ}$ satisfying
\be\label{eq48}
\sum^\infty_{k=-\infty} |k| \; \lf|V_k\rt|^2 < \infty \ ,
\ee
then
\be\label{eq49}
\lim_{n\to\infty} \frac{D_n\!\lf(e^V\rt)}{e^{nV_0}}=
e^{\sum^\infty_1 k\, V_k\, V_{-k}}.
\ee
\hfill $\Box$
\end{theorem}
Note that for $V(\eit)$ satisfying \eqref{eq48}, $e^{V(\eit)} \in L^p\!\lf(S^1, \frac{d\tht}{2\pi}\rt)$
for $1 \le p < \infty$ (see \cite{Sim2}).
The proof of Theorem \ref{theorem7} is given by Johansson \cite{Joh1}. In particular, his argument for the extension of the result to complex-valued $V(z)$ proceeds by considering the
analytic functions
$$
g_n (z) = D_n\!\lf(e^{\lf(\Re V + z \Im V\rt)}\rt)\bigl/ e^{n \lf[(\Re V)_0 + z (\Im V)_0\rt]}
$$
and applying Vitali's Theorem.  When $z$ is real, $\lim_{n\to\infty} g_n(z)$ is given by \eqref{eq47},
and to obtain \eqref{eq49} one just sets $z=i$.

Under the additional assumption that $V$ is in $L^\infty(S^1)$, Theorem \ref{theorem7}
was already established in \cite{Wid4}. We also note that condition \eqref{eq48} can be replaced by certain 
``non-symmetric'' requirements on the decay of the Fourier coefficients $V_k$ and $V_{-k}$, $k\ge  0$; see, 
for example, Theorem 10.32 and Corollary 10.42 of \cite{BottSilb3}.

\begin{remark}\label{remark3}
In \cite{Joh1}, Johansson observes the following very interesting consequence of Theorem 
\ref{theorem7} (see also \cite{Joh2}).
Let $V(\eit)$ be a real-valued function on $S^1$ with Fourier coefficients $\{V_k\}$ satisfying \eqref{eq48}
and suppose $V_0= \int V(\eit) \;\frac{d\tht}{2\pi} =0$.  As in the Bump-Diaconis method for SSLT
described above, equip the unitary group $\MCU(n)$ with Haar measure $dg$ and consider the random variable
$F_n(g)=\sum^{n-1}_{j=0} V\lf(e^{i\tht_j (g)}\rt)=\tr V(g)$, where  $\lf\{e^{i\tht_j (g)}\rt\}$ are the eigenvalues of
$g \in \MCU(n)$.  Then we have the following Central Limit Theorem:
\be\label{eq51}
F_n  \text{ converges in distribution to } N\!\lf(0, \sg^2\rt),\qquad \sg= \lf(2 \sum^\infty_{k=1}
\, k\lf|V_k\rt|^2\rt)^{\frac{1}{2}}
\ee
where $N\!\lf(0, \sg^2\rt)$ denotes the normal distribution with mean zero and variance $\sg$. Indeed, for any real $t$ we have
\be
\int_{\MCU(n)} e^{i\,t\,F_n(g)}\, dg = D_n \!\lf(e^{i\,t\,V}\rt)
\ee
and so by Theorem \ref{theorem7}, as $n\to\infty$
\be
\int_{\MCU(n)} e^{i\,t\,F_n(g)} \, dg \to e^{\sum^\infty_1\, k (itV)_k \: (itV)_{-k}}= 
e^{-\frac{t^2}{2} \lf[2\lf(\sum^\infty_1 k|V_k|^2\rt)\rt]}
\ee
from which we conclude \eqref{eq51}.   One notes that as opposed to the standard Central Limit Theorem for
independent random variables, the $\tht_j$'s are very strongly uniformly distributed and no $n^{-1/2}$
scaling is needed. A suitable scaling is needed in the more general ``interpolating'' situation considered by 
Duits and Johansson \cite{DJ} where $V(z)$ depends on $n$ (cf. Section \ref{sec-ds}).
In \cite{DJ}, for each $n\ge 1$, the authors consider sequences $\{ k_j(n)\}_{j\ge 1}$ of mutually distinct positive integers,
and they show that for each $m\ge 1$,  the random variables
$\frac{1}{\sqrt{\min(k_j(n),n)}}\tr g^{k_j(n)}$, $j=1,\dots, m$, converge to independent standard complex normals.
(This connects earlier results of Diaconis and Shahshahani for fixed $k_j$, and of Rains for $k_j(n)>n$: 
see \cite{DJ} for references.) Note an interesting difference in scaling for $k_j(n)\le n$ and for $k_j(n)>n$. In the 
latter case, the scaling is $n^{-1/2}$. 
The authors in \cite{DJ} obtained their results by proving an analogue of SSLT for symbols depending on $n$ in a special way,
as follows. If
\[
V(n;z)=\sum_{|j|>0}\frac{\gamma_j z^{k_j(n)}}{\sqrt{\min(|k_j(n)|,n)}},
\]
where $\sum_j|\gamma_j|<\infty$, the sequences $\{k_j(n)\}_{j\ge 1}$ are as above, and $k_{-j}(n)=-k_j(n)$, then
\[
\lim_{n\to\infty}D_n(e^{V(n;z)})=e^{\sum_{j=1}^{\infty}\gamma_j\gamma_{-j}}.
\]

\end{remark}

\section{Back to Ising}\label{section4}
After Onsager announced formula \eqref{eq15} for $M_0$ in 1948, and again in 1949, but did not
provide a proof, three years went by until a derivation of the formula was found by Yang in a tour
de force in 1952 \cite{Yan1}.  Yang's method is based on results of Kaufman and Onsager, but his approach
is different, and he does not explicitly use Toeplitz determinants. In his Selected Papers in 2005 
(\cite{Yan2}), Yang described the long and winding road that led him to the Kaufman-Onsager formula \eqref{eq15}. The story has much charm: 

``I was thus led to a long calculation, the longest in my career. Full of local, tactical tricks, the calculation proceeded by twists and turns. There were many obstructions. But always, after a few days, a new trick was somehow found that pointed to a new path. The trouble was that I soon felt I was in a maze and was not sure whether in fact, after so many turns, I was anywhere nearer the goal than when I began. This kind of strategic overview was very depressing, and several times I almost gave up. But each time something drew me back, usually a new tactical trick that brightened the scene, even though only locally.

Finally, after about six months of work off and on, all the pieces suddenly fitted together, producing miraculous cancellation, and I was staring at the amazingly simple final result.'' 

\begin{remark}\label{remark2}
Let $E_\La(\sg)$ denote the interaction energy in a rectangular box $\La$ of size $M\times N$ as
in \eqref{eq6} above, and let $E_{\La, h}(\sg)= E_\La(\sg) - h \sum_{(i,j)\in \La} \sg_{i,j}$ denote the
energy in the presence of a magnetic field of strength $h$.  Let $Z_{\La, h} = \sum_\sg e^{-E_{\La, h}
(\sg)/k_B\, T}$ denote the associated partition function, and let  $F_{\La, h} = -\frac{k_B \, T}{MN} \log
Z_{\La, h}$ denote the associated free energy in the box $\La$.  What Yang proved in \cite{Yan1} was
that, for $T<T_c$,
\be
\begin{aligned}
M_Y &\equiv -\lim_{h\downarrow 0}\; \lim_{|\La| \to \infty} \frac{F_{\La, h/N} - F_{\La, 0}}{h/N}\\
&= \lf(1-k^2_{\tons}\rt)^{\frac{1}{8}} = M_0 \equiv \lim_{n\to \infty} \lf\langle \sg_{1,1} \;
\sg_{1,n}\rt\rangle^{\frac{1}{2}}.
\end{aligned}
\ee
(In \cite{Yan1}, Yang only considered the case $J_1=J_2$.  Yang's method was extended to the general case
$J_1 \neq J_2$ by Chang \cite{Cha}.)
\end{remark}

As noted and emphasized by Schultz, Mattis and Lieb \cite{ScML} in 1964, neither $M_Y$ nor $M_0$ coincide,
a priori,  with the physically fundamental definition of the spontaneous magnetization $M_{\tph}$, viz.,
\be
\begin{aligned}
M_{\tph}  &\equiv -\lim_{h\downarrow 0} \; \lim_{|\La| \to \infty}
\frac{\p \, F_{\La, h}}{\p h}\\
&= \lim_{h \downarrow 0} \; \lim_{|\La|\to\infty} \frac{1}{MN}
\lf[ \frac{\sum_\sg \lf(\sum_{(i,j)\in \La} \sg_{i,j}\rt) e^{-\frac{1}{k_B\,T} \, E_{\La, h} (\sg)}}{
\sum_\sg e^{-\frac{1}{k_B\,T} \, E_{\La,h}(\sg)}} \rt].
\end{aligned}
\ee
Only some years later was it shown that indeed for $T< T_c$
\be
M_{\tph} = M_Y = M_0 = \lf(1-k^2_{\tons}\rt)^{\frac{1}{8}}
\ee
(see, in particular, \cite{LebM-Lof}, \cite{BGJ-LS}).

In 1955 an important technical advance was made by Potts and Ward \cite{PotWar} when they showed, in
particular, that the correlation function $\lf\langle \sg_{1,1}\, \sg_{1,n}\rt\rangle$ along a row for the
2-D Ising model could be expressed in terms of a \tit{single} Toeplitz determinant, rather than a sum
of two Toeplitz determinants as in \eqref{eq17}.  They were unable  at the time to verify that their
solution agreed with
\eqref{eq17}: This was done later in 1963 \cite{MPW} by Montroll, Potts and Ward, who also used recent ideas of
Kasteleyn on the Pfaffian approach to the dimer and Ising problems to give a new and simpler proof of the formula
in \cite{PotWar}.  The formula in \cite{PotWar} \cite{MPW} is the following
\be\label{eq50}
\lf\langle \sg_{1,1} \, \sg_{1, 1+n} \rt\rangle =  D_n\!\lf(\varphi_{\tons}\rt)
\ee
where $\varphi_{\tons}\!\lf(\eit\rt)$ is the Onsager function defined in \eqref{eq21}.  Through discussions with
Mark Kac, Montroll et al.\ were aware of SSLT. Furthermore they knew, as noted above, how to extend Szeg\H{o}'s
result via Kac's method to the  case where the symbol is complex.  For $T< T_c$, Montroll, Potts and Ward noted
that $0< \gamma_1, \gamma_2 <1$ by \eqref{eq26}, and so $\varphi_{\tons}^{\pm 1} \!\lf(\eit\rt) \neq 0$ 
and $\varphi_{\tons}$
has no winding on $S^1$.  As in the computation for $\varphi_{\tdg}\!\lf(\eit\rt)$ above, we find
\begin{align}
\lf(\log \varphi_{\tons}\rt)_\ell &= -\lf(\frac{\gamma_1^\ell}{2\ell} - \frac{\gamma_2^\ell}{2\ell}\rt), \qquad
\ell \ge 1\\
\lf(\log \varphi_{\tons}\rt)_{-\ell} &= \lf(\frac{\gamma_1^\ell}{2\ell} - \frac{\gamma_2^\ell}{2\ell}\rt), \qquad
\ell \ge 1
\end{align}
and $\lf(\log \varphi_{\tons}\rt)_0 =0$.
This implies $\sum^\infty_{\ell=1} \ell \lf(\log\varphi_{\tons}\rt)_\ell \, \lf(\log \varphi_{\tons}\rt)_{-\ell}
=\frac{1}{4} \, \log \frac{\lf(1-\gamma_1^2\rt)\lf(1-\gamma_2^2\rt)}{\lf(1-\gamma_1 \gamma_2\rt)^2}$ and after some
elementary algebra using \eqref{eq22}, we again obtain \eqref{eq15}, but now via \eqref{eq16} (cf. the discussion 
following \eqref{eq47}).

We now consider the spontaneous magnetization for the Ising model for temperatures $T \ge T_c$.  As
$k_{\tons} \uparrow 1$ when $T\uparrow T_c$, one anticipates from \eqref{eq15} that at $T=T_c$
the spontaneous magnetization $M_0$ is zero.  In \cite{KauOns}, Kaufman and Onsager proceed in the following
way.  They derive an approximate formula for $b_k$ in \eqref{eq19} ($b_k = \sum_{-k}$ in the notation
of \cite{KauOns}) at $T=T_c$
\be\label{eq52}
b_k \sim \frac{2}{\pi\,(2k+1)} \qquad\text{ as }\qquad k\to \infty.
\ee
They then take \eqref{eq52} as an equality for all $k\ge 0$, and substitute this relation into \eqref{eq17}
to obtain an approximate asymptotic formula for $\lf\langle \sg_{1,1}\, \sg_{1, 1+n}\rt\rangle$ in the form
of a linear combination $c^{*2}_2\, W_n - s^*_2\, \widetilde{W}_n$ of two Cauchy determinants $W_n$ and
$\widetilde{W}_n =W_{-n}$.  Such determinants can be evaluated explicitly in terms of gamma functions\footnote{
Formulae \eqref{eq53} \eqref{eq54} differ slightly from the corresponding formulae in \cite{KauOns}.}
\begin{align}
\lf|W_n\rt| &= \frac{2}{\pi} \prod^{n-1}_{q=1} \frac{\Gamma(q+1)\, \Gamma(q+1)}{\Gamma \lf(q+\frac{1}{2}\rt)
\Gamma \lf(q+\frac{3}{2}\rt)}, \label{eq53}\\
\lf|\widetilde{W}_n\rt| &= \frac{2}{3\pi} \prod^{n-1}_{q=1} \frac{\Gamma (q+1) \, \Gamma(q+1)}{
\Gamma \lf(q-\frac{1}{2}\rt)\Gamma\lf(q+\frac{5}{2}\rt)}\label{eq54}
\end{align}
and Kaufman and Onsager observe that both $W_n$ and $\widetilde{W}_n\to 0$ as $n\to\infty$, with
$\widetilde{W}_n$ decaying at a much faster rate than $W_n$.  They conclude that at $T=T_c$
\be\label{eq55}
\lf\langle \sg_{1,1}\, \sg_{1, 1+n}\rt\rangle \sim (-1)^n c^{*2}_2 W_n \to 0 \qquad \text{ as }\quad n\to \infty
\ee
and so $M_0=0$.  In 1959, Fisher \cite{Fis1} picks up on these calculations, and using Stirling's formula
in \eqref{eq53} he shows via \eqref{eq55} that as $n\to\infty$
\be\label{eq56}
\lf\langle \sg_{1,1} \, \sg_{1, 1+n}\rt\rangle \sim C\: n^{-1/4}, \qquad T=T_c
\ee
for some (undetermined) constant $C$.  This is the first derivation of an explicit decay rate for the
correlation function at $T_c$.   Fisher argues further that off the lattice axes, \eqref{eq56} may be
generalized (in the case $J_1=J_2$) to
\be\label{eq57}
\lf\langle \sg_{1,1}\, \sg_{1+ \ell, 1+m}\rt\rangle \sim C(\tht) \, k^{-1/4}, \qquad T=T_c
\ee
where $\tan\tht = \ell/m$ and $k=\lf(\ell^2 +m^2\rt)^{\frac{1}{2}} \to \infty$, and $C(\tht)$ depends
mildly on $\tht$ (see more below).  As Fisher learned later from Onsager himself (see reference 8
in \cite{Fis1}; also, cf.\  the reference preceding  \eqref{eq27} above
to  \cite{Dom2}), the analogue of formula \eqref{eq55} is in fact \tit{exact}
for correlations along the diagonal:
\be\label{Opure}
\lf\langle \sg_{1,1} \,\sg_{1+n, 1+n}\rt\rangle=(-1)^n W_n,\qquad T=T_c.
\ee 
Indeed from
\eqref{eq27} , \eqref{eq28}, for $k_{\tons}=1$ at $T_c$,
\begin{align}
\varphi_{\tdg}(\tht) &= \frac{1-e^{-i\tht}}{\lf|1-e^{i\tht}\rt|} = e^{i(\pi-\tht)/2}\ ,&\qquad &0<\tht <2\pi
\label{eq58}\\
\intertext{and so }
\lf(\varphi_{\tdg}\rt)_k &= \frac{2}{\pi} \; \frac{1}{2k+1}\ , &\qquad &k\in \ZZ
\label{eq59}
\end{align}
which agrees with \eqref{eq52} (note from \eqref{eq17} that $W_n$ is a Toeplitz determinant with shifted
entries $b_{k-1}$, and to observe \eqref{Opure} one can write $W_n$ as a determinant of the transposed matrix).   
If one repeats Fisher's calculation and keeps track of the constant,
one obtains the exact asymptotics at $T=T_c$,
\be\label{eq60}
\lf\langle \sg_{1,1}\, \sg_{1+n, 1+n}\rt\rangle = e^{\frac{1}{4}} \, A^{-3}\, 2^{\frac{1}{12}}\, n^{-1/4}
\;\lf(1-\frac{1}{64}\; \frac{1}{n^2} + \dots \rt) \quad\text{as}\quad n\to\infty
\ee
where 
\be\label{Glai}
A=e^{\frac{1}{12}}\, e^{-\zeta'(-1)}
\ee
is  Glaisher's constant, $\zeta=\text{Riemann's zeta function}$.
Formula \eqref{eq60} was first obtained by T. T. Wu in \cite{Wu} in serendipitous circumstances,
as we describe below. The full expansion of \eqref{eq60} to all orders was first given in \cite{AuYPer1}, together with a similar result, in the symmetric case $J_1=J_2$, for the next-to-diagonal correlation functions.

For $T>T_c$, Fisher notes in \cite{Fis1} that one can use the matrix viewpoint in \cite{Ons1} to show that
as $k=\sqrt{\ell^2 + m^2} \to \infty$
\be\label{eq61}
\lf\langle \sg_{1,1}\, \sg_{1+\ell, 1+m}\rt\rangle \sim B \,k^{-\frac{1}{4}} \; e^{-f_k(T)}
\ee
for some $\tht$-dependent constant $B=B(\tht)$, where to leading order as $k\to\infty$
\be\label{eq62}
f_k(T) \sim k\, g(T).
\ee
As $T\downarrow T_c$, $g(T) \downarrow 0$.
These calculations show, in particular, that $M_0=0$ for $T>T_c$.  If one includes logarithmic correction
terms in \eqref{eq62}, one obtains (see \cite{Fis2}) for $T>T_c$
\be\label{eq63}
\lf\langle \sg_{1,1}\, \sg_{1+\ell, 1+m}\rt\rangle \sim C(\tht, T) \;\frac{e^{-k\,g(T)}}{k^{1/2}}\;
\lf(1+O(1)\rt)
\ee
as $k\to\infty$.  This is the first derivation of an explicit decay rate for the correlation function at
$T>T_c$.  In \cite{KauOns} the authors show how the conclusion $M_0=0$ also follows from their determinantal
formulae, at least for $T>T_c$ sufficiently large.
In \cite{Kad}, Kadanoff obtains a decay rate for $T>T_c$ of the same form as in \eqref{eq63}, but only
in the double scaling limit $T-T_c \downarrow 0$, $k\to\infty$ such that $\lf|T-T_c\rt| \, k=c =
\text{ constant}$.  In \cite{Rya}, Ryazanov obtains the $k^{-1/4}$ decay rate at $T=T_c$, and he also
obtains the general form \eqref{eq63}, at least in the case that $T$ is large, $T>>T_c$.  The double
scaling limit, $T\to T_c$ and $k\to\infty$, for the Ising model is an important issue, to which we will
return at various points in this paper.  The importance of this issue rests in the fact that many
believe in universality  for spin models in the sense that the double scaling limit for the spin-spin
correlation function in the Ising model also gives precisely the same limit (modulo lattice-dependent
factors) for a wide class of ``non-integrable'' two-dimensional models with short range interactions.

\begin{remark}\label{remark4}
The result that $M_{\tph}$, the physically fundamental spontaneous magnetization defined above, is
in fact zero for all $T>T_c$, was proven rigorously only many years later by Lebowitz in \cite{Leb} in
(1972).  The proof in \cite{Leb} uses a variety of  techniques from statistical mechanics, but as shown
in \cite{LebM-Lof} \cite{BGJ-LS}, $M_{\tph} = M_0=\lim_{n\to\infty} \lf\langle \sg_{1,1}\, \sg_{1,1+n}
\rt\rangle^{\frac{1}{2}}$ for all $T>0$, and so the proof that
$M_{\tph}=0$ for $T>T_c$ (in fact, for $T\ge T_c$) follows, a posteriori,
from a proof that $\lim_{n\to \infty} D_n\!\lf(\varphi_{\tons}\rt)=0$, $T\ge T_c$.
\end{remark}

\begin{remark}\label{4+}
In the isotropic case, $J_1=J_2$, explicit expressions for the correlation functions 
$\lf\langle\sg_{1,1}\,\sg_{1+m, 1+\ell}\rt\rangle$ for values of $(m,\ell)$ other than $(0,\ell)$, $(m,0)$, or $(m,m)$,
were derived by Shrock and Ghosh \cite{ShGh}.  In particular, they obtained expressions
in the cases $(m,\ell)=(2,1)$, $(3,1)$, $(3,2)$, $(4,1)$, $(4,2)$, and $(4,3)$, in terms of complete elliptic integrals
$K$ and $E$. They also inferred a general structural formula for arbitrary $\lf\langle\sg_{1,1}\,\sg_{1+m, 1+\ell}\rt\rangle$
in terms of these elliptic integrals $K$ and $E$.

In further developments, Au-Yang, Jin, and Perk \cite{AuYPer2} \cite{AuYJPer} showed that the next-to-diagonal
correlations $\lf\langle\sg_{1,1}\,\sg_{1+n, n}\rt\rangle$ have a ``bordered'' Toeplitz determinant form, that is,
an $n\times n$ Toeplitz determinant with the last column replaced by a certain explicit vector 
$(b_{n-1},b_{n-2},\dots,b_0)^T$. This representation is then used to express   
$\lf\langle\sg_{1,1}\,\sg_{1+n, n}\rt\rangle$ in terms of a solution to the Painlev\' e VI equation by Witte in \cite{Wi}.
\end{remark}

\section{New questions for SSLT}
In view of the above calculations for $T\ge T_c$, the Ising problem raised new questions and challenges
in Toeplitz theory.  From \eqref{eq50}, $\lf\langle \sg_{1,1}\, \sg_{1, 1+n}\rt\rangle = D_n
\!\lf(\varphi_{\tons}\rt)$
(similarly, from \eqref{eq27}, $\langle \sg_{1,1}\, \sg_{1+n,\, 1+n}\rangle = D_n \lf(\varphi_{\tdg}\rt)$),
and so we want to know precisely how $D_n \!\lf(\varphi_{\tons}\rt)$ behaves as $n\to\infty$
for $T\ge T_c$.  But for $T>T_c$, we have $0< \gamma_1 <1 < \gamma_2$ (see \eqref{eq23}, \eqref{eq26}), and so in
contrast to the case $T<T_c$, the winding number of $\varphi_{\tons}$ is $-1 \neq 0$.  And at $T=T_c$, we
see from \eqref{eq25} that $\varphi_{\tons} \!\lf(\eit\rt)$ has a jump discontinuity  at $\tht=0$.
For such discontinuities the Fourier coefficients $\lf(\log \varphi_{\tons}\rt)_k$ decay as $k^{-1}$
and so $\sum^\infty_{k=1} k\lf|\lf(\log \varphi_{\tons}\rt)_k\rt|^2 =\infty$.  In both cases SSLT, or
more properly Theorem \ref{theorem7}, does not apply.

In 1966, in a tour de force in \cite{Wu}, Wu gave the first derivation of the precise asymptotics of
$\lf\langle \sg_{1,1},\, \sg_{1, 1+n} \rt\rangle = D_n\!\lf(\varphi_{\tons}\rt)$ as $n\to\infty$, for
$T\ge T_c$.  For $T>T_c$, Wu showed that for $n\to\infty$
\be\label{eq64}
\begin{aligned}
\lf\langle\sg_{1,1}\,\sg_{1, 1+n}\rt\rangle &= (\pi n)^{-\frac{1}{2}} \; \gamma_2^{-n} \lf(1-\gamma^2_1\rt)^{\frac{1}{4}}
\lf(1-\gamma^{-2}_2\rt)^{-\frac{1}{4}} \lf(1-\gamma_1 \gamma_2\rt)^{-\frac{1}{2}}\\
& \times \lf(1+\frac{A_1}{n} + \frac{A_2}{n^2} + \frac{A_3}{n^3} + \dots \rt)
\end{aligned}
\ee
with explicit expressions for $A_1$,  $A_2$ and $A_3$ in terms of $\gamma_1=z_1\,z_2^* $ and $\gamma_2 = z^*_2/z_1$,
defined in \eqref{eq22}.  As $T> T_c \Longleftrightarrow \gamma_2>1$, \eqref{eq64} agrees in the order
of decay with \eqref{eq63} for $\ell=0$ and $m=n$, but now we have the precise rate and constants.  For
$T=T_c$, Wu obtained
\be\label{eq65}
\begin{aligned}
\lf\langle \sg_{1,1}\, \sg_{1, 1+n}\rt\rangle &= e^{\frac{1}{4}}\: 2^{\frac{1}{12}}\:
A^{-3}\, n^{-\frac{1}{4}}\lf(1+ \gamma_1\rt)^{\frac{1}{4}} \, \lf(1-\gamma_1\rt)^{-\frac{1}{4}}\\
&\qquad \times \lf(1+B_1\, n^{-2} + O(n^{-3})\rt)
\end{aligned}
\ee
with an explicit expression for $B_1$ in terms of $\gamma_1$, and $A$ is again Glaisher's constant \eqref{Glai}.
For $T<T_c$, Wu also obtained higher order terms in \eqref{eq15} \eqref{eq16}, i.e., as $n\to\infty$
\be\label{eq66}
\lf\langle \sg_{1,1}\, \sg_{1, 1+n}\rt\rangle = \lf(1-k_{\tons}^2\rt)^{\frac{1}{4}}
\lf(
1+\lf(2\pi \, n^2\rt)^{-1}\:\gamma^{2n}_2 \lf(\gamma^{-1}_2- \gamma_2\rt)^{-2} \lf[1+\frac{c_1}{n}
+ \frac{c_2}{n^2} + \dots \rt]\rt)
\ee
with explicit expression for $c_1$ and $c_2$ in terms of $\gamma_1$ and $\gamma_2$.  In particular as
$\gamma_2<1$, we see that $\lf\langle \sg_{1,1}\, \sg_{1, 1+n}\rt\rangle \to \lf(1-k^2_{\tons}\rt)^{\frac{1}{4}}
=M_0$ exponentially fast.

There is an interesting twist to the story of Wu's computation of \eqref{eq65} at $T_c$.  As a warm-up
problem at $T_c$, $\sinh \frac{2J_1}{k_B \,T_c} \;\sinh \frac{2J_2}{k_B\, T_c} =1$, he considered the
limiting case when $J_1\to 0$ and $J_2\to\infty$ in such a way that
$\gamma_1\downarrow 0$ and $\gamma_2 \uparrow 1$ in \eqref{eq21}, and hence
$$
\varphi_{\tons}\!\lf(e^{i\tht}\rt) = \lf(\frac{1-e^{-i\tht}}{1-e^{i\tht}}\rt)^{\frac{1}{2}}=
e^{i(\pi-\tht)/2}, \qquad 0<\tht < 2\pi.
$$
Wu recognized $\det\!\lf(\lf(\varphi_{\tons}\rt)_{j-k}\rt)_{0\le j\:k\le n-1}$ as a Cauchy determinant and
computed the asymptotics which we have recorded above in \eqref{eq60}.  Wu was unaware at the time, just
like Fisher before him, that he had serendipitously computed the correlation function along the diagonal
at $T_c$!  Had Kaufman and Onsager published their formula for $\lf\langle \sg_{1,1}\, \sg_{1+n, 1+n}\rt\rangle$,
this comedy of errors would probably have been avoided.

In the symmetric situation, $J_1=J_2$, one finds that at $T_c$, $\gamma_1=3-2\sqrt{2}$ and so
$\lf(\frac{1+\gamma_1}{1-\gamma_1} \rt)^{\frac{1}{4}}=2^{\frac{1}{8}}$.  Substituting this relation into
\eqref{eq65}, we obtain as $n\to\infty$
\be\label{eq67}
\lf\langle \sg_{1,1} \, \sg_{1, 1+n}\rt\rangle = e^{\frac{1}{4}} \: 2^{5/24}\, A^{-3}\,n^{-\frac{1}{4}}
\lf(1+O\lf(\frac{1}{n^2}\rt)\rt),\qquad J_1=J_2.
\ee
This relation and \eqref{eq60} are consistent with the formula (cf.~\eqref{eq57})
\be\label{eq88p}
\lf\langle \sg_{1,1}\,\sg_{1+\ell, 1+m}\rt\rangle = e^{\frac{1}{4}} \: 2^{5/24}\, A^{-3}\:k^{-\frac{1}{4}}
\lf(1+O\lf(\frac{1}{k^2}\rt)\rt),\qquad J_1=J_2
\ee
where $k=\sqrt{\ell^2 +m^2} \to \infty$.  This asymptotic formula was derived by Cheng and Wu \cite{ChenWu}
in 1967 (they also consider $J_1\neq J_2$), but it was only established rigorously in the last year by
H.~Pinson \cite{Pin} using a Phragmen-Lindelof type argument to interpolate between \eqref{eq60} and
\eqref{eq67}.  
Recently two additional independent proofs of the rotational symmetry of 
$\lf\langle \sg_{1,1}\,\sg_{1+\ell, 1+m}\rt\rangle$ as $k=\sqrt{\ell^2 +m^2} \to \infty$ have been given by Dub\'edat 
\cite{Dub} and by Chelkak, Hongler, and Izyurov \cite{CheHonIzy}. In \cite{Dub} the author proceeds by equating the product of two Ising correlators with a free field (bosonic) correlator. In \cite{CheHonIzy}, the proof is based on convergence results for discrete holomorphic spinor variables pioneered by S.~Smirnov. In earlier work, Smirnov used such variables to prove conformal invariance and to detect Schramm-Loewner evolution within the Ising model, in an appropriate scaling limit (for references, see \cite{CheHonIzy}).
If $J_1\neq J_2$, one expects elliptical rather than rotational symmetry in the  analog
of \eqref{eq88p} (see \cite{ChenWu} and also \eqref{eq71} below, which is taken from \cite{WMTB}).

The full expansion of $\lf\langle \sg_{1,1}\, \sg_{1+\ell, 1+m}\rt\rangle $ as $\ell^2 +m^2 \to \infty$
at $T\neq T_c$ in terms of so-called form factors $f^{(j)}_{M,N}$, is discussed in detail in  \cite{McC1},
\begin{align}
\lf\langle \sg_{1,1}\, \sg_{1+m, 1+\ell}\rt\rangle &= (1-t)^{\frac{1}{4}} \lf(1+\sum^\infty_{n=1}
f^{(2n)}_{m,\ell} \rt)&\quad\text{for}\quad T<T_c \label{eq68}
\intertext{and}
\lf\langle\sg_{1,1}\,\sg_{1+m, 1+\ell}\rt\rangle &=(1-t)^{\frac{1}{4}} \sum^\infty_{n=0}
f^{(2n+1)}_{m,\ell}&\quad\text{for}\quad T>T_c. \label{eq69}
\end{align}
For $T< T_c$, $t=k^2_{\tons}$ and for $T>T_c$, $t=k^{-2}_{\tons}$.  These expansions were first
obtained by Wu, McCoy, Tracy and Barouch \cite{WMTB} in 1976, but explicit expressions for the form 
factors were only given in 1999 and 2000 by Nickel \cite{Nick1} \cite{Nick2}.

The derivation of the form factor expansions \eqref{eq68} \eqref{eq69} for 
$\lf\langle \sg_{1,1}\, \sg_{1+m, 1+\ell}\rt\rangle$ was obtained in \cite{WMTB} using the Fredholm determinant representations of Cheng and Wu \cite{ChenWu}. 
In the special cases $m=0$, $\ell\to\infty$, or $m=\ell\to\infty$, Lyborg and McCoy \cite{LyMcC} used the finite
Toeplitz determinant representations to derive form factor expansions as in \eqref{eq68} \eqref{eq69}, but the detailed form of the expansions they obtained differs from the general expansions obtained in \cite{WMTB}.
It is an interesting, and still unresolved, question to reconcile the form factor expansions obtained in these two different ways. For recent results on the form factor expansion of the diagonal correlations  $\lf\langle \sg_{1,1}\, \sg_{1+m, 1+m}\rt\rangle$
which utilize results from \cite{BorOk} and \cite{GerCase}, see \cite{WiFo}.

\section{Fisher-Hartwig singularities. Fisher-Hartwig conjecture}\label{sectionFH}
In addition to the new questions and challenges raised in Toeplitz theory by the Ising model for $T\ge T_c$,
other kinds of singularities also began to appear.  For example, in 1964, in his work on the ground state
of the one-dimensional system of impenetrable bosons \cite{Len1} Lenard expressed the one-particle density
matrix for the system as a Toeplitz matrix with a symbol $\varphi\!\lf(e^{i\tht}\rt)$ that vanished at some points on $S^1$
(see \eqref{eq81} and particularly Remark \ref{Lenard} below for further discussion).
In 1968, in a major step in the development and generalization of SSLT, Fisher and Hartwig \cite{FisHart1}
introduced a class of singular symbols for Toeplitz determinants, that allowed for zeros, (integrable)
singularities, discontinuities and non-zero winding numbers.

The symbols of Fisher-Hartwig (FH) class have the following form (we use the notation in \cite{DIK1})
\be\label{eq74}
f(z) =e^{V(z)} \, z^{\sum^m_{j=0} \, \beta_j} \prod^m_{j=0} \lf|z-z_j\rt|^{2\al_j} \, g_{z_j,\, \beta_j}
\, (z) \, z_j^{-\beta_j}, \qquad  z=e^{i\tht},\quad 0 \le \tht <  2\pi,
\ee
for some $m=0, 1, 2, \dots,$  where
\begin{align}
&z_j = e^{i\,\tht_j}, \qquad j=0,1,\dots, m, \qquad 0=\tht_0 < \tht_1 < \dots < \tht_m < 2\pi,
\label{eq75}\\
&g_{z_j\, \beta_j} (z) \equiv g_{\beta_j}(z) = \lf\{
\begin{aligned}
e^{i\pi\beta_j}, \qquad &0 \le \arg \;z< \tht_j \ ,\\
e^{-i\pi\beta_j}, \qquad &\tht_j \le \arg \;z < 2\pi,
\end{aligned} \rt. 
\label{eq76}\\
&\Re\al_j >-\frac{1}{2}, \qquad \beta_j \in \CC, \qquad j=0, 1, \dots, m,
\label{eq77}
\end{align}
and $V\!\lf(e^{i\tht}\rt)$ is a sufficiently smooth function on $S^1$.  Here the condition on $\Re\al_j$
insures integrability.  Note that a Fisher-Hartwig singularity at $z_j, \;\;j=1, \dots, m$, consists
of a root-type singularity
\be\label{eq78}
\lf|z-z_j \rt|^{2\al_j} = \lf| 2\sin \frac{\tht-\tht_j}{2}\rt|^{2\al_j}
\ee
and a jump singularity $z^{\beta_j}\, g_{\beta_j} (z)$ at $z_j$ (observe that 
$z^{\beta_j} \, g_{\beta_j}(z)$ is continuous at
$z=1$ for $j\neq 0$).  For $j=0$, $z_0=1$, in addition to the root-type singularity \eqref{eq78} at $\tht_0=0$, one has
$g_{z_0,\, \beta_0}(z) =e^{-i\pi \beta_0}$, $0\le \tht < 2\pi$ and so $z^{\beta_0}\; g_{z_0, \, \beta_0}(z)
= e^{i \,(\tht-\pi)\,\beta_0}, \; 0 \le \tht< 2\pi$, has a jump at $z=z_0=1$.  The factors $z_j^{-\beta_j}$ are
singled out to simplify comparison with the existing literature, in particular,  \cite{FisHart1}.

Here are some examples of FH symbols:
\begin{itemize}
\item
If $f\!\lf(e^{i\tht}\rt) =-i$ for $0 \le \tht <\pi$ and $+ i$
for $\pi \le \tht < 2\pi$ (cf \cite{BasTr}),
then $z_0=1$, $z_1=-1=e^{i\pi}$, $\al_0 = \al_1=0$, $\beta_0 =\frac{1}{2}$ and $\beta_1=-\frac{1}{2}$, and
\be
\begin{aligned}
f\!\lf(e^{i\tht}\rt) &= z^{\beta_0 + \beta_1} \; g_{1,\, \frac{1}{2}} (z) \;g_{-1, -\frac{1}{2}}(z)
\; z_0^{-\beta_0} \; z_1^{-\beta_1} \\
&= g_{1,\, \frac{1}{2}}(z)\; g_{-1, -\frac{1}{2}}(z) \; e^{i\pi/2}
\end{aligned}
\ee

\item
If $f\!\lf(e^{i\tht}\rt)$ is a sufficiently smooth function on $S^1$ (in particular, a trigonometric polynomial) with  
two zeros at $0 < \tht_1 < \tht_2 < 2\pi$
(cf \cite{BGrM}, \cite{DIK2}), then  with $z_1=e^{i\tht_1}$, $z_2=e^{i\tht_2}$, $\al_1=\al_2=1/2$,
$\beta_1=\frac{1}{2}$, $\beta_2=-\frac{1}{2}$,
\be\label{eqBG}
\begin{aligned}
f(z) &= e^{V(z)} \; z^{\beta_1 + \beta_2}\; \lf|z-z_1\rt| \lf|z-z_2\rt|\, g_{z_1,\; \frac{1}{2}} (z)
\; g_{z_2, \, -\frac{1}{2}}(z) \lf(\frac{z_1}{z_2}\rt)^{-\frac{1}{2}} \\
&= e^{V(z)} \,\lf|z-z_1\rt| \lf|z-z_2\rt| \, g_{z_1, \,\frac{1}{2}}(z) \; g_{z_2,\,-\frac{1}{2}}(z)
\lf(z_1/z_2\rt)^{-\frac{1}{2}}
\end{aligned}
\ee
for a suitable function $V\!\!\lf(e^{i\tht}\rt)$ on $S^1$.

\item
The following symbol arose in Lenard's work in \cite{Len1} on impenetrable bosons mentioned above,
\be\label{eq81}
f\!\lf(e^{i\tht}\rt) = \lf|z-e^{i\tht_1}\rt| \lf|z-e^{i\tht_2}\rt|
\ee
where $\tht_1 \neq \tht_2\,(\text{mod } 2\pi)$.  Here $z_1=e^{i\tht_1}$, $z_2=e^{i\tht_2}$, $\al_1=\al_2=1/2$ and
$\beta_1=\beta_2=0$.

\item
For the Ising model at $T< T_c$, the winding number of $\varphi_{\tons}\!\lf(\eit\rt)$ is zero and so $\varphi_{\tons}
(z) = e^{V(z)}$ for a suitable smooth function $V\!\!\lf(\eit\rt)$ on $S^1$.  For $T=T_c$, we have from \eqref{eq25} with
$0< \gamma_1 < 1$,
\be\label{eq82}
\begin{aligned}
\varphi_{\tons}\! \lf(\eit\rt) &=  \frac{1-\gamma_1 \, \eit}{\lf|1-\gamma_1 \, e^{-i\tht}\rt|} \quad i\,e^{-i\tht/2},
\qquad  0\le \tht < 2\pi,\\
&= e^{V(z)} \; z^{\beta_0} \; g_{z_0, \, \beta_0}(z) \, z_0^{-\beta_0}
\end{aligned}
\ee
where $z_0=1$, $\al_0=0$, $\beta_0=-\frac{1}{2}$ and $V\!\lf(\eit\rt)$ is smooth on $S^1$.
For $T>T_c$, we have from \eqref{eq24} with $0 < \gamma_1 <1< \gamma_2$,
\be\label{106}
\begin{aligned}
\varphi_{\tons} \!\lf(\eit\rt) &= \frac{\lf(1-\gamma_1\, \eit\rt) \lf(\gamma_2-\eit\rt)}{\lf|1-\gamma_1 \,e^{-i\tht}\rt|
\lf|1-\gamma_2 \,\eit\rt|} \lf(-e^{-i\tht}\rt) \\
&= e^{V(z)} \, z^{\beta_0}\, g_{z_0, \, \beta_0}(z) \, z_0^{-\beta_0}
\end{aligned}
\ee
where $z_0=1$, $\al_0=0$, $\beta_0=-1$ and again $V(z)$ is smooth on $S^1$.
\end{itemize}

From the point of view of Toeplitz determinants with FH singularities, Wu was the  first to obtain \cite{Wu}
detailed asymptotic results for symbols with discontinuities as in \eqref{eq82} and for symbols with non-zero winding 
as in \eqref{106}.  Lenard was the first to obtain \cite{Len1} \cite{Len2} detailed asymptotic results for
symbols which have zeros on $S^1$ as in \eqref{eq81}.  For symbols \eqref{eq74} with $\beta_j=0$ and
$\al_j$ real, $\al_j>-\frac{1}{2}$, $j=0, \dots, m$, Lenard made a general conjecture \cite{Len2}\footnote{
Note that although \cite{Len2} was published in 1972, a preliminary unpublished version of the paper
was already in existence in 1968 (see \cite{FisHart1}).}
about the asymptotic form of $D_n(f)$,
\be\label{eq83} 
D_n(f) \sim E\lf(e^V, \, \al_0, \dots, \al_m, \, \tht_0, \dots, \tht_m\rt) n^{\sum^m_{j=0} \al_j^2}
\, e^{n V_0}
\ee
as $n\to\infty$, where $V_0=\frac{1}{2\pi} \int^{2\pi}_0 V\!\lf(\eit\rt)d\tht$ as in
\eqref{eq32} \eqref{eq49}, and $E\!\lf(f, \, \al_0, \dots, \al_m, \tht_0, \dots, \tht_m\rt)$ is some (unspecified) constant.
Moreover he was able to verify this conjecture for the case $f\!\lf(\eit\rt)=\lf|z-e^{i\,\tht_1}\rt|
\lf|z-e^{i\tht_2}\rt|$ in \eqref{eq81}.  

\begin{remark}\label{Lenard}
As noted above, the symbol \eqref{eq81} arose in Lenard's study of impenetrable bosons \cite{Len1}.  More precisely,
Lenard considered $N$ free bosons in one dimension confined to a box of length $L$ with periodic boundary 
conditions and lying in the (symmetric) ground state of the system. According to the criterion of Penrose and Onsager
\cite{PenOns},
Bose-Einstein condensation takes place if the largest eigenvalue $\lambda_{N,L}^{(\mathrm{max})}$ of
the density matrix $\rho_{N,L}(x,x')$ for the system is of order $N$ as $N\sim L\to\infty$.  
For the system at hand $\rho_{N,L}$ is translation invariant,  $\rho_{N,L}(x,x')=\rho_{N,L}(x-x')$,
so that the eigenvalues of $\rho_{N,L}$ are just the Fourier coefficients $\rho_{N,L}^{(n)}$ of $\rho_{N,L}(\cdot)$.
A simple calculation shows that  $\lambda_{N,L}^{(\mathrm{max})}=\rho_{N,L}^{(0)}$, the zeroth Fourier coefficient.
Writing $\rho_{N,L}(\xi)=\frac{1}{L}R_N(2\pi\xi/L)$, we have that
\be\label{r34-1}
\frac{\lambda_{N,L}^{(\mathrm{max})}}{N} =\frac{1}{2\pi N}\int_{-\pi}^{\pi}R_N(t)dt.
\ee
(Note that in this case $\lambda_{N,L}^{(\mathrm{max})}$ is actually independent of $L$.) Now it turns out that 
\be\label{34p3}
R_N(t)=D_{N-1}(f_t)
\ee
where $f_t$ is precisely a symbol of type \eqref{eq81},
\be\label{Lensym}
f_t\lf(\eit\rt)=\lf|\eit-e^{it/2}\rt|\,\lf|\eit-e^{-it/2}\rt|.
\ee
(Formula \eqref{34p3} was also obtained independently by Dyson, see \cite{Len1}.)
At that point (April 10, 1963, according to a letter Lenard has kindly made available to the authors)
Lenard asked Szeg\H o for help in evaluating the Toeplitz determinant $D_{N-1}(f_t)$. On July 10, 
three months later, Szeg\H o wrote back to Lenard with the information that, although he had not been able to 
obtain precise asymptotics, he had obtained the bound
\be\label{r34-2}
\lf| R_N(t)\rt| \le \lf|\frac{eN}{\sin(t/2)}\rt|^{1/2}.
\ee
Substituting \eqref{r34-2} into \eqref{r34-1}, we see that
\be
\frac{\lambda_{N,L}^{(\mathrm{max})}}{N} =O\lf(N^{-1/2}\rt)\to 0,\quad\mbox{as}\quad N\to\infty.
\ee
In a letter of thanks to Szeg\H o on July 23, Lenard describes the situation as follows:
``This has great physical interest: It may be expressed as saying that for the model considered there is no
Bose-Einstein condensation in momentum space''.\footnote{This confirms an earlier result of Schultz \cite{Sch} 
proving the absence of condensation. Schultz used other methods which produced the (weaker) bound $O(N^{-4/\pi^2})$.}

Following a further suggestion contained in the above correspondence with Szeg\H o in 1963, Lenard 
considered, in particular, symbols of the form
\be\label{symJac}
f^{(\lb,\mu)}(\eit)=|\eit-e^{i\pi/2}|^\lb |\eit-e^{-i\pi/2}|^\mu,\qquad \lb,\mu>0.
\ee
He started with the general formula for positive symbols $f\!\lf(\eit\rt) >0$
\be\label{93}
D_n (f) = \prod^{n-1}_{k=0} \chi^{-2}_k
\ee
where $\chi_k>0$ is the leading coefficient of the orthonormal polynomial $p_k (z) =\chi_k\, z^k +
\dots, k \ge 0$, with respect to the weight  of $d\mu=f\!\lf(\eit\rt) d\theta/(2\pi)$ on  $S^1$ as in \eqref{eq35}.
For $f^{(\lb,\mu)}(\eit)$ in \eqref{symJac}, the polynomials $p_k(z)$
can be expressed in terms of classical  Jacobi polynomials on the interval $[-1,1]$.  This then leads to
explicit formulae for the $\chi_k$'s in terms of Euler's gamma functions. 

After Lenard published the bound \eqref{r34-2} for the impenetrable boson problem (a Plasma Physics Laboratory preprint of 
\cite{Len1} appeared in 1963), Dyson \cite{Dy3} took up the problem of computing the asymptotics of 
\eqref{r34-1} precisely as $N\to\infty$. 
In order to control the singularity in the symbol $f_t(\eit)=|\eit-e^{it/2}||\eit-e^{-it/2}|$,
Dyson resorted to his interpretation of $D_{N-1}(f_t)$ as a Coulomb gas at inverse temperature $\beta=2$ with $N-1$ unit charges free to move around on the unit circle, but now with two additional half-charges. The probability of finding the half-charges at an angular distance $t$ apart is then proportional to 
$P(t)=|2\sin{t\over 2}|^{1/2}R_N(t)$. 
On physical grounds, Dyson then argued that for distances $t$ much greater than the ``Debye length'' 
$\lb=2\pi/N$,
the two half-charges are effectively screened off from each other, and so $P(t)=\const$ for $t>>\lb$. In other words,
$R_N(t)=k_N\lf(\sin{t\over 2}\rt)^{-1/2}$ for $t>>N^{-1}$, where $k_N$ is a constant. However, clearly $k_N=R_N(\pi)$,
which can be evaluated explicitly as $N\to\infty$ since
\[
f_t(\eit)|_{t=\pi}=f^{(\lb,\mu)}(\eit)|_{\lb=\mu=1}
\]
(see discussion of \eqref{symJac} above). Substituting the result into \eqref{r34-1} and neglecting the contribution from 
$|t|<N^{-1}$, Dyson obtained, as $N\to\infty$,
\be\label{20.1}
\frac{\lambda_{N,L}^{(\mathrm{max})}}{N} \sim C N^{-1/2},\qquad C=\lf({e\over\pi}\rt)^{1/2}2^{-5/6}A^{-6}
\Gamma\lf({1\over 4}\rt)^2
\ee
where $A$ is again Glaisher's constant \eqref{Glai}. The crux of the above calculation is clearly the explicit evaluation of $R_N(\pi)$ as $N\to\infty$.

More detailed analysis of \eqref{r34-1} requires further study of the double-scaling limit
as $t\to 1$, $N\to\infty$. We return to this issue in Section \ref{sec-ds}.


\end{remark}

Inspired by the work of Wu in \cite{Wu}, and perhaps also influenced by the calculations 
of Lenard in \cite{Len1} \cite{Len2}, Fisher and Hartwig \cite{FisHart1} made a conjecture about the
asymptotic form of $D_n(f)$ for general symbols \eqref{eq74},
\be\label{eq84}
\begin{aligned}
D_n(f)  &= E\lf(e^V, \, \al_0, \dots, \al_m, \, \beta_0, \dots, \beta_m,\, \tht_0, \dots, \tht_m\rt)\\
&\quad \times \;n^{\sum^m_{j=0} \lf(\al^2_j-\beta^2_j\rt)} \, e^{n V_0}\lf(1+o(1)\rt),\qquad 
V_0=\frac{1}{2\pi} \int^{2\pi}_0 V\!\lf(\eit\rt)d\tht
\end{aligned}
\ee
as $n\to\infty$, together with a conjecture on the general form of the constant  
$E$. In \cite{FisHart2}, Hartwig and Fisher
present various examples,  heuristics and theorems in support of their conjecture.

In the remarkable paper \cite{Wid2}, representing ``a jump of several quanta in depth and sophistication''
(see MathSciNet MR0331107), Widom verified Lenard's conjecture for complex $\al_j$ with $\Re \al_j > -\frac{1}{2},
\; j=0, \dots, m$, and $V\!\lf(\eit\rt)$ suitably smooth, and obtained an explicit expression for the constant
\be\label{eq85}
\begin{aligned}
E\lf(e^V, \, \al_0, \dots, \al_m, \, \tht_0, \dots, \tht_m\rt)
=&  E\!\lf(e^V\rt) \prod_{0\le j< k \le m} \lf| e^{i\tht_j} -e^{i\tht_k} \rt|^{-2\al_j\al_k} \\
& \times \prod^m_{j=0} e^{-\al_j \widehat{V}\lf(e^{i\tht_j}\rt)} \times \prod^m_{j=0} E_{\al_j}
\end{aligned}
\ee
where
\begin{alignat}{3}
&E\!\lf(e^V\rt) = e^{\sum^\infty_{k=1} k\, V_k\, V_{-k}}, \qquad &V_k &= \text{ Fourier coefficient of } V\!
\lf(\eit\rt),\label{eq86}\\
&\widehat{V}\!\lf(\eit\rt) = V\!\lf(\eit\rt) - V_0&& \label{eq87}\\
\intertext{and}
&E_\al = G^2 (1+\al)/G(1+2\al), \qquad &G(z) &= \text{ Barnes G-function (see \cite{Bar}).} \label{eq88}
\end{alignat}
Formula \eqref{eq85} is consistent with and confirms the general form of the conjecture
\cite{FisHart1} for the constant term in the case $\beta_0 = \beta_1 = \dots = \beta_m =0$.

In \cite{Bas1} Basor considered complex $\al_j$ with $\Re  \al_j >-\frac{1}{2}$, $\beta_j$ pure imaginary,
$j=0, \dots, m$, and $V\!\lf(\eit\rt)$ suitably smooth, and verified \eqref{eq84} with constant term, again
consistent with the general conjecture in \cite{FisHart1},
\be\label{eq89}
\begin{aligned}
&E \lf(e^V, \al_0, \dots, \al_m, \beta_0, \dots, \beta_m, \tht_0, \dots, \tht_m\rt)\\
&\quad =  E(e^V) \prod^m_{j=0} \lf[b_+ (z_j)^{-\al_j + \beta_j}\, b_-(z_j)^{-\al_j-\beta_j} \rt] \\
& \qquad \qquad \times \prod_{0\le j < k \le m} \lf[ \lf|z_j-z_k\rt|^{2\lf(\beta_j\beta_k-\al_j\al_k\rt)}
\lf(\frac{z_k}{z_j \, e^{i\pi}}\rt)^{\al_j \beta_k-\al_k\beta_j}\rt] \\
&\qquad \qquad\times \prod^m_{j=0} \frac{G\lf(1+\al_j + \beta_j\rt)  G \lf(1+ \al_j - \beta_j\rt)}{G\lf(1+
2\al_j\rt)} \ .
\end{aligned}
\ee
with $E\!\lf(e^V\rt)$ and $G(z)$ as above, and
\be\label{eq90}
b_+(z) =e^{\sum^\infty_{k=1}  V_k\,z^k}, \qquad b_-(z)=e^{\sum^{-\infty}_{k=-1} V_k \, z^k}\ .
\ee
The branches in \eqref{eq89} are determined naturally, $b_+(z_j)^{-\al_j +\beta_j} =
e^{\lf(-\al_j + \beta_j\rt)\sum^\infty_{k=1} V_k z_j^k}$ etc, and 
\[
\lf(\frac{z_k}{z_j\,e^{i\pi}}\rt)^{
\lf(\al_j\beta_k-\al_k\beta_j\rt)} = e^{i\lf(\tht_k-\tht_j-\pi\rt)\lf(\al_j \beta_k-\al_k \beta_j\rt)},
\qquad 0 \le \tht_j< \tht_k < 2\pi.
\]

In the paper \cite{BasHelt} discussed above, Basor and Helton also introduce the first example
of a so-called ``separation'' theorem of the following form: For symbols $\varphi$ and $\psi$ 
of a certain prescribed type, including FH symbols with disjoint singularities, they show that
\be
\lim_{n\to\infty}\frac{D_n(\varphi\psi)}{D_n(\varphi)D_n(\psi)}\equiv L(\varphi,\psi)
\ee
exists. Moreover, they provide an explicit form for $L(\varphi,\psi)$. Iterating this result, they are able to reduce
the problem for general FH symbols to adding in pure FH singularities of the form
$z^\beta |z-\widehat z|^{2\al}g_{\widehat z,\beta}(z)\widehat z^{-\beta}$, $|\widehat z|=1$, one at a time.
They analyze the asymptotics of Toeplitz determinants with such pure (i.e., with $V\equiv 0$) 
FH singularities in the following way.
In addition to the results described above, in \cite{Wid2} Widom also verified \eqref{eq84} for a single FH
singularity with $|\Re\al |<\frac{1}{2}$ and $|\Re\beta |<\frac{1}{2}$, but without determining the exact value of the 
constant $E$. But then Basor's formula \eqref{eq89} for $E$ obtained in \cite{Bas1} for $\Re\al_j>-\frac{1}{2}$
and $\beta_j$ pure imaginary, together with Vitali's theorem, directly yields the desired explicit asymptotics for
a single pure FH singularity with  $|\Re\al |<\frac{1}{2}$ and $|\Re\beta |<\frac{1}{2}$. (In fact, as discovered later
by B\"ottcher and Silbermann, see \eqref{eqBS} below, there exists an exact formula for a Toeplitz determinant with
a single pure FH singularity in terms of the Barnes G-functions, which easily yields the asymptotics).
In the end, Basor and Helton are able to verify \eqref{eq84} \eqref{eq89} for $\al_j$, $\beta_k$ in the open set
$\{\max_{j,k}(|\al_j|,|\beta_k|)<\delta\}\subset\CC^{2n}$ for some (small) $\de>0$. 

In \cite{Bas3}, Basor showed using a separation theorem that \eqref{eq84} \eqref{eq89} also hold in the case $\al_j=0$, 
$|\Re\beta_j|<\frac{1}{2}$, $j=0,\dots,m$.

B\"ottcher and Silbermann \cite{BottSilb5} finally proved \eqref{eq84} with Basor's constant \eqref{eq89} under the 
assumption that  $|\Re\al_j|<1/2$ and $|\Re\beta_j|<1/2$ for all $j$. This is the most general setting in which the conjecture \eqref{eq84} is true in its original form, that is, with the exponent of $n$ conjectured to be 
$\sum(\alpha_j^2-\beta_j^2)$,
without the need to take care of degenerate situations and of the ambiguity in the $\beta$'s (see \eqref{eq97} below).
We note that in \cite{BottSilb6} it had already been 
shown that \eqref{eq84} does not hold in general if $\al_j+\beta_j$ and $\al_j-\beta_j$ are nonnegative integers for all 
$j=0,\dots,m$, and a corrected version of the conjecture was proved in this case (this is a particular case of a formula 
conjectured by Basor and Tracy: see next section). 
The developments up to the 1990's are described in \cite{Bott1} and \cite{Ehr} and also in the books \cite{BottSilb2} \cite{BottSilb3}. We note that certain separation theorems also play a key role in the work of B\"ottcher and Silbermann. Another important element in their proofs is an explicit formula, noted above, for a Toeplitz determinant with a single, pure FH singularity found by the authors in \cite{BottSilb5} (see also \cite{BottWid2} \cite{BasChe} for later alternative derivations). Namely, if $V(z)\equiv 0$, $m=0$, then
with $\alpha_0\equiv\alpha$, $\beta_0\equiv\beta$,
\be\label{eqBS}
\begin{aligned}
D_n(f)=&\frac{G(1+\al+\beta)G(1+\al-\beta)}{G(1+2\al)}
\frac{G(n+1)G(n+1+2\al)}{G(n+1+\al+\beta)G(n+1+\al-\beta)},\\
&\Re\al>-{1\over 2},\qquad\al\pm\beta\neq -1,-2,\dots,\qquad n\ge 1.
\end{aligned}
\ee
The condition on $\al\pm\beta$ is needed because the $G$-function has zeros at $z=0,-1,-2,\dots$
The asymptotics for \eqref{eqBS}
follow from the asymptotics of the $G$-function (see \cite{Bar})
\be\label{eq92} 
\begin{aligned}
\log G(t+a+1) &= \frac{1}{12} -\log A - \frac{3t^2}{4} - at + \frac{t+a}{2} \: \log (2\pi)\\
&\quad + \lf( \frac{t^2}{2} + at + \frac{a^2}{2} - \frac{1}{12}\rt) \log t+o(1), \qquad \text{as } \; t\to\infty
\end{aligned}
\ee
where again  $A$ is Glaisher's constant \eqref{Glai}. The resulting asymptotics are consistent with
\eqref{eq84} \eqref{eq89}.

In view of the Heine representation \eqref{eq40}, it is not surprising that the formula \eqref{eqBS} is 
related to the Selberg integral. In fact, independently of \cite{BottSilb5}, a more general formula 
was obtained by Selberg himself:
\be\label{Selb}\begin{aligned}
\int^{2\pi}_0 \dots \int^{2\pi}_0 &\prod_{0\le j < k \le n-1}
\lf|e^{i\tht_j} - e^{i\tht_k} \rt|^{2\gamma} \prod^{n-1}_{j=0} e^{i(\tht_j-\pi)\beta}\lf|1-e^{i\tht_j}\rt|^{2\al}
\,\frac{d\tht_j}{2\pi}\\
=&
\prod_{j=0}^{n-1}\frac{\Gamma(1+2\al+j\gamma)\Gamma(1+(j+1)\gamma)}
{\Gamma(1+\al+\beta+j\gamma)\Gamma(1+\al-\beta+j\gamma)\Gamma(1+\gamma)}
\end{aligned}
\ee
where
\[
\Re\al>-{1\over 2},\qquad \Re\gamma>-\min\lf\{\frac{1}{n},\frac{\Re(2\al+1)}{n-1}\rt\}.
\]
In a published form, this formula first appeared in the thesis of Morris in 1982 
(see \cite{ForWar} for a historic account).
The interest in \eqref{Selb} was due to the efforts to prove one of Dyson's conjectures in the theory of
random matrices in \cite{Dy1}.
If we set $\gamma=1$ in \eqref{Selb}, and use the relation
$G(z+1)=\Gamma(z)\, G(z)$,
we recover \eqref{eqBS} from \eqref{eq40}.

The Barnes' function $G(z)$ is a basic object in analysis, bearing the same relation to the gamma function as the gamma function bears to the function $z$, $\Gamma(z+1)
=z \, \Gamma(z)$.  Just as $\log \Gamma(z)$ is a sum of values of $\log z$,  so $\log G(z)$ is the
iterated sum.  
It is of interest to delineate some of the analytic pathways by which $G(z)$ appears in the
asymptotics of $D_n(f)$.  

The most direct example, a particular case of \eqref{eqBS}, arises for the diagonal correlation
function $\lf\langle\sg_{1,1}\, \sg_{1+n, 1+n}\rt\rangle$ at $T_c$, where $\varphi_{\tdg}(\tht)=e^{i\lf(\pi-\tht\rt)/2}$, i.e. $V\equiv 0$, $m=0$, $\al_0=0$, $\beta_0=-1/2$,
 in \eqref{eq58} gives rise to the Cauchy determinant \eqref{Opure} which can be evaluated as in \eqref{eq53}
in terms of a product of gamma functions, and hence in terms of $G$-functions, 
$\prod^{n-1}_{q=1} \Gamma(q+1)=G(n+1)/G(2)$, etc.  


For a system of polynomials orthonormal with respect to a function $f$ on the unit circle, 
the leading coefficients $\chi_k$ are related to $D_n(f)$ by the product formula \eqref{93}.
As we discussed above, in the case of $f$ given by \eqref{symJac} the coefficients  $\chi_k$
are expressed in terms of gamma functions, and therefore $G$-functions appear in the product.
More generally, there is a certain universality in the asymptotic behavior of orthogonal
polynomials associated with FH symbols on $S^1$ (see \cite{DIK1} \cite{DIK3}), and, mutatis mutandis, 
for general FH symbols $f$, the $G$-functions enter the arena through the same stage door.

As noted above in \eqref{eq82}, at $T=T_c$, $\varphi_{\tons}\!\lf(\eit\rt)$ has a single FH singularity at
$z_0=1$ with $\al_0=0$ and $\beta_0=-\frac{1}{2}$.  The first confirmation of \eqref{eq84} \eqref{eq89}
which allowed for these values was given in \cite{BottSilb4}, where the authors proved the result for
$m=0$, $\Re \al_0 \ge 0$, $\Re\lf(\al_0 + \beta_0\rt) >-1$ and  $\Re\lf(\al_0-\beta_0\rt) > -1$.  It is
instructive to recover Wu's asymptotic formula \eqref{eq65} from \eqref{eq84} \eqref{eq89}.
Here $V\!\lf(\eit\rt)= \frac{1}{2} \log \lf(\frac{1-\gamma_1\,\eit}{1-\gamma_1 \, e^{-i\tht}}\rt)$
and we find $V_0=0$ and $V_k=-\gamma_1^{|k|}/2k$
for $k \neq 0$.  Thus $E\lf(e^V\rt)=e^{-\frac{1}{4} \sum^\infty_{k=1} \, k \,\gamma_1^{2k}/k^2} = \lf(1-\gamma^2_1
\rt)^{\frac{1}{4}}$.  On the other hand   $\lf(b_+(1)\rt)^{-\frac{1}{2}} \lf(b_-(1)\rt)^{\frac{1}{2}}=
e^{\frac{1}{4} \sum^\infty_{k=1}\gamma^k_1/ k} \, e^{\frac{1}{4} \sum^{-\infty}_{k=-1} \gamma_1^{|k|}/|k|} =
\lf(1-\gamma_1\rt)^{-\frac{1}{2}}$.  But $G\lf(\frac{1}{2}\rt)=2^{\frac{1}{24}} \, e^{\frac{1}{8}}
\,\pi^{-1/4} \, A^{-3/2}$, where again $A$ is Glaisher's constant (see \cite{Bar}), and
we have $G\lf(\frac{3}{2}\rt) = \Gamma \lf(\frac{1}{2}\rt)
G\lf(\frac{1}{2}\rt) =2^{\frac{1}{24}} \, e^{\frac{1}{8}}\, \pi^{\frac{1}{4}}\, A^{-3/2}$.  Thus $G\lf(\frac{1}{2}
\rt) G\lf(\frac{3}{2}\rt)=2^{\frac{1}{12}} \, e^{\frac{1}{4}}\, A^{-3}$.  Substituting these relations into
\eqref{eq89}, we obtain $E\!\lf(e^V, 0, \frac{1}{2}, \tht_0=0\rt)=2^{\frac{1}{12}} \, e^{\frac{1}{4}}\, A^{-3}
\lf(1+\gamma_1\rt)^{\frac{1}{4}} \lf(1-\gamma_1\rt)^{-\frac{1}{4}}$, and as $n^{-\beta^2_0} = n^{-\frac{1}{4}}$,
we arrive at precisely  the leading order term in \eqref{eq65}.

The {\bf definitive result} for the leading order asymptotics of Toeplitz determinants with FH singularities was
obtained by Ehrhardt in his PhD thesis at TU Chemnitz in 1997 (see \cite{Ehr}).  Let $\Bl$ denote the
seminorm
\be\label{eq91}
\Bl = \max_{j,k} \lf|\Re \beta_j - \Re \beta_k\rt|
\ee
where $1\le j, \, k \le m$ if $\al_0= \beta_0=0$, and $0\le j, \, k\le m$ otherwise.  If $m=0$, set $\Bl=0$.
Note that in the case of a single singularity, we always have $\Bl=0$.
If $V\!\lf(\eit\rt)$ is $C^\infty$ on $S^1$, $\Bl <1$, $\Re \al_k>-\frac{1}{2}$ and $\al_j \pm \beta_j
\neq -1, -2, \dots $ for $j,\, k=0,1, \dots, m$, then Ehrhardt proved the Fisher-Hartwig conjecture
\eqref{eq84} with Basor's constant \eqref{eq89}.  As above, the conditions on $\al_j \pm \beta_j$ ensure that 
Basor's constant is non-zero and \eqref{eq89}
indeed provides the leading order asymptotics for $D_n(f)$ as $n\to \infty$. 
A key element in Ehrhardt's proof is a suitably generalized version of the separation theorem
in \cite{BasHelt}.\footnote{Note that the RHS of \eqref{eq84} and \eqref{eq89} make sense for any $\al_j$,
$2\al_j\neq -1,-2,\dots$. In fact, \eqref{eq84} \eqref{eq89} still hold in this case. This extension to the left of
the lines $\Re\al_j=-1/2$ is explained in \cite{Ehr} and corresponds to replacing the symbol $f$ by an appropriate
distribution.}

An independent proof of
Ehrhardt's result was given recently in \cite{DIK3} where only a finite degree of smoothness is needed for
$V\!\lf(\eit\rt)$. In \cite{DIK3}, estimates for the error term in the asymptotics are also provided. 
If $V(z)\in C^{\infty}$ on $S^1$, the asymptotics \eqref{eq84} \eqref{eq89}
hold with error term
\be
o(1)=O(n^{\Bl-1}).
\ee
Note that the smallness of the error term breaks down if $\Bl\ge 1$. We will consider this situation in detail in 
the next section. The main steps of the method of \cite{DIK3} are as follows. First, the authors derive  
differential identities for the logarithmic derivatives of $D_n(f)$, in particular, for 
\be
{\partial\over\partial\al_j}\log D_n(f),\qquad {\partial\over\partial\beta_j}\log D_n(f),
\ee
in terms of the associated OPUC's of degrees $n$ and $n+1$ only.
Second, the authors obtain the asymptotics of these OPUC's using their Riemann-Hilbert representation and the 
steepest descent analysis of the Riemann-Hilbert problem for large $n$. 
Third, the asymptotic formulas so obtained are substituted 
into the RHS of the differential identities, and the final result \eqref{eq84} \eqref{eq89}
is obtained by integration.\footnote{
Note that if one just used the formula \eqref{93} instead of the differential identities, then the 
product $\prod_{k=0}^{n_0}\chi_k^{-2}$ for small $n_0$, and hence the
constant factor $E$ in \eqref{eq84}, would remain undetermined.}
In \cite{DIK1}, \cite{DIK3} the authors also note various uniformity properties of the large $n$ asymptotics for Toeplitz determinants in the parameters, such as $\al_j$, $\beta_j$, on which $f$ depends. This issue becomes
important in discussions of double-scaling limits (cf Section \ref{sec-ds}). 

For the Ising model at $T> T_c$, we see from \eqref{106} that $z_0=1,\, \al_0=0$ and $\beta_0=-1$,
and so $\al_0 + \beta_0=-1$.  Thus \eqref{eq84} \eqref{eq89} do not provide the leading order asymptotics for
$\lf\langle\sg_{1,1}\,\sg_{1,1+n}\rt\rangle$ in this case.  One must analyze lower order terms, in fact go
``beyond all orders'' in the language of Kruskal, to capture the exponential decay in Wu's result \eqref{eq64}
(see Remark \ref{degen} below).

\section{Basor-Tracy conjecture}
As already observed in \cite{BottSilb5} \cite{BottSilb6},  
Ehrhardt's condition on the seminorm, $\Bl <1$, is not just a technical matter.  In a critical calculation in
1991, Basor and Tracy \cite{BasTr} considered a particular symbol $f\!\lf(\eit\rt)$ which has two FH
singularities at $z_0=1$ and $z_1=e^{i\,\pi}=-1$ respectively
\be\label{eq94}
f\!\lf(\eit\rt) =  g_{1, \frac{1}{2}}(z) \, g_{-1, -\frac{1}{2}}(z)\, e^{i\,\pi/2}
\ee
with $\beta_0= \frac{1}{2}, \, \beta_1 = -\frac{1}{2}$ so that $\Bl=1$.  By direct computation they found
that as $n\to \infty$
\be\label{eq95}
D_n(f) = \frac{1+(-1)^n}{2} \sqrt{\frac{2}{n}} \; G\lf(\frac{1}{2}\rt)^2  G\lf(\frac{3}{2}\rt)^2 \lf(1+O(1)\rt)
\ee
which is not of the general FH asymptotic form.  Basor and Tracy observed, however, that $f\!\lf(\eit\rt)$
had another FH representation with $\beta_0=-\frac{1}{2}$ and $\beta_1= \frac{1}{2}$,
\be\label{eq96}
f\!\lf(\eit\rt) =e^{i\,\pi} \:g_{1, -\frac{1}{2}} (z) \: g_{-1, \frac{1}{2}}(z) \, e^{-i\,\pi/2} \ .
\ee
If one substitutes the values 
$V=0$, $\al_0=\al_1=0$, $\beta_0=\half$, $\beta_1=-\half$ into \eqref{eq84} \eqref{eq89},
one obtains for \eqref{eq94}
\begin{align*}
&n^{-\lf(\frac{1}{4} + \frac{1}{4}\rt)} \,\lf|1-e^{i\pi}\rt|^{-\half} \:\frac{G\lf(\frac{3}{2}\rt) G\lf(\half\rt)
}{G(1)} \:\frac{G\lf(\half\rt) G\lf(\frac{3}{2}\rt)}{G(1)} \; \lf(1+o(1)\rt)\\
&\qquad \qquad = (2n)^{-\half} \: G\lf(\tfrac{3}{2}\rt)^2 \: G\lf(\tfrac{1}{2}\rt)^2 \lf(1+ o(1)\rt)
\end{align*}
as $G(1)=1$.  On the other hand, if one substitutes the values $V=i\,\pi$, $\al_0=\al_1=0$, $\beta_0=-\half$,
$\beta_1=\half$, one obtains for \eqref{eq96}
$$
e^{in\pi} \,(2n)^{-\half} \, G^2 \lf(\tfrac{1}{2}\rt) G^2\lf(\tfrac{3}{2}\rt) \lf(1+ o(1)\rt).
$$
Just adding these two asymptotic forms, we obtain precisely the Basor-Tracy result \eqref{eq95}.  So we see that if
a symbol $f\!\lf(\eit\rt)$ has more than one FH representation, then the asymptotics of $D_n(f)$ is given,
at least in this case, by some combination of the Fisher-Hartwig asymptotics for each of them.  Based on this
example, Basor and Tracy \cite{BasTr} made the following bold, general conjecture.  First note the following.
Let $f\!\lf(\eit\rt)$ be a Fisher-Hartwig symbol as in \eqref{eq74}.  Let $n_0, n_1, \dots, n_m$ be a set
of integers with $\sum^m_{j=0} n_j=0$, and let $f\lf(z; \wbe\rt)$ denote the FH symbol obtained by
replacing $\beta_j$ with $\wbe_j= \beta_j+n_j$, $0\le j\le m$, and $e^{V(z)}$  by $\lf(\prod^n_{j=0}
z_j^{n_i}\rt) e^V = e^{V+i \sum^m_{j=0} \, n_j\, \tht_j}$.  Then  $f(\wbe)=f\lf(z;\wbe\rt)$
gives another FH representation for $f(z)$,
\be\label{eq97}
f(z) =f\lf( z;\wbe\rt), \qquad \wbe_j = \beta_j +n_j, \quad \sum^m_{j=0} n_j=0 \ .
\ee
Moreover, it is easy to see that all FH representations of $f(z)$ arise in this way.  Given $\beta$, we call
\be\label{eq98}
O_\beta = \lf\{ \wbe : \wbe_j = \beta_j + n_j, \quad \sum^n_{j=0} n_j=0 \rt\}
\ee
the \tit{orbit} of $\beta$.  We consider the discrete minimization problem
\be\label{eq99}  
F_\beta = \min_{\wbe \in O_\beta} \lf(\sum^m_{j=0} \lf(\Re \wbe_j\rt)^2\rt).
\ee
It is clear that the minimum is indeed obtained in $O_\beta$.  Let $\CM_\beta=\lf\{\wbe\in O_\beta
: \sum^n_{j=0} \lf(\Re \wbe_j\rt)^2 = F_\beta\rt\}$.  We say $\CM_\beta$ is \tit{non-degenerate}
if $\al_i \pm \wbe_j \neq -1, -2, \dots$ for all $j=0, \dots, m$ and all $\wbe \in \CM_\beta$.

\begin{conject}\cite{BasTr}\label{conject8}
Let $f\!\lf(\eit\rt)$ be given as in \eqref{eq74} with $V\!\lf(\eit\rt)$ sufficiently smooth on $S^1$ and $\Re \al_j
>-\frac{1}{2}$, $0\le j\le m$.  Suppose $\CM_\beta$ is non-degenerate.  Then as $n\to\infty$,
\be\label{eq100}
D_n(f) = \sum_{\wbe\in \CM_\beta} \lf[R_n \lf(f(\wbe)\rt) (1+o(1))\rt].
\ee
Each $R_n\lf(f(\wbe)\rt)$ stands for the FH asymptotic form \eqref{eq84} \eqref{eq89} with $\wbe$
and $\lf(\prod^m_{j=0} z_j^{n_j} \rt) e^V$ in place of $\beta$ and $e^V$ respectively.
\end{conject}

Note that by the definition of $\CM_\beta$, all the terms $R_n\lf(f(\wbe)\rt)$ in \eqref{eq100} have the
same order of magnitude.  Indeed for all $\wbe \in \CM_\beta$, $\Re \sum^m_{j=0}\, \wbe_j^2=
F_\beta-\sum^m_{j=0} \lf(\Im \beta_j\rt)^2$.

To see how this conjecture works for $f\!\lf(\eit\rt)$ in \eqref{eq94}, where $\beta_0 = \half$ and $\beta_1=
-\half$, and $n_0 + n_1 =0$, consider
\[
\sum^1_{j=0} \lf(\Re  \beta_j + n_j \rt)^2 = \lf(\half + n_0 \rt)^2 + \lf(-\half-n_0 \rt)^2  = 2 \lf(\half + n_0\rt)^2
\]
which clearly achieves its minimum at $n_0=0=-n_1$ and also at $n_0=-1=-n_1$.  These values correspond to
$\wbe = \beta=\lf(\half, -\half\rt)$ and $\wbe = \lf(\half -1, -\half + 1\rt)= \lf(-\half, \half\rt)$,
respectively.  Thus the leading asymptotics of $D_n(f)$ is a sum of contributions from $f(z) = f(z;\beta)$ and
$f(z) = f\lf(z;\wbe\rt)$, as in \eqref{eq95}.

Conjecture \ref{conject8} was  proved recently in \cite{DIK1}.

\begin{theorem}\label{BT} \cite{DIK1}
Conjecture \ref{conject8} holds.
\end{theorem}

The relation of the semi-norm $|||\beta|||$ to
the minimization problem is given by the following elementary result.

\begin{lemma}\label{lem9}\cite{DIK1}
Given $\beta$, there are only two, mutually exclusive, possibilities: Either
\begin{enumerate}
\item[(i)]  There exists $\wbe \in O_\beta$ such that $|||\wbe ||| <1$.  Then such a $\wbe$
is unique and it is the unique element of $\CM_\beta=\lf\{\wbe\rt\}$
\end{enumerate}
or
\begin{enumerate}
\item[(ii)] There exists $\wbe \in O_\beta$ such that $|||\wbe|||=1$.  Then there are at least
two such $\wbe$'s and all of them are obtained from each other by a repeated application of the following
rule:  add 1 to a $\wbe_j$ with the smallest real part, $\Re\wbe_j = \min_k \, \Re  \wbe_k$,
and subtract 1 from a $\wbe_{j'}$ with the largest real part, $\Re\wbe_{j'} = \max_k \:\Re  \wbe_k$.
Moreover $\CM_\beta=\lf\{ \wbe \in O_\beta: |||\wbe ||| =1 \rt\}$.
\end{enumerate}
\end{lemma}

Clearly for $f\!\lf(\eit\rt)$ in \eqref{eq94}, $|||\wbe|||=1$, if and only if $\wbe=\beta$ or $\wbe=\beta
+(-1, 1)$.  Also we see from (i) that the case $|||\beta ||| <1$ in Ehrhardt's Theorem corresponds
to a unique minimizer in \eqref{eq99}, and hence to a single term in \eqref{eq100}.

The proof of Conjecture \ref{conject8} in \cite{DIK1} proceeds as follows.  Given $\beta$, if there exists
$\wbe \in O_\beta$ with $|||\wbe|||<1$, then $\wbe$ is unique and as $n\to\infty$, $D_n(f) = R_n \lf(f(\wbe)\rt)
(1+o(1))$  (Ehrhardt's Theorem).  On the other hand, if there exists $\wbe \in  O_\beta$ with $|||\wbe|||=1$,
then $\CM_\beta = \lf\{\wbe  \in O_\beta: ||| \wbe||| =1\rt\}$ has at least two elements.
Then for some $b$,
\be\label{eq101}
b \le \Re  \wbe_j \le b+1
\ee
for all $\wbe_j$, and for some $p,\, \ell >0$ there are $p$ values of $j$ such that $\Re\wbe_j=b$
and $\ell$ values for which $\Re \wbe_j = b+1$.  Let $\tf \!\lf(\eit\rt)$ be the FH symbol
(\tit{not} a FH-representation of $f\!\lf(\eit\rt)$) with $\al_j$, $V$ as before but $\wtbe_j=
\wbe_j$ if $\Re \wbe_j < b+1$ and $\wtbe_j = \wbe_j-1$ if $\Re \wbe_j = b+1$.
Then clearly
\be\label{eq102}
f(z) = c\, z^\ell\, \tf(z)
\ee
for some simple constant $c$, and by a general relation
\be\label{eq103}
D_n(f) = \frac{c^n (-1)^{\ell n} \, F_n\, D_n (\tf)}{\prod^{\ell-1}_{j=0} j!}
\ee
where
\be\label{eq104}
F_n = \det \lf(\frac{d^{j}\widetilde{\pi}_{n+k}}{dz^j}\,(0) \rt)_{0\le j, k\le \ell-1}\
\ee
and $\wpi_q(z) = z^q + \dots,  q \ge 0$, are the monic orthogonal polynomials with respect to
the (in  general non-real) weight $\tf\!\lf(\eit\rt) \frac{d\tht}{2\pi}$ on $S^1$,
\be\label{eq105}
\int^{2\pi}_0 \wpi_q \lf(\eit\rt) e^{-ir\tht} \: \tf\!\lf(\eit\rt) \frac{d\tht}{2\pi} =0,
\qquad 0 \le r< q.
\ee
(A part of the proof of Conjecture \ref{conject8} consists in showing that the $\wpi_q$'s exist
and are unique for $q$ sufficiently large, $q\ge q_{\tf}$.)  Now clearly $|||\wtbe|||<1$, and
so we can evaluate $D_n(\tf)$ as $n\to\infty$ using Ehrhardt's Theorem.  The proof of the conjecture
then reduces by \eqref{eq103} \eqref{eq104} to evaluating the polynomials $\wpi_n(z),  \dots,
\wpi_{n+\ell-1}(z)$ and their derivatives at $z=0$, as $n\to\infty$.  This is effected by using
the steepest-descent method of Deift and Zhou for Riemann-Hilbert problems mentioned in Section 
\ref{secSzego} above
(see also, for example, \cite{DKMVZ1}, \cite{DKMVZ2}).  The fact that orthogonal polynomials
with  respect to a weight on the line,  can be expressed in terms of a Riemann-Hilbert problem is due to
Fokas, Its and Kitaev in \cite{FIK}: the extension of this result to orthogonal polynomials on the circle
is given in \cite{BDJ}.

\begin{remark}\label{degen}
As noted above, for the Ising model, at $T>T_c$, $\varphi_{\tons}\!\lf(\eit\rt)$ in \eqref{106} corresponds
to a degenerate case for Ehrhardt's Theorem.  However $\varphi_{\tons}\!\lf(\eit\rt)=-e^{-i\tht}
\;\tvph\!\lf(\eit\rt)$ where $\tvph\!\lf(\eit\rt)= \frac{\lf(1-\gamma_1\, \eit\rt) \lf(\gamma_2-\eit\rt)}{
\lf|1-\gamma_1\,\eit\rt| \lf|1-\gamma_2 \, \eit\rt|}$ is smooth, nonzero, and has no winding on $S^1$.  Hence $D_n(\tvph)$
can be computed as $n\to\infty$ using SSLT as in \eqref{eq47}.  Analogous to \eqref{eq102}, \eqref{eq103},
we have the formula
\be\label{eq106}
D_n (\varphi_{\tons}) = \whpi_n(0) \, D_n(\tvph)
\ee
where $\whpi_n(z)=z^n + \dots, n\ge 0$ are the complementary orthogonal polynomials to the $\wpi_n$'s
(see \cite{DIK1}) 
\be\label{eq107}
\int^{2\pi}_0 \whpi_n \!\lf(e^{-i\tht}\rt) e^{i\,j\tht}\: \tvph\!\lf(\eit\rt) \frac{d\tht}{2\pi} =0,
\qquad 0 \le j < n.
\ee
Applying the steepest descent method to the Riemann-Hilbert Problem for $\whpi_n$, and taking into
account the analyticity of $\tvph(z)$ in an annulus around $\{|z|=1\}$,  we obtain the precise decay rate of
$\whpi_n(0)$ as $n\to\infty$, and hence via \eqref{eq106}, the precise leading order exponential decay rate
$c n^{-1/2}\gamma^{-n}_2 \lf(1+o(1)\rt)$ as in \eqref{eq64}.  The steepest descent
method also yields the higher order terms in the expression, to any desired order.
An alternative, and very elegant, proof of \eqref{eq64} using a contour deformation argument, is given in \cite{FF}.
For $\varphi_{\rm diag}(e^{i\tht})$ and the diagonal correlation function, the situation is similar, and one obtains:
\be
\lf\langle\sg_{1,1}\,\sg_{1+n, 1+n}\rt\rangle=
\frac{k_{\rm ons}^{-n}}{\sqrt{n}}
\frac{\sqrt{\pi}}{(1-k_{\rm ons}^{-2})^{1/4}}(1+o(1)),\qquad T>T_c.
\ee
We also note that in Section 10.7 of \cite{BottSilb3} and in the papers \cite{BottWid1} \cite{Wid7}, 
the authors also analyze the Toeplitz determinants arising in the case $T>T_c$.
\end{remark}

\begin{remark}\label{remark6}
Formulae \eqref{eq103} and \eqref{eq106} should be compared respectively with formulae \eqref{eq61} and 
\eqref{eq52} in \cite{FisHart1}. 
\end{remark}

\section{Hankel and Toeplitz+Hankel determinants}
In \cite{DIK1} the authors also consider Hankel and Toeplitz + Hankel determinants with FH-type singularities.
For a weight $w=w(x) \ge 0 $ on $[-1,1]$, the associated Hankel determinant is given by \cite{Sz5}
\be\label{eq108}
D_n^{(H)}(w) = \det \lf(\int^1_{-1} x^{j+k}\, w(x) \, dx\rt)^{n-1}_{j, k=0} \ .
\ee
The Hankel determinant $D_n^{(H)}(w)$ is related to the polynomials orthogonal with respect to the weight $w$
on the interval $[-1,1]$ in a similar way as a Toeplitz determinant is related to OPUC's.   

For fixed $r=0,1,2, \dots$ we consider $w(x)$ of the following form:
\be\label{eq109}
w(x) =e^{U(x)} \prod^{r+1}_{j=0} \lf|x-\lb_j\rt|^{2\al_i}\, \omega_j(x)
\ee
where
\begin{align*}
&1=\lb_0 > \lb_1 > \dots > \lb_{r+1} =-1, \\
&\omega_j(x) =  \lf\{
\begin{array}{lll}
e^{i\pi\,\beta_j}, &\qquad  -1\le x\le \lb_j, &\qquad \Re \beta_j \in \lf(-\half, \half\rt] \\
e^{-i\pi\,\beta_j}, &\qquad  1 \ge x > \lb_j
\end{array}\rt. \\
&\qquad \beta_0 = \beta_{r+1} =0, \qquad\qquad  \Re \al_j >-\half, \qquad j=0,1, \dots, r+1
\end{align*}
and where $U(x)$ is a sufficiently smooth function on $[-1,1]$.

\begin{remark}\label{remark7}
Note that there is no loss of generality in setting $\beta_0 = \beta_{r+1}=0$ and $\Re\beta_j \in \lf(-\half,
\half\rt]$ as the functions $\omega_0(x), \omega_{r+1}(x)$ 
are just constant on $(-1,1)$ and $\omega_j \lf(x; \beta_j + k_j\rt)
= (-1)^{k_j} \, \omega_j \lf(x; \beta_j\rt)$ for $k_j \in \ZZ$.
\end{remark}

In \cite{DIK1} the authors prove that for $w$ as in \eqref{eq109} with $\Re\beta_j \in \lf(-\half, \half\rt),
\:j=1, \dots, r$, as $n\to\infty$,
\be\label{eq110}
D_n^{(H)}(w) =2^{-n^2} \, G^n_H \, n^{-\frac{1}{4}+2 \lf(\al^2_0 + \al^2_r\rt) + \sum^n_{j=1} \lf(\al^2_j-\beta^2_j\rt)}
\, E_H \lf(1+o(1)\rt)
\ee
for certain explicit constants $G_H$ and $E_H$ depending on $U(x)$, $\lb_j$ and $\beta_j$ (see (1.37) in
\cite{DIK1}).
This result is derived by introducing the function
\be\label{eq111}
f\!\lf(\eit\rt) = f\!\lf(e^{i(2\pi-\tht)}\rt) \equiv w(\cos \tht) \, |\sin \tht|, \qquad\qquad 0\le \tht < 2\pi
\ee
on $S^1$.  With $w$ as in \eqref{eq109}, $f$ has FH singularities of the form
\be\label{eq112}
\begin{aligned}
f\!\lf(\eit\rt) &= c\, e^{V\lf(\eit\rt)} \lf|z-1\rt|^{4\al_0 +1} \lf|z+1\rt|^{4\al_{r+1} +1}
\, \prod^r_{j=1} \lf|z-z_j\rt|^{2\al_j}  \lf|z-z'_j\rt|^{2\al_j} \\
&\qquad\times \prod^r_{j=1} \lf[g_{z_j, -\beta_j} (z) \, z^{\beta_j}_j \; g_{z'_j,\, \beta_j} (z)\lf(z'_j\rt)^{-\beta_j}\rt],
\qquad 0 \le \tht <2\pi
\end{aligned}
\ee
where
$$
\lf\{\begin{aligned}
\cos\,\tht_j &=\lb_j,\quad z_j = e^{i\tht_j},\quad j=0,1,\dots, r+1,\quad 0=\tht_0< \tht_1 < \dots < \tht_{r+1} =\pi\\
z'_j &= e^{i\lf(2\pi -\tht_j\rt)}, \qquad j=0,1,\dots, r+1\\
c &= 2^{-2\lf(\sum^{r+1}_{j=0} \,\al_j\rt)-1} \; e^{2i \sum^r_{j=1}\, \beta_j \, \arcsin \lb_j}
\end{aligned}
\rt.
$$
and
\be\label{eq113}
V\!\lf(\eit\rt) = U(\cos \tht).
\ee
The Hankel determinant $D_n^{(H)}(w)$ and the Toeplitz determinant $D_{2n}(f)$ are 
related in the following way:
\be\label{eq114}
D_n^{(H)}(w)^2 = \frac{\pi^{2n}}{4^{(n-1)^2}} \; \frac{\lf(\pi_{2n} (0) + 1 \rt)^2}{\pi_{2n} (1) \, \pi_{2n}(-1)}
\; D_{2n} (f)
\ee
where $\pi_k(z)=z^k + \dots, k\ge 0$, are the monic orthogonal polynomials with respect to the weight
$f\!\lf(\eit\rt) \frac{d\tht}{2\pi}$ on $S^1$ as in \eqref{eq105},
$$
\int^{2\pi}_0 \pi_k \lf(\eit\rt) e^{-ij\tht} \, f\!\lf(\eit\rt) \frac{d\tht}{2\pi} =0,
\qquad  0 \le j < k.
$$
As the semi-norm $|||\beta||| <1$ for $f\!\lf(\eit\rt)$, we may use Ehrhardt's result to evaluate $D_{2n}(f)$,
$n\to\infty$, and the result \eqref{eq110} follows from \eqref{eq114} by computing the asymptotics for $\pi_{2n}
(z)$ for $z=0, \pm 1$.  Note that if $\Re \beta_j = \half$ for some $j \in \{0, 1, \dots, r+1\}$, then
$\lf|\Re \beta_j - \Re \lf(-\beta_j\rt)\rt|=1$ and so $||| \beta |||=1$.  In this case we would have had to use
the Basor-Tracy form \eqref{eq100} for $D_{2n}(f)$: we have restricted our attention to the case
$\Re \beta_j \in \lf(-\half, \half\rt)$ only for simplicity.

\begin{remark}\label{remark8}
For a discussion of results of various authors, particularly Basor and Ehrhardt, related to \eqref{eq110},
see \cite{DIK1}, Remark 1.24.
\end{remark}

For even FH symbols $f$, i.e., $f\!\lf(\eit\rt) = f\!\lf(e^{i(2\pi-\tht)}\rt)$, the
authors in \cite{DIK1} also consider Toeplitz + Hankel determinants of four types defined in terms of the
Fourier coefficients $f_k$ of $f$,
\be\label{eq115} 
\det \lf(f_{j-k} + f_{j+k}\rt)^{n-1}_{j,\, k=0}, \quad\det \lf(f_{j-k}- f_{j+k+2}\rt)^{n-1}_{j,\, k=0},
\quad\det \lf(f_{j-k}\pm f_{j+k+1}\rt)^{n-1}_{j,\,k=0}.
\ee
Such matrices arise, in particular, in  the theory of classical groups and its application to random matrix
theory and statistical mechanics, and there are simple relations between these determinants and Hankel
determinants on $[-1,1]$ with added singularities at the end-points (see \cite{DIK1} and the references
therein).  For example, for $f\!\lf(\eit\rt)=f\lf(e^{i(2\pi-\tht)}\rt)$,
\be\label{eq116}
\det \lf(f_{j-k} + f_{j+k}\rt)^{n-1}_{j,\,k=0}=\frac{2^{n^2-2n +2}}{\pi^n} \: D_n^{(H)}(v)
\ee
where $D_n^{(H)}(v)$ is the Hankel determinant with symbol $v(x) = f\!\lf(e^{i\tht (x)}\rt)/\sqrt{1-x^2}$,
$\tht(x)= \arccos x$, on $[-1,1]$, with similar formulae for the other three determinants.
Let $f\!\lf(\eit\rt) = f\lf(e^{i(2\pi-\tht)}\rt)$ with FH singularities at $0=\tht_0 < \tht_1 < \dots
< \tht_r < \tht_{r+1} =\pi$ and complementary singularities at $2\pi-\tht_j$, $1\le j \le r$.  Suppose
that $\Re \beta_j\in (-\half, \half)$, $j=1, \dots, r$, and $\beta_0 = \beta_{r+1}=0$.  Then utilizing
\eqref{eq116} and \eqref{eq110}, the authors show that as $n\to\infty$
\be\label{eq117}
\begin{aligned}
\det \lf(f_{j-k} + f_{j+k}\rt)^{n-1}_{j,\,k=0} &= G^n_{T+H} \,  n^{\sum^{r}_{j=1} \lf(\al^2_j-\beta^2_j\rt)
+\half\lf(\al_0^2+\al_{r+1}^2-\al_0-\al_{r+1}\rt)} \\
& \quad \times E_{T+H} \lf(1+o(1)\rt)
\end{aligned}
\ee
for certain explicit constants $G_{T+H}$, $E_{T+H}$, with similar results for the other three Toeplitz+Hankel
determinants in \eqref{eq115}.

\begin{remark}\label{remark9}
For a discussion of the related results of Basor and Ehrhardt for determinants of type \eqref{eq115}, see
\cite{DIK1}, Remark 1.27.
\end{remark}

\begin{remark}
In \cite{FF} the authors discuss a large number of conjectures and results for Toeplitz, Hankel, and Toeplitz+Hankel determinants with FH-type singularities. The determinants arise in turn from problems in statistical physics. Many of the conjectures in \cite{FF} can now be resolved using results from \cite{DIK1} such as \eqref{eq110} and \eqref{eq117}.
\end{remark}

\section{Some applications}\label{sec9}
As advertised at the outset, the goal of this paper has been to show how SSLT developed and was generalized in
response to questions arising in the analysis of the Ising model.  Although the case $|||\beta|||=1$ in
\cite{DIK1} does not arise in the Ising model (the same is true for  the Hankel and Toeplitz + Hankel
determinants \eqref{eq108} and \eqref{eq115} respectively), this case does arise in many other problems
in statistical physics.  For example, the probability $P_E(n)$ of a string of length $n$
of ferromagnetically aligned spins
in the antiferromagnetic ground state in the $XY$ spin chain, is given, for a certain range of parameters,
by a Toeplitz determinant with $2$ $\beta$-singularities such that $|||\beta|||=1$.  Thus \eqref{eq100}
verifies the result on $P_E(n)$, $n\to\infty$, presented in \cite{FrAb}, which the authors based on the
Basor-Tracy conjecture.  In a similar way, the results in \cite{DIK1} can be used to justify the asymptotic
results for correlations arising in the theory of the impenetrable Bose gas that were obtained in the
work of Ovchinnikov \cite{Ov} on the basis of the Basor-Tracy conjecture.
In another direction in \cite{GGM}, the authors use \eqref{eq100} to evaluate the long-time behavior of a
variety of non-equilibrium 1D many-body problems.  In particular, in the Fermi edge singularity and
tunneling spectroscopy problems considered in \cite{GGM}, the authors show that the various contributions
to \eqref{eq100} have a very transparent physical meaning:  they correspond to contributions to the
Green's functions for these problems from multiple Fermi edges.

Theorem \ref{BT} can also be used to analyze the asymptotic behavior of the eigenvalues
$\lb^{(n)}_1 \le \lb^{(n)}_2 \le \dots \le \lb^{(n)}_n$, $n\to\infty$, of a Toeplitz matrix
$T_n(\varphi)$ with a real symbol $\varphi$.  At the beginning of our story (see Theorem \ref{theorem2})
we interpreted Szeg\H{o}'s early results as saying that the $\lb^{(n)}_j$'s were equidistributed as $n\to\infty$.
But what about the behavior of individual eigenvalues, $\lb^{(n)}_k$, $n\to\infty$?  Early work on this
question considered the so-called extreme eigenvalues $\lb^{(n)}_k$ and $\lb^{(n)}_{n-k}$ with $k$
fixed as $n\to\infty$ (see \cite{GreSz}, \cite{BGrM} and the references therein).  The important case for
$\lb^{(n)}_k$ with $n\to\infty$ and $k/n\to x \in (0,1)$ was only considered for the first time recently
in \cite{BGrM}.
Recall that by general principles (see e.g. \cite{GreSz}), for bounded, real-valued functions $f\!\lf(\eit\rt)$, the
eigenvalues $\{\lb^{(n)}_k\}$ lie in the interval $[L, M]$, $L=\inf \, f\!\lf(\eit\rt)$, $M=\sup
\, f\!\lf(\eit\rt)$, and as $n\to\infty$ the $\lb^{(n)}_k$'s fill out the interval in the sense that if
$\lb\in [L, M]$, then there exist $k_n(\lb)$ such that $\lb^{(n)}_{k_n(\lb)} \to \lb$.
In \cite{BGrM} the authors consider function $f\!\lf(\eit\rt)$ that are unimodal, real-valued trigonometric polynomials,
i.e., for some $0<\tht_0 < 2\pi$, $\frac{d}{d\tht} \: f\!\lf(\eit\rt)>0$ for $0< \tht <\tht_0$
and $\frac{d}{d\tht} \:f\!\lf(\eit\rt) <0$ for $\tht_0 < \tht < 2\pi$, with the additional property that
$\frac{d^2}{d\tht^2} \:f\!\lf(\eit\rt) \neq 0$ at $\tht=0$ and $\tht=\tht_0$.  For any $\lb \in (L,M)$,
there exist unique
$0 < \tht_1(\lb) < \tht_0 < \tht_2(\lb) < 2\pi$ such that $f\!\lf(e^{i\tht_j(\lb)}\rt) - \lb =0, \;
j=1,2$.  The function
$$
\psi(\lb) = \half \lf[\tht_1(\lb) - \tht_2(\lb)\rt] + \pi
$$
extends to a homeomorphism from $[L,M]$ onto $[0,\pi]$.  In \cite{BGrM} the authors prove, in particular,
that for any $x\in (0,1)$
\be\label{eq118}
\lb^{(n)}_k = \lb_x + o(1)
\ee
as $n\to\infty$ and $k/n \to x$, where $\psi (\lb_x) = \pi x$, and
\be\label{eq119}
\lb^{(n)}_{k+1} - \lb^{(n)}_k = \frac{\pi}{\psi'(\lb_x)} \:\frac{1}{n+1} + o\lf(\frac{1}{n}\rt)
\ee
if $\lf| k/n - x\rt|= O\lf(\frac{1}{n}\rt)$ as $n\to\infty$.

The estimate \eqref{eq119} brings precision to the equidistribution result in Theorem \ref{theorem2}.  In
\cite{DIK2} the authors extend the result in \cite{BGrM} to $C^\infty$ unimodal, real-valued functions on $S^1$,
and they also consider piecewise--constant symbols of the form
\be\label{eq120}
\begin{aligned}
f\!\lf(\eit\rt) &= 1, \qquad & \tht &\in \lf[0, \wht_1 \rt) \cup \lf[\wht_2, 2\pi\rt) \\
&= e^{2\pi\gamma}, \qquad &\tht&\in \lf[ \wht_1, \wht_2\rt)
\end{aligned}
\ee
where $0\le \wht_1 < \wht_2 < 2\pi$ and $\gamma >0$ are given.  For symbols such as
\eqref{eq120}, the eigenvalues $\lb^{(n)}_k$ cannot be spaced $1/n$ apart as in \eqref{eq119}.
This is clear, for example, from \eqref{eq4} in Theorem \ref{theorem2}, if we let
$F(\lb)$ be any non-negative continuous function with support in
$\lf(1, e^{2\pi\gamma}\rt)$.   And indeed, the general result in \cite{Bas2} asserts that  in the case
\eqref{eq120}, for any small $\ep >0$, the interval $\lf(1+\ep, e^{2\pi\gamma}-\ep\rt)$ contains
order $\log n$ eigenvalues $\lb^{(n)}_k$, spaced at distance $1/\log n$ apart, as $n\to\infty$: The bulk
of the eigenvalues accumulate at $1$ and at $e^{2\pi\gamma}$.  On the basis of numerical computations,
the authors in \cite{LeSh} (see also \cite{LeSoSo}) conjectured that in the case
\be
\wht_2 - \wht_1 = 2\pi \, p/q, \qquad p,\, q\in \ZZ, \qquad 0 < p < q
\ee
there is a ``near periodicity'' in the spectrum of $T_n(f)$.  In \cite{DIK2} the authors prove this
conjecture in the following sense: for $n$ sufficiently large, if $\lb^{(n)}_k$ in an eigenvalue of $T_n(f)$
in $\lf(1+\ep,\, e^{2\pi\,\gamma}-\ep \rt)$, then there is an eigenvalue $\lb^{(n+q)}_j$ of $T_{n+q}$
such that $\lf|\lb^{(n)}_k - \lb^{(n+q)}_j \rt| \le c (n\log n)^{-1}$.  Of course
$1/(n \, \log  n) << 1/\log n$, the spacing of the $\lb^{(n)}_k$'s, $k=1, \dots, n$.

The proof of the above results in \cite{DIK2} rests on the simple fact that $\det \lf(T_n(f) - \lb\rt)=
\det \lf(T_n (f-\lb)\rt)$, so that the eigenvalue equation for $T_n(f)$ is equivalent to the vanishing
of the Toeplitz determinant $D_n(f-\lb)$.  Now if $f\!\lf(\eit\rt)$ is a smooth, unimodal function on
$S^1$, then $f\!\lf(\eit\rt)-\lb$ is a FH symbol with two singularities at $z_1=e^{i\tht_1 (\lb)},
\; z_2=e^{i\tht_2(\lb)}$ and $\al_1 = \al_2 =\half, \; \beta_1=-\beta_2 =\half$ (see \eqref{eqBG})
\be\label{eq121}
f(z) -\lb =e^{V(z)} \lf|z-z_1\rt| \lf|z-z_2\rt| g_{z_1,\, \half} (z) \:g_{z_2,\, -\half}(z)
\lf(\frac{z_1}{z_2}\rt)^{-\half}
\ee
where $V\!\lf(\eit\rt)$ is smooth on $S^1$.  
This is a case of a FH symbol with $|||\beta|||=1$. 
Therefore we need to use Theorem \ref{BT} to obtain the asymptotics for $D_n(f-\lb)$.  
The set $M_{\beta}$ here contains 2 elements, and the condition for an eigenvalue
is the vanishing of the sum in \eqref{eq100}:
\be
D_n(f-\lb)=R_n(f(\beta_1,\beta_2)-\lb)(1+o(1))+R_n(f(\beta_1-1,\beta_2+1)-\lb)(1+o(1))=0
\ee
The asymptotic results on the eigenvalues follow from this formula. 


The Toeplitz matrix $T_n(f)$ with the symbol $f\!\lf(\eit\rt)$ in \eqref{eq120} is handled similarly.
The function $f(z)-\lambda$ in this case has 2 $\beta$-type singularities whose positions on the circle are fixed,
but the imaginary parts of $\beta$'s depend on $\lambda$:
$\beta_1=-\beta_2=\frac{1}{2}+i\gamma(\lambda)$. Note that we again have $|||\beta|||=1$, and the condition for 
the eigenvalues is thus given by Theorem \ref{BT}.

Finally we describe an example where the results in \cite{DIK1} are needed for a Toeplitz+Hankel determinant
of type \eqref{eq115}.  The following problem arises in the framework of the random matrix approach to the
theory of the Riemann zeta-function and other $L$-functions (see \cite{Keat} for a recent survey).
Define
\be\label{eq123}
\phi(\eit) = \lf|2 \sin \frac{\tht}{2}\rt|^{2k} \, e^{V(\eit)}, \qquad k\in \NN
\ee
where
\begin{align*}
V\!\lf(\eit\rt) &= 2k \lf[\int^e_1 u(y)
\lf\{\sum^\infty_{j=-\infty} \mathrm{Ci}\, \lf(  \lf|\tht+2\pi j\rt|
(\log y) (\log X)\rt)\rt\}dy- \log \lf|2\sin \frac{\tht}{2}\rt| \rt], \\
&\quad \mathrm{Ci}\,(z) =- \int^\infty_z \frac{\cos t}{t} \; dt
\end{align*}
and $u(y)$ is a smooth, non-negative function supported on  $\lf[e^{1-X^{-1}}, e\rt]$ and of total mass
one.  Consider the following average over the orthogonal group $SO(2n)$,
\be
\EE_{SO(2n)} \lf(\prod^n_{j=1} \phi \lf(e^{i\tht_i}\rt)\rt).
\ee
Here $e^{\pm i\tht_1}, \dots e^{\pm i\tht_n}$ are the eigenvalues of a random matrix in $SO(2n)$: note that
$\phi$ is even, $\phi\!\lf(\eit\rt)=\phi\!\lf(e^{-i\tht}\rt)$.  This expectation was introduced by Bui and Keating
in \cite{BuKe} (following \cite{GHK}) as the random matrix counterpart of a key term contributing to
the mean values of certain Dirichlet $L$-functions in the Katz-Sarnak orthogonal family.  The issue
is the large $n$ and large $X$ behavior of this average.  The question can be resolved with the help of
\eqref{eq117}.  One observes that
\be
\EE_{SO(2n)} \lf(\prod^n_{j=1} \phi \lf(e^{i\tht_j}\rt)\rt)=\half \,\det \lf( \phi_{j-k}+
\phi_{j+k}\rt)^{n-1}_{j,\,k=0}
\ee
and that the symbol \eqref{eq123} is of FH type with a single  $\al$-singularity at $z_0=1$ and
$\al_0=k$.  Directly from (the precise version of) the asymptotic formula \eqref{eq117} (see Theorem
1.25 in \cite{DIK1}), one obtains, for $X$ large, as $n\to\infty$
\be\label{eq124}
\EE_{SO(2n)} \lf(\prod^n_{j=1} \phi\lf(e^{i\tht_j}\rt)\rt) \sim G(1+k)
\lf( \frac{\Gamma (1+2k)}{G(1+2k) \;\Gamma(1+k)} \rt)^{\!\half} \lf(\frac{2n}{e^{\gamma_E} \,\log X}
\rt)^{\!\frac{k(k-1)}{2}}
\ee
where $\gamma_E$ is Euler's constant.  Formula \eqref{eq124} is precisely the asymptotic form conjectured
by Bui and Keating in \cite{BuKe}.

\section{Block Toeplitz matrices}
We now consider block Toeplitz matrices $D_n(\varphi)$ of the form $D_n(\varphi) =
\det\!\lf(\varphi_{j-k}\rt)^{n-1}_{j, k=0}$ where $\varphi_q \in M(r,\CC)$ are the Fourier
coefficients of an $r\times r$-matrix-valued integrable function $\varphi\!\lf(\eit\rt)$ on $S^1$.
The earliest asymptotic result for such determinants is the analog of Szeg\H{o}'s Theorem 1 due to
Gyires in 1956 \cite{Gy}, who showed that if $\varphi(z)$ is continuous and positive definite for $|z|=1$, then
\be\label{eq125}
\lim_{n\to\infty} \,\frac{1}{n}\, \log D_n(\varphi) = \frac{1}{2\pi}
\int^{2\pi}_0   \log\det \, \varphi\!\lf(\eit\rt)d\tht.
\ee
Block Toeplitz matrices arose in the Ising model in the following way.  In their famous paper
\cite{LeYa} on the Lee-Yang Theorem, asserting that the zeros of the partition function $Z_{\Lambda,h}$
(cf.\ Remark \ref{remark2} in Section \ref{section4}), 
for the Ising model in the presence of a magnetic field, all lie on the imaginary
$h$-axis, Lee and Yang also proposed, but did not prove, expressions for the free energy  $F_h=
\lim_{|\Lambda| \to\infty} - k_B\, T\, \log Z_{\Lambda, h}$ and the magnetization 
$M_h= \lim_{n\to\infty} \lf\langle \sg_{1,1}\, \sg_{1, 1+n}\rt\rangle^{\half}$ at the
particular imaginary value of $h$, $\frac{h}{k_B\,T}= \half \, i\, \pi$.  In terms of the variable
$\rho =e^{-2h/k_B\,T}$, the zeros of $Z_{\Lambda, h}$ lie on the circle $|\rho|=1$.  As noted in
\cite{LeYa}, $M_h$ is precisely the density of the zeros of $Z_{\Lambda, h}$ on $\{|\rho|=1\}$ at
$\rho=-1$, in the limit as $|\Lambda|\to\infty$.  In \cite{McWu2}, McCoy and Wu present a derivation
of these expressions for $F_h$ and $M_h$, $h/k_B\,T=\half\, i\pi$.  Regarding $M_h$, for this value of
$h$, they first show how to express the correlation function $\lf\langle \sg_{1,1}\,\sg_{1, 1+n}\rt\rangle$
in terms of a $2\times 2$ block Toeplitz determinant $B_n$.  Then they show by certain ingenious manipulations
how to evaluate $B_n$ in terms of a determinant which is a product of two scalar Toeplitz determinants to
which SSLT applies.  In this way they are able to compute $M_h$ explicitly and recover the formula in \cite{LeYa}.  
This is the first example of an SSLT-type result for block Toeplitz matrices.

The general theory of SSLT for block Toeplitz matrices begins with Widom's paper \cite{Wid5}.
Under certain smoothness assumptions on
the $r\times r$-matrix-valued function $\varphi\!\lf(\eit\rt)$, Widom shows that
\be\label{eq126}
Y(\varphi) \equiv \lim_{n\to\infty} D_n(\varphi)/\exp\lf\{\frac{n}{2\pi} \int^{2\pi}_0 \,\log  \det
\,\varphi \!\lf(\eit\rt)d\tht\rt\}
\ee
exists, and in \cite{Wid3} he shows that
\be\label{eq127}
Y(\varphi) =\det \lf[T(\varphi) \, T\bigl(\varphi^{-1}\bigr)\rt]
\ee
where $T(\varphi)$, $T\bigl(\varphi^{-1}\bigr)$ are the Toeplitz operators associated with
$\varphi$ and $\varphi^{-1}$ respectively (cf.~\eqref{eq46-3}).   Part of the proof of \eqref{eq127}
consists in showing that $T(\varphi)\: T(\varphi^{-1})-1$ is trace-class in $\ell^2_+$,
so that the RHS is well-defined.   In the scalar case, as mentioned above, Widom uses the Helton-Howe-Pincus formula to show that \eqref{eq127} reduces to the standard result $Y(\varphi)=
e^{\sum^\infty_{k=1} \, k(\log\varphi)_k \, (\log\varphi)_{-k}}$.  As noted above (see \eqref{eq46}
et seq), in \cite{Wid4}, Widom introduced a ``standard'' pathway to SSLT.  The method extends with
only a change of notation to the block Toeplitz case, and so Widom obtains a new direct proof of
\eqref{eq126}, \eqref{eq127} with no more effort than in the scalar case. The same is true for the method of Basor and Helton in \cite{BasHelt}.

Block Toeplitz matrices occur in many physical situations, for example, in the study of the Ising
model with next-nearest-neighbor interactions,  and also, more recently, in the theory of the
dimer model  \cite{BasEhr} and  in the study of the entanglement spectrum in quantum systems (see  \cite{VLRK},
\cite{JK}). 
The use  of the general SSLT for block Toeplitz matrices in physical
applications  poses new challenges vis a vis the scalar case. The new difficulties are due to
the fact that formula (\ref{eq127}), as elegant as it is, is hard to evaluate in concrete calculations in the case
of a matrix symbol $\varphi$. The determinant on the right hand side of (\ref{eq127}) is the Fredholm determinant of an infinite matrix, and even for  $2 \times 2$ matrix functions $\varphi$ an
effective evaluation  of  $Y(\varphi)$   is  a highly nontrivial enterprise. Indeed, the authors  are aware of only three situations when such effective evaluations have been made so far. The first situation arises  when, as in the work of McCoy and Wu mentioned above, it is possible, by using an ad hoc technique,  to express the  block Toeplitz  determinant in question in terms of scalar Toeplitz determinants. The work of  Tanaka, Morita, and Hiroike 
\cite{TMH}, and of B\"ottcher \cite{Bott2} on  the Ising model with the next-nearest-neighbor interactions provides another example of this sort. In this example, the relevant matrix 
symbol $\varphi$ is described by the formula 
\begin{equation}\label{eqBott1}
\varphi(z) = \left(I - zR\right)a(z)\left(I - z^{-1}S\right)
\end{equation}
where $a(z)$ is a diagonal matrix function,
\begin{equation}\label{eqBott2}
a(z) = \begin{pmatrix}\lambda(z)&0\cr0&\lambda(z^{-1})\end{pmatrix}
\end{equation}
with $\lambda(z)$ being an elementary algebraic function, 
and the constant matrices $R$ and $S$ are expressed as follows:
$$
R= U\begin{pmatrix}\alpha&0\cr0&\beta\end{pmatrix}U^{-1},
\quad 
S= U\begin{pmatrix}\beta&0\cr0&\alpha\end{pmatrix}U^{-1}.
$$
In these expressions, $U$ is an invertible $2\times 2$ matrix, 
and $\alpha$, $\beta$ are complex numbers, satisfying the conditions: $|\alpha| < 1$,
$|\beta| < 1$. Note that these inequalities guarantee that the first
matrix factor in \eqref{eqBott1} is holomorphic and invertible inside
the unit disc while the last factor is holomorphic and 
invertible outside of the unit disc. The function $\lambda(z)$ (an analog of 
Onsager's symbol \eqref{eq21}), the matrix $U$, and the numbers
$\alpha$ and $\beta$ are determined by the physical parameters of the problem,
the temperature $T$ and the anisotropy  coefficients. As in the usual Ising model,
there exists a critical temperature $T_c$, and the behavior of the determinant 
$D_n(\varphi)$ depends on whether $T > T_c$ or $T<T_c$. B\"ottcher shows that in the case $T>T_c$ (the paramagnetic phase) the Wiener-Hopf factorization of the matrix symbol
$\varphi^{-1}(z)$ has nontrivial partial indices which implies that the  operator $T(\varphi^{-1})$ is not
invertible. In view of \eqref{eq127}, this implies that  the factor  $Y(\varphi)$ is zero.
Moreover, it is easily seen that $\int_0^{2\pi}\log \det \varphi \Bigl(e^{i\theta}\Bigr)d\theta = 0$,
and hence for $T > T_c$ 
$$
\lim_{n \to \infty}D_n(\varphi) =0
$$
which has the physical interpretation that there is no spontaneous magnetization in the paramagnetic phase.
In the ferromagnetic phase, that is when  $T < T_c$,  an efficient  direct analysis 
of Widom's factor $Y(\varphi)$ is not apparent, and B\"ottcher develops  an ad-hoc
approach which uses  the specifics of  the symbol \eqref{eqBott1}. Indeed, in his preprint, by means of ingenious operator-algebra arguments, 
B\"ottcher obtains the formula
\begin{equation}\label{eqBott3}
\lim_{n \to \infty}\frac{D_{n-1}(\varphi)}{D_n(a)} = (\det U)^{-2} (1-\alpha\beta)^{-2}
\left[\frac{U_{11}U_{22}}{\lambda_+(\beta)\lambda_-(1/\alpha)}
-\frac{U_{12}U_{21}}{\lambda_+(\alpha)\lambda_-(1/\beta)}\right]^{2}
\end{equation}
where $\lambda = \lambda_+\lambda_-$ is a Wiener-Hopf factorization of $\lambda$
(which, as B\"ottcher shows,  is possible in the case $T < T_c$). 
Hence the evaluation of $D_{n-1}(\varphi)$ as $n\to\infty$ is reduced to the asymptotic evaluation of the scalar Toeplitz determinant $D_n(\lambda)$, for which the standard (scalar)  SSLT applies.
In the critical case, $T=T_c$, the function $\lambda(z)$ acquires a Fisher-Hartwig
type singularity. B\"ottcher then provides a modification of his previous  analysis for this case 
and, using in addition various results obtained earlier with Silbermann (see Section \ref{sectionFH}), he then 
derives the leading asymptotics of $D_n(\varphi)$ in the critical case as well.

The two other types of  block Toeplitz determinants for which 
the Widom pre-factor $Y(\varphi)$ can be evaluated are not reducible to the scalar case.
The first one corresponds to the matrix symbols $\varphi$ with a Fourier series that is truncated on at 
least one side.
This class was singled out by Widom himself in \cite{Wid5}. For such matrices  $\varphi$,
the evaluation of $Y(\varphi)$ is reduced to the evaluation of a finite dimensional
determinant (in fact, a block Toeplitz determinant of a finite fixed size). 
This result of Widom has been used in the 
recent paper \cite{BasEhr} of  Basor and  Ehrhardt devoted 
to the dimer model. The matrix function  $\varphi$ appearing in the dimer model
is not originally of the truncated form, but it is truncated up to a scalar algebraic
factor. Basor and  Ehrhardt show that this factor can be accounted for 
by means of skillful algebraic manipulations which transform 
the original quantity $Y(\varphi)$ to one with $\varphi$ replaced by its
truncated part. 

The second type of  block Toeplitz determinants that are not reducible to the scalar case, but for which the large $n$ limit is still computable, 
are determinants with algebraic matrix symbols.  
Such determinants appear in the theory of quantum entanglement when one evaluates
the {\it von Neumann  entropy}, which is a fundamental measure of the entanglement of a subsystem
with the rest of the quantum system.  For instance,
as shown in  \cite{VLRK}, \cite{JK}, the  evaluation of the von Neumann  entropy of a subsystem
of neighboring spins in the XY quantum spin chain is  equivalent to the asymptotic evaluation of the block Toepltz determinant $D_n(\varphi)$
whose matrix symbol $\varphi$  is given  by the formula (for more details and
references see  \cite{JK}),
\begin{equation}\label{phient}
\varphi(z) = \begin{pmatrix}i\mu & \psi(z)\cr
               -\psi^{-1}(z)& i\mu
               \end{pmatrix},
\quad \mbox{and}\quad 
 \psi(z)=\sqrt{\frac{(z-z_1)(z-z_2)}{(1-z_{1}z)(1-z_{2}z)}}\;.
\end{equation}
Here,  $z_1 \neq z_2$ are complex nonzero numbers not lying on the unit circle and
determined by the physical data of the model (the anisotropy parameter and the magnetic field) and
$\mu$ is a free complex  parameter lying outside of the interval $[-1, 1]$. The asymptotic evaluation
of $D_n(\varphi)$ is performed in \cite{IJK} and it is based on yet another formula of Widom
which he obtained in \cite{Wid5}. This formula is concerned with block Toeplitz determinants 
depending on an external parameter, say $\mu$, and is used to evaluate the leading asymptotics of the 
$\mu$-logarithmic
derivative of the determinant in terms of the Wiener-Hopf factorizations of the
matrix function $\varphi^{-1}$, i.e. in terms of the $\pm$ - matrix functions $u_{\pm}$ and
$v_{\pm}$  defined by the equations
\begin{equation}\label{wienerblock}
\varphi^{-1}(z) = u_+(z)u_-(z) = v_-(z)v_+(z)
\quad u_-(\infty) = v_-(\infty) = I,
\end{equation} 
where, as usual, $f_+$ ($f_-$) means that the matrix function is invertible and analytic 
inside (outside) of the unit circle, respectively. Under some natural
smoothness assumptions (see \cite{Wid5}, Theorem 4.1), Widom's formula reads:
\be
\begin{aligned}\label{our1}
\frac{d}{d\mu}\log
  D_n(\varphi) &= \frac{n}{2\pi i}\int_{S^1}\tr
  \left(\varphi^{-1}(z)\frac{\partial \varphi(z)}{\partial \mu}\right)
  \frac{dz}{z}\\
& + \frac{i}{2\pi } \int_{S^1}\tr\left(\Bigl(u_{+}'(z)u_{-}(z)
-v_{-}'(z)v_{+}(z)\Bigr)\frac{\partial \varphi(z)}{\partial \mu}\right)dz
+ o(1), \quad n \to \infty.
\end{aligned}
\ee
(In the case of symbols analytic in an annulus including the unit circle,
formula (\ref{our1}) was re-derived with a quantitative estimate
of the error term in \cite{IJK} via the Riemann-Hilbert scheme (see also \cite{IMM}).)
Widom used this asymptotic formula in the proof of his main result, vis., the asymptotic
relation (\ref{eq126}). A curious fact is that, up until the work \cite{IJK},  equation (\ref{our1}),
apparently, had never been used directly for the asymptotic evaluation of
the determinants of block Toeplitz matrices. And for a good reason:
the efficiency of formula (\ref{our1}) relies on the  effectiveness of  
the Wiener-Hopf factorizations (\ref{wienerblock}) of a given matrix valued
function $\varphi(z)$.  Usually, this is a transcendental problem equivalent to the solution of
a system of singular integral equations on the unit circle which can seldom
be effected in terms of known special functions. However, as  was observed 
in \cite{IJK}, in the case of algebraic matrix functions $\varphi(z)$ one
can take advantage of the algebro-geometric method that had been
developed in the late 70s to the early 90s in soliton theory (see e.g. \cite{BBEIM}).
Using the Riemann-Hilbert  version of this method exploited    in \cite{DIZ}
(see also \cite{DKMVZ1}, \cite{DKMVZ2}), the authors of \cite{IJK} were
able to construct the Wiener-Hopf factorizations of the symbol (\ref{phient})
in terms of Jacobi theta-functions. This in turn,  with the help of Widom's formula
(\ref{our1}), yielded the following asymptotics of the corresponding 
Toeplitz determinant: 
\begin{equation}\label{block3}
D_n(\varphi) \sim \frac{ \theta_{3}\left( \frac{1}{2\pi i}\log\frac{\mu+1}{\mu-1}+
\frac{\sigma\tau}{2}\right) \theta_{3}\left(\frac{1}{2\pi i}\log\frac{\mu+1}{\mu-1} -
\frac{\sigma\tau}{2}\right)}{\theta^{2}_{3}\left(\frac{\sigma\tau}{2}\right)}
(1 - \mu^2)^{n}, \quad n \to \infty.
\end{equation}
Here 
$$
\theta_3(s)\equiv \theta(s;\tau) = \sum_{k=-\infty}^{\infty} e^{\pi i \tau k^2+2\pi i s
k}
$$
is the Jacobi theta-function, $\sigma = \pm 1$, and   $\tau$ is the modulus of the elliptic curve
\begin{equation}\label{curve}
w^{2}(z) = (z-z_1)
(z-z_2)(z-z_2^{-1})(z-z_1^{-1}).
\end{equation}
The concrete choice of  $\sigma$  depends on the way the  branch points
$z^{\pm 1}_{1,2}$ are located with respect to the unit circle (see \cite{IJK} for details).
We refer the reader to the
survey \cite{IK} for more on these results and on their applications in
the theory of quantum entanglement.
The asymptotic formula (\ref{block3}) was  extended in \cite{IMM}  to the case of more
general integrable spin chains introduced in \cite{KM}. The relevant symbol has the same matrix structure as 
(\ref{phient}), but with scalar
function $\psi(z)$ now given by 
\begin{equation}
  \label{eq:g_def}
   \psi(z) = \sqrt{\frac{p(z)}{z^{2m}p(1/z)}}
\end{equation}
where  $p(z)$ is a polynomial of degree $2m$.  The analog of the formula (\ref{block3}) 
in the case $m >1$ involves a hyperelliptic curve
instead of an elliptic curve, and instead of the Jacobi
theta-function one needs $2m -1$ dimensional Riemann theta-functions. 
Independently of   \cite{IJK} and  \cite{IMM},  similar results for block Toeplitz determinants with algebraic symbols appearing in the theory of the Gelfand-Dickey hierarchy, 
were  obtained, again with the help of the algebro-geometric
approach, in  \cite{Caf}.  It is worth noticing that the method of \cite{IJK} can be
applied to  the problems considered above by
B\"ottcher and Basor-Ehrhardt. Indeed, the symbol 
of the determinant appearing in the work of  Basor  and Ehrhardt is in fact algebraic 
and hence can be factorized within the algebro-geometric approach, while the symbol  of B\"ottcher's
determinant is already in the form $u_+u_-$. The reverse factorization, i.e. representation of the
symbol (\ref{eqBott1}) as $v_-v_+$ involves an elementary algebraic operation. 
We note finally that Widom's formula (\ref{our1})
can be used to extend B\"ottcher's results to symbols of the form
$\varphi(z) = R(z)a(z)S(z^{-1})$, where $a(z)$ is a diagonal matrix function and $R(z)$ and $S(z)$
are polynomial matrix functions invertible for all $z$ inside the unit circle.

\section{Double-scaling limits}\label{sec-ds}
Of all the challenges in Toeplitz theory that arise from the Ising model, the deepest is the issue
of the double-scaling limit (or in physicists' parlance, simply the scaling limit),
$T\to T_c$ and $n\to\infty$.  We have already mentioned the importance
of this problem in physics.  At the purely mathematical level, the issue is the following:
Suppose one has a Toeplitz determinant $D_n(\varphi)$ with a symbol that depends on some external parameter,
say $t$, $\varphi\!\lf(\eit\rt) = \varphi_t \!\lf(\eit\rt)$.  For $t>0$, $\varphi_t$ is regular, but
at $t=0$, $\varphi_t$ has a (Fisher-Hartwig) singularity, and hence, as we know, the nature of the
asymptotics of $D_n(\varphi_t)$ as $n\to\infty$ is different for $t>0$ and $t=0$.  This raises the
general question: How does one describe the transition from the one kind of asymptotics to the other
as $n\to\infty$ and $t\downarrow 0$? 


In \cite{WMTB}, Wu, McCoy, Tracy, and Barouch discovered a
remarkable scaling relation for the Ising model.  They defined the scaling functions
\be\label{eq70}
G_\pm (r) = \lim_{\ell, m\to \infty, t\to\mp 0} \; \left|1 - e^{-2t}\right|^{-\frac{1}{4}} \lf\langle \sg_{1,1}\,
\sg_{1+\ell, 1+m}\rt\rangle, \quad t = \log \left(\sinh\frac{2J_1}{k_BT}\sinh\frac{2J_2}{k_BT}\right)
\ee
in the double-scaling limit
\be\label{eq71}
\lf(\frac{\sinh \lf(2 \,J_1/k_B\,T_c\rt)\ell^2 + \sinh \lf(2\,J_2/k_B\,T_c\rt)m^2} 
{\sinh\lf(2 \,J_1/k_B\,T_c\rt)+\sinh\lf(2 \,J_2/k_B\,T_c\rt)}
\rt)^{\frac{1}{2}} |t|
\equiv r, \text{ fixed }.
\ee
Here $\pm$ refer to $T> T_c$ and $T< T_c$ respectively.  They then showed that $G_\pm (r)$ could  be
expressed in terms of a solution of the Painlev\'{e} III equation
\be\label{eq72} 
\frac{d^2\eta}{d\tht^2} = \frac{1}{\eta} \lf(\frac{d\eta}{d\tht}\rt)^2 - \frac{1}{\tht} \frac{d\eta}{d\tht}
+ \eta^3-\eta^{-1}
\ee
as
\be\label{eq73}
G_\pm (r) = \frac{1 \mp \eta\lf(\frac{r}{2}\rt)}{2\eta \lf(\frac{r}{2}\rt)^{1/2}}
\exp \lf[ \frac{1}{4} \int^\infty_{r/2} \tht \frac{\lf(1-\eta^2\rt)^2 -\lf(\eta'\rt)^2}{\eta^2}\;d\tht
\rt]
\ee
with the boundary condition $\eta(\tht) \sim 1-\frac{2}{\pi}\;K_0 (2\tht)$ as $\tht\to\infty$, where $K_0$
denotes the modified Bessel function.  
The above calculation was the first derivation of an explicit scaling law for a two-spin correlation function for
any model in statistical physics.  As noted before, the scaling functions $G_\pm (r)$ are believed to be
universal for a wide class of two-dimensional models with short range interactions.

A crucial issue for \cite{WMTB} scaling theory
was to show that formula (\ref{eq70}) matches the critical $1/4$ 
behavior (\ref{eq60}) as $r \to 0$. In order to do this, the authors needed to know the 
$ \infty \leftrightarrow 0$ {\it connection formulae} for the Painlev\'e 
functions $\eta(\theta)$.  In \cite{WMTB} they extracted this information from the unpublished thesis of Myers \cite{My} (see below). In a later paper in 1977, McCoy, Tracy, and Wu \cite{MTW} derived  $ \infty \leftrightarrow 0$ connection formulae for a two-parameter class of bounded solutions $\eta(\theta;\nu,\lambda)$ of the general Painlev\'e III equation
\[
\frac{d^2\eta}{d\tht^2} = \frac{1}{\eta} \lf(\frac{d\eta}{d\tht}\rt)^2 - \frac{1}{\tht} \frac{d\eta}{d\tht}
+ \frac{2\nu}{\tht}\lf(\eta^2-1\rt)+\eta^3-\eta^{-1}.
\]
We now describe in more detail what was done in \cite{MTW}, restricting our discussion to the case $\nu=0$
which is relevant to \cite{WMTB}.

Consider the one parameter family $\eta(\tht)=\eta(\theta;0,\lambda)$
of the solutions of equation (\ref{eq72}) defined by the asymptotic
condition
\begin{equation}\label{mccoy1}
\eta(\theta) \sim  1 - 2\lambda K_0(2\theta), \quad \theta \to \infty.
\end{equation}
In \cite{MTW}, it was shown that, for $0\le\lambda<1/\pi$, as $r \to 0$,
\begin{equation}\label{mccoy2}
\eta(r/2) = Br^{\sigma}\left(1 - \frac{1}{16}B^{-2}(1-\sigma)^{-2}r^{2-2\sigma} + O(r^2)\right)
\end{equation}
where the constants $B$ and $\sigma$ are functions of $\lambda$
given by the explicit formulae:
\begin{equation}\label{mccoy3}
\sigma=\sigma(\lambda) = \frac{2}{\pi}\arcsin(\pi\lambda), \quad 
B = B(\sigma) = 2^{-3\sigma}\frac{\Gamma((1-\sigma)/2)}{\Gamma((1+\sigma)/2)}.
\end{equation}
In \cite{My}, Myers derived the asymptotics for $\eta(r/2;0,1/\pi)$ as $r\downarrow 0$,
\be\label{mccoy4}
\eta\lf(\frac{r}{2};0,\frac{1}{\pi}\rt)\sim -\frac{r}{2}\lf(\log\frac{r}{8}+\gamma_E\rt)
\ee
where $\gamma_E$ is Euler's constant.
In \cite{MTW} the authors recovered \eqref{mccoy4} by analyzing formally the limit $\lambda\to 1/\pi$
in \eqref{mccoy2}. Proofs of \eqref{mccoy4} were given later in \cite{Wid8} and, using the Riemann-Hilbert method, in \cite{Niles}.  
In turn, the asymptotic formula (\ref{mccoy4})
yields the estimate:
\begin{equation}\label{mccoy5}
G_{\pm} \sim\mbox{const}\,r^{-1/4}\left[1 \pm \frac{r}{2}\left(\ln\frac{r}{8} + \gamma_E\right)
\right].
\end{equation}
In order to see that this estimate matches equation (\ref{eq60}), it is enough to notice that  (\ref{eq70}) and (\ref{eq71})
imply that
\begin{equation}\label{mccoy6}
\lf\langle\sigma_{1,1}\,\sigma_{1+n,1+n}\rt\rangle \sim 2^{1/4}G_{\pm}(r)|t|^{1/4} = 
2^{1/4}G_{\pm}(r)r^{1/4}n^{-1/4}.
\end{equation}
The main term in equation (\ref{eq60}), up to the multiplicative constant,
follows now directly from (\ref{mccoy5}). The evaluation of the constant in (\ref{mccoy5})
that matches the transcendental constant $A$ in (\ref{eq60}) was done by Tracy \cite{T} in 1991.

The Painlev\'e functions, PI, \dots, PVI,  (see \cite{In}) were discovered and analyzed intensively at the
beginning of the 20$^{th}$ century  until the First World War, and then fell into a period of latency.  They
reappeared in theoretical physics quite unexpectedly in 1965 in the work of Myers \cite{My} mentioned above, who showed that the
scattering of electromagnetic radiation from a strip in the plane could be expressed in terms of a solution
of the Painlev\'e III equation. The paper \cite{WMTB} was a seminal event in mathematical physics, followed soon after
by the landmark paper of Jimbo et al \cite{JMMS} on the impenetrable Bose gas, and the Painlev\'e equations
have emerged over the years as the core of modern special function theory, with applications across the
board in mathematics, physics and engineering (see, e.g., \cite{FIKN}). In  PDE's, in
particular, the Painlev\'e equations (specifically, Painlev\'e II) appeared for the first time in the
analysis of Ablowitz and Segur \cite{AS} of self-similar solutions of the Korteweg-de Vries equation.
This pioneering development opened the door to understanding the crucial
role that the Painlev\'e equations play in the theory of integrable dynamical systems (see \cite{AC}, \cite{FIKN}).
The work \cite{MTW} demonstrated for the first time the possibility to
derive explicit connection formulae for families of solutions to a Painlev\'e equation,
a property which had been previously known only for the linear ODEs of hypergeometric type.
Similar connection formulae were obtained for the second Painlev\'e equation in \cite{AS} (only
formulae for the amplitude) and in \cite{SA} (complete formulae for the amplitude and the phase). 
The papers  \cite{MTW},
\cite{AS}, and \cite{SA} generated a surge of activity in the global asymptotic
analysis of the Painlev\'e equations (see again \cite{FIKN} for more history and the state of the art in the area). We view these developments as one more example of the way in which questions arising in the analysis of the Ising model gave rise to developments in Toeplitz theory, and in this case, through Toeplitz theory, to the global theory of the Painlev\'e equations. 

An alternative to the representation (\ref{eq73}) of 
$G_{\pm}(r)$ was obtained in 1980 by Jimbo and Miwa \cite{JM}.
Their formula for $G_-(r)$ reads:
\begin{equation}\label{p5G}
G_-(r) = \exp\left[-\int_{2r}^{\infty}\frac{\sigma(x)}{x}dx\right],
\end{equation}
where $\sigma(x)$ is the unique solution of the Painlev\'e V equation
\begin{equation}\label{p5}
\left(x\frac{d^2\sigma}{dx^2}\right)^2 = \left(\sigma
-x\frac{d\sigma}{dx} +2\left(\frac{d\sigma}{dx}\right)^2 
\right)^2 
-4\left(\frac{d\sigma}{dx}\right)^2\left(\left(\frac{d\sigma}{dx}\right)^2 -\frac{1}{4}\right)
\end{equation} 
satisfying the boundary conditions
\be\label{westintro0}
\sigma(x)=\begin{cases}
-1/4+O(x\log x),& x\to
0,\cr 
\frac{1}{2\pi}x^{-1}e^{-x}
 \left(1 +
O\left(\frac{1}{x}\right)\right),& x\to +\infty.
\end{cases}
\ee
The equivalence of (\ref{p5G}) to the original PIII-formula (\ref{eq73}) was established
in \cite{MP}. 

Recently in \cite{CIK} the following generalization of the Jimbo-Miwa result was obtained.
The authors in \cite{CIK} considered the double-scaling limit for the Toeplitz
determinant $D_n(t)$ with the symbol:
\begin{equation}\label{ciksymb}
\varphi(z) = (z-e^{t})^{\alpha + \beta} (z-e^{-t})^{\alpha - \beta}z^{-\alpha +\beta}e^{-i\pi(\alpha +\beta)}e^{V(z)},
\end{equation}
where $t\geq 0$,  $\alpha \pm \beta \neq -1, -2, ....$, $\Re\alpha > -1/2$, and  $V(z)$ is
analytic in an annulus containing the unit circle. The powers in (\ref{ciksymb}) are 
defined by requiring the arguments to lie between zero and $2\pi$. The symbol $\varphi(z)$ 
is analytic in ${\mathbb C}\setminus \left([0, e^{-t}]\cup[e^{t}, +\infty)\right)$. 
It is regular  for $t>0$ and has a Fisher-Hartwig singularity for
$t=0$ at $z_0=1$, $\varphi(z) = |z-1|^{2\alpha}z^{\beta}e^{-i\pi \beta}e^{V(z)}$.
The main result in \cite{CIK} is the following:  There exists a finite (perhaps empty) set
$\De \subset (0,\infty)$ and a (small) number $t_0>0$, such that for $n\to\infty$ and for all
$0< t < t_0$, $\log D_n(t)$ has an  expansion of the form (see \cite{CIK} Theorem 1.4)
\begin{multline}\label{expansionDn}
\log D_n(t)=n V_0+(\alpha+\beta)nt+
\sum_{k=1}^{\infty}k\left[V_k-(\alpha+\beta)\frac{e^{-tk}}{k}\right]
\left[V_{-k}-(\alpha-\beta)\frac{e^{-tk}}{k}\right]\\
+\int_0^{2nt}\frac{\sigma(x)-\alpha^2+\beta^2}{x}dx+(\alpha^2-\beta^2)\log
2nt
+\log\frac{G(1+\alpha+\beta)G(1+\alpha-\beta)}{G(1+2\alpha)}+o(1),
\end{multline}
with the error term being  uniform for $0\leq t < t_0$ provided $\text{dist}\lf(2nt, \De\rt) \ge \de >0$ for some $\de >0$, and the path of integration in \eqref{expansionDn}
does not intersect with $\De$. (In fact, the asymptotics \eqref{expansionDn} hold uniformly in a sector
of the complex plane
$-\pi/2+\ep<\arg x<\pi/2-\ep$, $0<\ep<\pi/2$, away from a finite number of points. These points, which include
the points of $\De$, are possible
poles of the function $\sigma(x)$. A choice of the integration contour corresponds to a branch of the logarithm.)

The function $\sigma(x)$ in (\ref{expansionDn}) is the unique solution of the Painlev\'e V equation
\begin{multline}\label{p5ab}
\left(x\frac{d^2\sigma}{dx^2}\right)^2 = \left(\sigma
-x\frac{d\sigma}{dx} +2\left(\frac{d\sigma}{dx}\right)^2 + 2\alpha\frac{d\sigma}{dx}
\right)^2 \\
-4\left(\frac{d\sigma}{dx}\right)^2\left(\frac{d\sigma}{dx} +\alpha + \beta\right)\left(\frac{d\sigma}{dx} +\alpha - \beta\right)
\end{multline}
which is real analytic in $(0, +\infty)$  and satisfies  the boundary conditions
 \be\label{westintro1}
\sigma(x)=\begin{cases}\alpha^2-\beta^2+
\frac{\alpha^2-\beta^2}{2\alpha}\{x-x^{1+2\al}C(\alpha,
\beta)\}(1+O(x)),& x\to 0,\quad 2\alpha\notin {\mathbb Z}\cr
\alpha^2-\beta^2+O(x)+O(x^{1+2\alpha})+O(x^{1+2\alpha}\log x),& x\to
0,\quad 2\alpha\in {\mathbb Z}\cr x^{-1+2\alpha}e^{-x}
\frac{-1}{\Gamma(\alpha-\beta)\Gamma(\alpha+\beta)} \left(1 +
O\left(\frac{1}{x}\right)\right),& x\to +\infty,
\end{cases}
\ee
with
\be\label{Cab}
C(\alpha,\beta)=
\frac{\Gamma(1+\alpha+\beta)\Gamma(1+\alpha-\beta)}{\Gamma(1-\alpha+\beta)\Gamma(1-\alpha-\beta)}
\frac{\Gamma(1-2\alpha)}{\Gamma(1+2\alpha)^2}\frac{1}{1+2\alpha}.
\ee
For fixed $t >0$, equation (\ref{expansionDn}) yields the standard Szeg\H{o} large $n$  asymptotics
for the regular symbol (\ref{ciksymb}). 
It is interesting to note that this implies the following integral identity for the Painlev\'e
function $\sigma(x)$ (cf. (1.31) in \cite{CIK}): for any $x_0>0$,
\begin{multline}\label{identitycik}
\int_0^{x_0}\frac{\sigma(x)-\alpha^2+\beta^2}{x}dx 
+\int_{x_0}^{\infty}\frac{\sigma(x)}{x}dx +(\alpha^2 - \beta^2)\log x_0\\
= -\log\frac{G(1+\alpha+\beta)G(1+\alpha-\beta)}{G(1+2\alpha)}.
\end{multline}
In the case $t=0$, the equation (\ref{expansionDn}) transforms to a single-point  Fisher-Hartwig
formulae: To see this, one utilizes the identity $\sum_{k=1}^{\infty}e^{-2kt}/k = -\ln\left(1 - e^{-2t}\right)$.   

In the special case  $\alpha =0 $, $\beta = -\frac{1}{2}$,  $V(z) \equiv 0$, and with the identification,
$e^{-t}=k_{\tons}$, the symbol (\ref{ciksymb}) becomes $e^{-t/2}\varphi_\tdg(z)$. Therefore, in this case,
 $D_n(t)= e^{-nt/2}\lf\langle \sg_{1,1}\, \sg_{1+n, 1+n}\rt\rangle$ with  $t\downarrow 0$ corresponding to $T\uparrow
T_c$. Set $2nt \equiv 2r$. Then  formula (\ref{expansionDn}), together with the identity (\ref{identitycik})  
and  the definition (\ref{eq70}) of the scaling
functions $G_{\pm}(r)$, yields  the Jimbo-Miwa relation (\ref{p5G}). With the equivalence  
of (\ref{p5G}) and \eqref{eq73}, this
means that the {\it uniform} asymptotics  
(\ref{expansionDn}) include the double-scaling limit in \cite{WMTB} as a special case 
corresponding to $2nt \equiv 2r$ fixed 
(away from $\De$), $n\to\infty$, and
$\alpha =0 $, $\beta = -\frac{1}{2}$, $V(z) \equiv 0$.

\begin{remark}\label{remarkJM} 
In  \cite{JM}, the authors used a remarkable relation (which they themselves discovered) 
between the diagonal Ising 
correlations and isomonodromy deformations of a certain $2\times2$ Fuchsian system.
This discovery laid the foundation for the Riemann-Hilbert method in the theory of correlation
functions, random matrices, and Toeplitz and Hankel determinants (see \cite{DIZ} for more on  the history of these developments). Using this link to Fuchsian systems, it was shown in \cite{JM}
that the diagonal correlation function $\lf\langle \sg_{1,1}\, \sg_{1+n, 1+n}\rt\rangle$,
as a function of the variable 
$$
s = e^{-2t} \equiv \left(\sinh\frac{2J_1}{k_BT}\sinh\frac{2J_2}{k_BT}\right)^{-2}
$$
is expressed in terms of a solution to the Painlev\'e VI equation. The derivation
of (\ref{p5G}) in \cite{JM} is done by formally performing  a scaling transformation
taking the inverse monodromy problem associated with this Painlev\'e VI   to the inverse monodromy
(Riemann-Hilbert) problem associated with the fifth Painlev\'e function $\sigma(x)$. 
In the rigorous setting of \cite{CIK}, the latter problem appears as a parametrix 
during the implementation of the nonlinear steepest-descent method. 
\end{remark}
\begin{remark}\label{remarkJM3}
The  Painlev\'e VI description of the diagonal correlations in \cite{JM} was the first example 
of the appearance of the Painlev\'e functions in the theory of  Toeplitz
determinants before the large $n$ limit is taken. One realizes now that
nonlinear differential equations of Painlev\'e type are always present in 
situations where the  symbol $\varphi(z)$ satisfies the condition:
\begin{equation}\label{quasiclassymb}
\frac{d\log\varphi(z)}{dz} = \mbox{rational function}.
\end{equation}
This fact is quite easy to see within the Riemann-Hilbert formalism (see, e.g., discussion
in \cite{ITW}). These ``finite $n$'' Painlev\'e representations are not universal:
the type of the nonlinear ODE depends strongly on the structure of the rational right
hand side of (\ref{quasiclassymb}), i.e.\ on the symbol. 
\end{remark}

Double-scaling problems for Toeplitz determinants occur in many different areas.  For example, in
combinatorics, let $\pi = \lf(\pi(1), \pi(2), \dots, \pi(N)\rt)$ be a permutation of the numbers
$1, 2, \dots, N$.  If $1\le i_1 < i_2 < \dots < i_k \le N$ and $\pi(i_1) < \pi(i_2) < \dots < \pi(i_k)$,
we say that $\pi(i_1), \dots, \pi(i_k)$ is an \tit{increasing subsequence} in $\pi$ of length $k$.  Let
$\ell_N(\pi)$ denote the maximum length of all increasing subsequences in $\pi$.  With uniform distribution on
the permutations, one is
interested (see \cite{BDJ}) in the distribution function
$$
p_N(n)= \text{Prob } \lf\{ \ell_N(\pi) \le n\rt\} = \frac{\#\{ \pi:\ell_N(\pi)\le n\}}{N!}
$$
as $N, \, n\to\infty$.  It turns out \cite{Ges} that the Poissonized version $\phi_n(\lb)$ of the
process $p_N(n)$
\be\label{eq129}
\phi_n(\lb) = \sum^\infty_{N=0} e^{-\lb} \frac{\lb^N}{N!} \; p_N(n), \qquad \lb>0
\ee
can be expressed in terms of a Toeplitz determinant
\be\label{eq130}
\phi_n(\lb) = e^{-\lb} \, D_n(\varphi_\lb)
\ee
where $\phi_\lb(\tht)=e^{2\sqrt{\lb} \cos \tht}$.  In order to recover $p_N(n)$ as $N,\, n\to\infty$
from \eqref{eq129} \eqref{eq130} one must analyze the double-scaling limit for the determinant as $\lb$,
$n\to\infty$.  The main result in \cite{BDJ} is the following:
\be\label{eq131}
\lim_{N\to\infty} \text{Prob } \lf(\frac{\ell_N-2\sqrt{N}}{N^{1/6}} \le t\rt) = F_2(t)
\ee
where $F_2(t)$ is the Tracy-Widom distribution \cite{TW} for the largest eigenvalue of a random matrix
in the Gaussian Unitary Ensemble.  The distribution $F_2(t)$ can be expressed in terms of Painlev\'e
functions as follows:
\be
F_2(t) =e^{-\int^\infty_t (s-t)\, u^2(s)\,ds}
\ee
where $u(s)$ is the unique solution of $PII$,
\be
u'' (s) =2 u^3(s) +s u(s)
\ee
with
$$
u(s) \sim -\mathrm{Ai}\,(s)\quad \text{ as }\quad s\to +\infty.
$$
Here $\mathrm{Ai}\,(s)$ is the standard Airy function. On the other hand, $F_2(t)$ can be expressed as a Fredholm determinant
\be\label{TWairydet}
F_2(t)=\det(I-K^{(s)}_{\mathrm{Airy}})
\ee
where $K^{(s)}_{\mathrm{Airy}}$ is the trace-class operator with kernel
\[
K^{(s)}_{\mathrm{Airy}}(x,y)=\frac{\mathrm{Ai}\,(x){d\over dy}\mathrm{Ai}\,(y)-
\mathrm{Ai}\,(y){d\over dx}\mathrm{Ai}\,(x)}{x-y}
\]
acting on $L^2(s,\infty)$. In random matrix theory,
this determinant is the probability that the interval $(s,\infty)$ contains no eigenvalues in the edge scaling limit
in the Gaussian Unitary Ensemble.

The result in \cite{BDJ} was the first of many results linking problems in
pure and applied mathematics and in physics to the Tracy-Widom distribution $F_2(t)$ and its symplectic
and orthogonal analogs (see, for example \cite{Dei2}).

As already noted on several occasions,
in his work on impenetrable bosons \cite{Len1}, Lenard was led to consider Toeplitz determinants with
symbols of type $\lf|\eit - e^{i\tht_1}\rt| \lf|\eit-e^{i\tht_2}\rt|$ (see \eqref{eq81} above) that vanished
on $S^1$.  In \cite{Len2} he then considered symbols with any finite number of such $\al$-type FH
singularities.  Such symbols vanish at only a finite number of points, and in \cite{Wid6} Widom began
considering symbols which vanished on a full interval.  Let $\varphi_\mu\!\lf(\eit\rt)$ denote the
characteristic function of the interval $\lf(\mu, 2\pi-\mu\rt)$, $0<\mu <\pi$. In a virtuoso calculation
Widom showed that, for $\mu$ fixed, as $n\to\infty$
\be\label{eq132}
\log D_n(\varphi_\mu) = n^2 \log\cos \frac{\mu}{2} - \frac{1}{4} \log \lf(n\sin \frac{\mu}{2}\rt)
+ c_0 + o(1)
\ee
where $c_0 = \frac{1}{12} \log 2 + 3 \zeta'(-1)$ and $\zeta(z)$ is the Riemann zeta-function.
A few years later in 1976, Dyson \cite{Dy4} returned to the problem 
(cf. discussion following \eqref{6.1} above) of computing
the probability $P_s$ that there are no eigenvalues for a random
matrix in the interval $(0,2s/\pi)$, in the bulk scaling limit for the Gaussian Unitary Ensemble.  He
showed that as $s\to\infty$, $\log P_s$ has a full asymptotic expansion
\be\label{eq133}
\log P_s =  -\frac{s^2}{2} - \frac{1}{4} \log s + a_0 + \frac{a_1}{s} + \frac{a_2}{s^2}
+ \dots \ ,
\ee
Dyson identified all the constants $a_0, a_1, a_2, \dots ,$ and it turns out that $a_0$, which is of
particular interest, is precisely the constant $c_0$ arising in Widom's expansion \eqref{eq132} above.
Dyson arrived at the identification $a_0=c_0$ by noting first that the probability $P_s$ is given by (see, e.g., 
\cite{Meh}, \cite{Dei3})
\be\label{12.1}
P_s=\det(I-K_s)
\ee
where $K_s$ is the trace-class operator with kernel
\be\label{sinekernel}
K_s(x,y)=\frac{\sin(x-y)}{\pi(x-y)}
\ee
acting on $L^2(-s,s)$. A simple calculation shows that 
\be
D_n(\varphi_\mu)=\det\lf(\delta_{jk}-\frac{\sin\mu(j-k)}{\pi(j-k)}\rt)_{0\le j,k\le n-1}
\ee
and so for fixed $s>0$,
\be\label{eq134}
\lim_{n\to\infty} D_n \!\lf(\varphi_{\frac{2s}{n}}\rt) =P_s
\ee
and hence, if the error term $o(1)$ in \eqref{eq132} was uniform in the double scaling limit, $n\to\infty$,
$\mu\to 0$ such that $\mu n=2s$, one could conclude from \eqref{eq132} \eqref{eq134} that $a_0$ is
indeed equal to $c_0$.  The uniformity of the error term $o(1)$, however, remained open.  There are now
three proofs that $a_0=c_0$.  The first two were given simultaneously and independently
by Ehrhardt and Krasovsky, and the third was given a little later in \cite{DIKZ}.   The proofs by
Krasovsky and by Deift, Its, Krasovsky, and Zhou proceed by following Dyson's observation \eqref{eq134} and verifying that the $o(1)$ error term is indeed uniform.  We refer the reader to
\cite{DIKZ} for some history of the problem $P_s$,
$s\to\infty$, and, in particular, for the relevant references to the work of Krasovsky and Ehrhardt.
A similar asymptotic problem arises when instead of the sine kernel determinant \eqref{12.1} one considers
the general so-called confluent hypergeometric kernel determinant which depends on two extra parameters
(and is related to the so-called  Bessel kernel determinant and also to the determinant \eqref{12.1}  for particular 
choices of the parameters); another related problem is the asymptotics for the Tracy-Widom distribution \eqref{TWairydet} as $s\to-\infty$. 
These problems can also be represented as double scaling problems for Toeplitz
(and also Hankel) determinants, and such representations were crucial for finding complete solutions
for the asymptotic problems at hand.
We refer the reader to \cite{Kr} 
for a discussion of the results and for the references.

Dyson's identification of $a_1,a_2,\dots$ in \cite{Dy4} involves an ingenious application of inverse scattering theory for Schr\"odinger operators on the half-line $0\le x<\infty$. He begins by noting that $P_s=\det(I-K_s)$ can be factored
\be\label{13.1}
\det(I-K_s)=D_+(s)D_-(s)
\ee
where $D_\pm(s)=\det(I-f_\pm)_{L^2(0,s)}$ are Toeplitz+Hankel type determinants on $L^2(0,s)$ with kernels
\be\label{13.2}
f_\pm(x,y)={1\over\pi}\lf[\frac{\sin(x-y)}{x-y}\pm \frac{\sin(x+y)}{x+y}\rt].
\ee
He then considers the functions 
\be\label{13.3}
W_\pm(s)=-2{d^2\over ds^2}\log D_\pm(s)-1
\ee
as potentials for Schr\"odinger operators $H_\pm=-{d^2\over ds^2}+W_\pm(s)$ on $L^2(0,\infty)$ with appropriate Robin boundary conditions at $s=0$. Using the Gelfand-Levitan formalism (see, e.g., \cite{Fad}), Dyson is able to identify the spectral measures $\rho_\pm(\lb)d\lb$ for the operators $H_\pm$. 
He then obtains the scattering phase functions $e^{i\eta_\pm}$.
In the Marchenko formalism (see, e.g., \cite{Fad}) one utilizes the functions $e^{i\eta_\pm}$ to express the potentials $W_\pm(s)$ in terms of Fredholm determinants of the form 
\be\label{14.1}
W_\pm(s)=-{1\over 4}\lf(s\pm\half\rt)^{-2}\!\!-2{d^2\over ds^2}\log \Delta_\pm(s),\qquad s>\half
\ee
where $\Delta_\pm(s)=\det(I-F_\pm)_{L^2(s,\infty)}$ and $F_\pm(x,y)$ are certain explicit kernels which
decay as $s\to\infty$. Equating \eqref{13.3} and \eqref{14.1}, and integrating twice, one obtains the following identity for $s>\half$:
\be\label{14.2}
\log P_s=\log (I-K_s)=\log D_+D_-=-\half s^2-{1\over 8}\log\lf( s^2-{1\over 4}\rt)+a_0+\log\lf(\Delta_+(s)\Delta_-(s)\rt)
\ee
(the terms linear in $s$ drop out). As $F_\pm(x,y)$ decay as $x,y\to\infty$, it is clear that \eqref{14.2} gives rise to a full expansion of the form \eqref{eq133}, where the coefficients can be computed  in terms of the trace powers
$\tr F_\pm^m$.

Dyson is able to compute $e^{i\eta_\pm}$, and hence $F_\pm$, explicitly, because of the following fundamental fact from inverse scattering theory. The spectral weights have the form $\rho_\pm(\lb)=\sqrt{\lb}/(\pi |\varphi_\pm(\sqrt{\lb})|^2)$,
where $\varphi_\pm(k)$, $k=\sqrt{\lb}$, are certain functions constructed from the Jost solutions $J_\pm(s,k)$ for 
$H_\pm$, $H_\pm J_\pm(s,k)=k^2J_\pm(s,k)$, $J_\pm(s,k)\sim e^{iks}$ as $s\to\infty$. The functions $\varphi_\pm(k)$
are analytic in $\{\Im k>0\}$ with prescribed asymptotics as $k\to\infty$, and in the cases at hand, have no zeros in
$\{\Im k>0\}$. Hence, by standard arguments, $\varphi_\pm(k)$ are determined by their absolute value $|\varphi(k)|$ for
$\Im k=0$. Thus $\rho_\pm(\lb)$ determine $\varphi_\pm(k)$. However, $e^{i\eta_\pm}$ can be expressed in terms of 
$\varphi_\pm$ as $e^{i\eta_\pm(k)}=i\overline{\varphi_\pm(k)}/|\varphi_\pm(k)|$, and hence $\rho_\pm(\lb)$ determine the scattering phase functions, and hence, eventually, $F_\pm$.

These considerations bring to mind the Borodin-Okounkov-Geronimo-Case 
formula \eqref{eq38} in which the Szeg\H o function
 ${\cal D}(z) = \exp \lf(\frac{1}{4\pi} \int^\pi_{-\pi}  \frac{e^{i\tht}+z}{e^{i\tht}-z}
\;\log \varphi\!\lf(\eit\rt)d\tht\rt)$, which is analytic in $|z|<1$, plays the role of the functions $\varphi_\pm(k)$. Indeed, 
$|{\cal D}(\eit)|^2=\varphi(\eit)$, which is the spectral measure for the associated Cantero-Moral-Velasquez (CMV)
unitary operator acting on $\ell_+^2$ (see \cite{Sim2}). Thus we may think of the Toeplitz determinant $D_n(\varphi)$
as an expression in the ``Gelfand-Levitan formalism''. On the other hand, we may think of the Fredholm determinant on the RHS of \eqref{eq38} as an expression in the ``Marchenko formalism''. Indeed, from \eqref{eq39}, 
$b=\overline{\cal D}^2/|{\cal D}|^2$, which we may view as a scattering phase function $b=e^{2i\eta}$ for the CMV operator, and $c$ is just $\overline{b}=e^{-2i\eta}$. These considerations are very suggestive and invite further investigation. One notes that the original proof of \eqref{eq38} by Geronimo and Case in \cite{GerCase} was developed in the context of scattering theory.

Another example of a double-scaling limit is the following.
Recall that the density matrix for the impenetrable bosons is given by 
$\rho_{N,L}(\xi)=\frac{1}{L}R_N(2\pi\xi/L)$, where $R_N$ is the Toeplitz determinant \eqref{34p3} with symbol
\eqref{Lensym}. For definiteness, fix the length scale so that $L=N$.
As we discussed in Remark \ref{Lenard} in Section \ref{sectionFH}, of particular interest is the 
double-scaling limit $N\to\infty$, $\xi/ N\to 0$.
In \cite{Len1} Lenard considered this limit 
$\rho(\xi)=\lim_{N\to\infty}{1\over N}R_N(2\pi\xi/N)$ with $\xi$ fixed.
He obtained an expression for $\rho(\xi)$ in terms of a Fredholm determinant with explicit kernel related
to the sine kernel \eqref{sinekernel}. He also obtained another representation for $\rho(\xi)$ in terms of
a Fredholm minor for the kernel $2\frac{\sin\pi |\xi|(x-y)}{\pi(x-y)}$. 

In order to complete his analysis of the limiting momentum distribution in the gas, Lenard needed 
to know the large $\xi$ behavior of $\rho(\xi)$, but he did not derive this behavior in \cite{Len1}.
The derivation of this limiting behavior was carried out by Vaidya and Tracy \cite{VaiTr} who found that
\be\label{22.1}
\rho(\xi)=\frac{\rho_\infty}{(\pi\xi)^{1/2}}
\lf(1+{1\over 8(\pi\xi)^2}\lf(\cos2\pi\xi -{1\over 4}\rt)+\cdots\rt),\qquad \rho_\infty=\pi e^{1/2} 2^{-1/3}A^{-6},
\ee
as $\xi\to\infty$,
where $A$ is again Glaisher's constant \eqref{Glai}. Note that 
if one evaluates  $\rho_{N,N}(\xi)=\frac{1}{N}R_N(2\pi\xi/N)$, 
$N, \xi=\frac{tN}{2\pi}\to\infty$ using \eqref{eq84} \eqref{eq85}, and assumes that the formulae remain valid as 
$t=2\pi\xi/N\to 0$, then one obtains the leading term in \eqref{22.1}.

The double-scaling limit $N\to\infty$, $\xi/ N\to 0$ corresponds to the Toeplitz determinant
with two FH singularities merging (at $z=1$)
with $\al_1=\al_2=\half$, $\beta_1=\beta_2=0$.
A general analysis of the double-scaling asymptotics of a Toeplitz determinant with two merging FH singularities (in particular, of Lenard's symbol \eqref{34p3}) is given in \cite{CK}.

Dyson also considered \cite{Dy5} the following scaling problem associated with Toeplitz determinants.
Let 
\be
\De(z,t)=\det(I-z\wh K_t)
\ee
where $0<z<1$ and $\wh K_t$ acts on $L^2(-t,t)$ with kernel $\wh K_t(x,y)=\frac{\sin\pi(x-y)}{\pi(x-y)}$. With the extra 
$\pi$ in the sine function (cf. $K_s$ in \eqref{sinekernel}), $\De(1,t)$ is now the probability that there are no eigenvalues in the interval $(0,2t)$. Of interest here is the double scaling limit for $\De(z,t)$ as $t\to\infty$ and $z\uparrow 1$. Dyson analyses the problem by interpreting $\De(z,t)$ in terms of a Coulomb gas at inverse temperature $\beta=2$. In this interpretation, $z$ corresponds to an external potential $v$ via the formula $1-z=e^{-\beta v}$. Then
\be
\De(z,t)=\frac{Z_2(v,2t)}{Z_2(0,2t)}
\ee
where
\be\label{23.2}
Z_\beta=\sum e^{-\beta (W+vn)}
\ee
is the partition function of the gas with external potential $v$ applied to charges in a fixed interval of length $2t$.   
The sum in \eqref{23.2} is to be interpreted as an infinite-dimensional integral over all configurations of the gas, $W$ represents the sum of all Coulomb interactions $W_{jk}=-\log|x_j-x_k|$, and for any configuration, $n$ is the number of charges in a fixed interval of length $2t$. Dyson does not define the sum rigorously, but he uses it as a heuristic guide in his calculations. The main result in the paper is the calculation of oscillatory factors in $\De(z,t)$ as $t\to\infty$: These factors turn out to be 
so-called ``genuinely non-linear oscillations'' and are
described by
Jacobi elliptic functions. For $0<z<1$, $0<\lb<1$, the symbol 
\be
f_{z,\lb}(\eit)=\chi_{[\pi\lb,\pi(2-\lb)]}+(1-z)\chi_{[-\pi\lb,\pi\lb]}
\ee
gives rise to the Toeplitz determinant
\be
D_n(f_{z,\lb})=\det\lf(\de_{jk}-z\frac{\sin\lb\pi(j-k)}{\pi(j-k)}\rt)_{0\le j,k\le n-1}
\ee
and so for fixed $z$, $t$
\be
\lim_{n\to\infty}D_n\lf(f_{z,{2t\over n}}\rt)=\De(z,t).
\ee
Thus Dyson's double scaling limit $z\uparrow 1$, $t\to\infty$ can be viewed as a triple scaling limit for the Toeplitz determinant $D_n(f_{z,2t/n})$ with $z\uparrow 1$, $t\to\infty$, and $n\to\infty$. Note finally that $f_{z,\lb}$ is precisely the kind of symbol that arises in the Toeplitz eigenvalue problem \eqref{eq120}. 

\bigskip

Lenard ends his 1972 paper \cite{Len2} with a hope and a prophecy:  ``It is the author's hope that
a rigorous analysis will someday carry the results to the point where the true role of the zeros of the
generating function will be understood.  When that day comes a capstone will have been put on  a
beautiful edifice to whose construction many contributed and whose foundations lie in the studies of
Gabor Szeg\H{o} half a century ago''.

Almost 100 years have passed since Szeg\H{o} published his first paper on the asymptotics of Toeplitz
determinants.  We certainly know more now than in 1972, but many new problems  continue to arise, and
it is still too early, alas, to set the capstone.

\section*{Acknowledgements:}
The work of the first two authors was supported in part by NSF grant DMS--1001886 and DMS--1001777 respectively.
The work of the third author was supported in part by EPSRC grant \#EP/E022928/1.  The authors would like 
to thank Estelle Basor, Rodney Baxter, Albrecht B\"ottcher, Freeman Dyson,
Michael Fisher, Peter Forrester, Helge Holden, Leo Kadanoff, Donald Knuth, Andrew Lenard,
Barry McCoy, John Palmer, Bob Shrock, Barry Simon, Tom Spencer, Craig Tracy, and Elias Wegert
for their help in preparing this manuscript.   Estelle Basor,  Albrecht B\"ottcher, 
Michael Fisher, Andrew Lenard, and Barry McCoy, in particular, graciously 
responded to many emails, and provided many references to the literature going back more than 60 years.

\newpage

\frenchspacing
\bibliographystyle{plain}

\end{document}